\documentclass{siamart0516}
\usepackage[margin=1in]{geometry}

\usepackage{tikz}
\usetikzlibrary{calc}
\usetikzlibrary{matrix}
\usepackage{etoolbox}
\usepackage{amsmath,amssymb,amsfonts} % Typical maths resource packages
\usepackage{hyperref}                 % For creating hyperlinks in cross references
\usepackage[normal]{subfigure}

\renewcommand{\eqref}[1]{(\ref{#1})}
\usepackage[numbers,sort&compress]{natbib}

\newcommand{\bigoh}{\mathcal{O}}

 \def \M{\mathbb{M}}
 \def \J{\mathbb{J}}

\newcommand{\eps}{\varepsilon}

\newcommand{\bsub}{\begin{subequations}}
\newcommand{\esub}{\end{subequations}$\!$}

\newcommand{\bx}{\mathbf{x}}
\newcommand{\by}{\mathbf{y}}
\newcommand{\bX}{\mathbf{X}}
\newcommand{\bU}{\mathbf{U}}
\newcommand{\bV}{\mathbf{V}}
\newcommand{\bxi}{\bs{\xi}}
\newcommand{\skel}{\mathcal{S}_{\Omega}}
\newcommand{\sat}{t_{\mathcal{S}}}

\newcommand{\ds}[0]{\displaystyle}
\newcommand{\bs}[0]{\boldsymbol}

\newcommand{\bey}{\begin{eqnarray}}
\newcommand{\eey}{\end{eqnarray}}

\newcommand{\beq}{\begin{equation}}
\newcommand{\eeq}{\end{equation}}
\theoremstyle{plain}% default

\theoremstyle{definition}

\theoremstyle{remark}
\newtheorem{exam}{\hspace{1mm}Example}[section]

\renewcommand{\abstractname}{Abstract}

\title{Moving Mesh simulation of contact sets in two dimensional models of elastic-electrostatic deflection problems}
\date{\today}

\author{
Kelsey L. DiPietro\thanks{Department~of Applied and Computational Mathematics \& Statistics, University of Notre Dame, Notre Dame, Indiana 46656, USA. {\tt Kelsey.L.Dipietro.5@nd.edu}} \and
Ronald D. Haynes\thanks{Department~of Mathematics and Statistics, Memorial University of Newfoundland, St.~John's, NL A1C 5S7, Canada  {\tt rhaynes@mun.ca}} \and
Weizhang Huang\thanks{Department~of Mathematics, University of Kansas, Lawrence, Kansas 66045, USA {\tt whuang@ku.edu}} \and
Alan E. Lindsay\thanks{Department~of Applied and Computational Mathematics \& Statistics, University of Notre Dame, Notre Dame, Indiana, 46656, USA. {\tt a.lindsay@nd.edu}} \and
Yufei Yu\thanks{Department~of Mathematics, University of Kansas, Lawrence, Kansas 66045, USA {\tt y920y782@ku.edu}}
}

\begin{document}

\baselineskip=16pt

\maketitle

\begin{abstract}
Numerical and analytical methods are developed for the investigation of contact sets in electrostatic-elastic deflections modeling micro-electro mechanical systems. The model for the membrane deflection is a fourth-order semi-linear partial differential equation and the contact events occur in this system as finite time singularities. Primary research interest is in the dependence of the contact set on model parameters and the geometry of the domain. An adaptive numerical strategy is developed based on a moving mesh partial differential equation to dynamically relocate a fixed number of mesh points to increase density where the solution has fine scale detail, particularly in the vicinity of forming singularities. To complement this computational tool, a singular perturbation analysis is used to develop a geometric theory for predicting the possible contact sets. The validity of these two approaches are demonstrated with a variety of test cases.
\end{abstract}

\section{Introduction}\label{sec:introduction}

\begin{figure}[htbp]
\centering
\subfigure[Schematic diagram.]{\includegraphics[width = 0.45\textwidth]{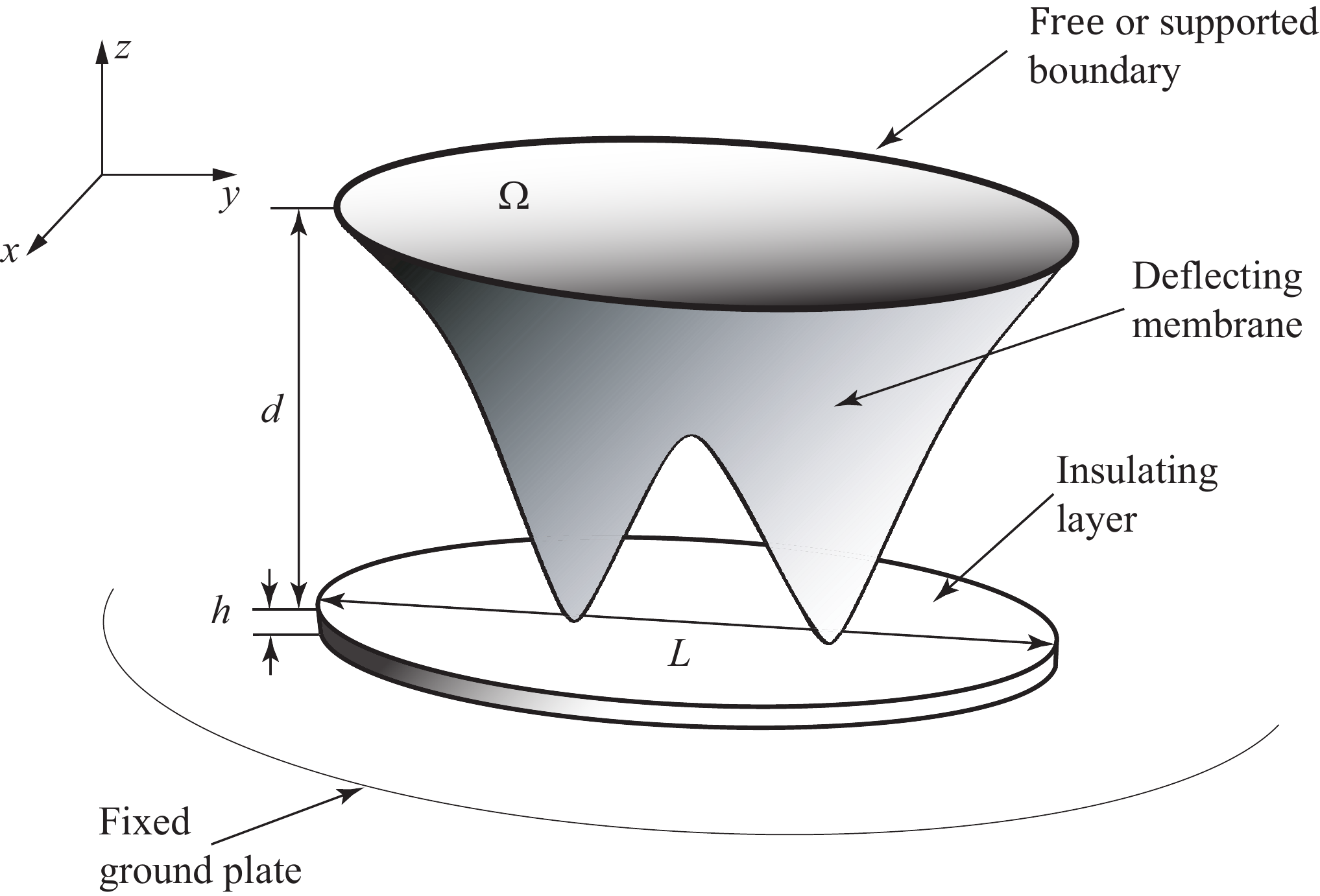} \label{fig:intro_diagram_A}}
\qquad
\subfigure[A MEMS device (source: \cite{sandia})]{\includegraphics[width = 0.45\textwidth]{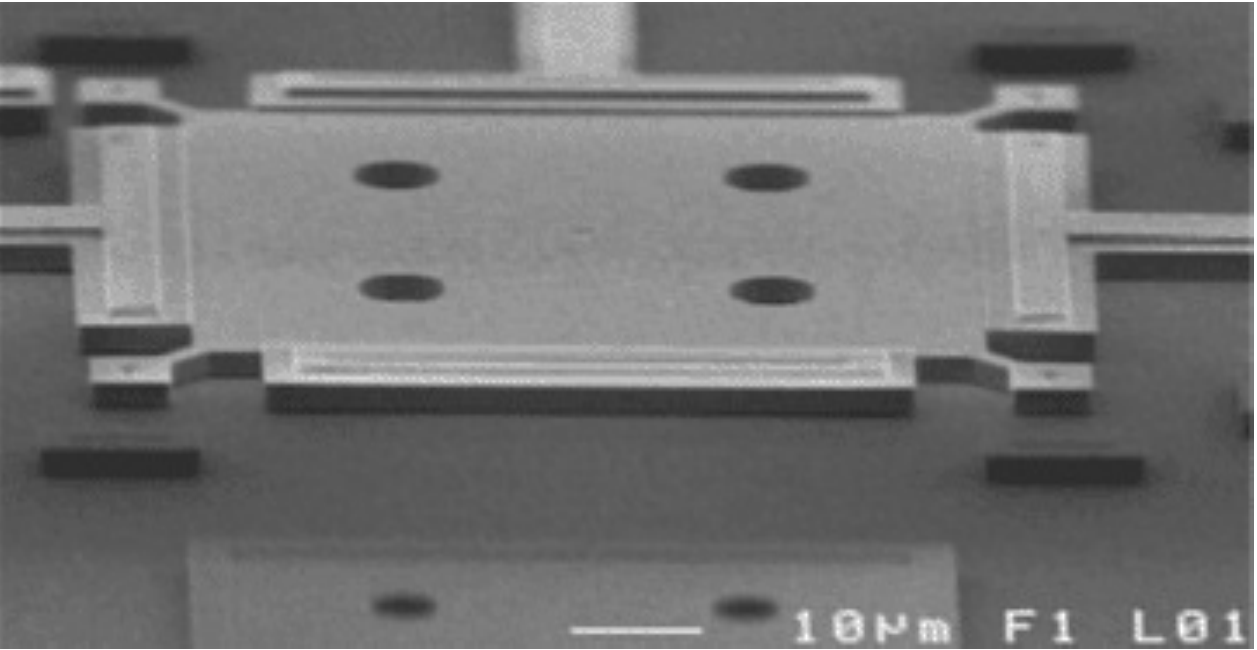} \label{fig:intro_diagram_B}}
\caption{A MEMS device (right) and a schematic (left) around which models are formulated. \label{fig:intro_diagram}}
\end{figure}

The present work is concerned with numerical simulation and singular perturbation analysis of the initial-boundary value problem of a fourth-order parabolic partial differential equation (PDE)
\begin{equation}\label{eq:intro}
\begin{cases}
u_t = - \eps^2 \Delta^2 u  - \ds\frac{1}{(1+u)^2}, & (\bx,t) \in\Omega \times(0,T);\\[5pt]
u = \Delta u = 0, & (\bx,t)\in\partial\Omega \times(0,T); \\[5pt]
u(\bx,0) = 0, & \bx \in \Omega;
\end{cases}
\end{equation}
in a variety of bounded two dimensional geometries $\Omega$. The system \eqref{eq:intro} models the non-dimensional vertical deflection $z = u(\bx)$ for $\bx = (x,y)\in\Omega$ of a Micro-electro mechanical systems (MEMS) capacitor \cite{Equilibrium2014,PB,Pelesko2002}. The MEMS capacitor is a key component of modern nanotechnology \cite{Tsai07,Beek2012,BatraReview} that features a deformable elastic membrane held fixed above a rigid substrate (see Fig.~\ref{fig:intro_diagram_A}). When an electric potential is applied between the deflecting plates, the top surface deforms towards the substrate. In equation \eqref{eq:intro}, the parameter $\eps$ quantifies the relative importance of electrostatic and elastic forces in the system. If the restorative elastic forces are too weak, the attractive Coulomb forces between the two surfaces will bring them into physical contact. This event, called {\em touchdown} or snap-through, can be useful or deleterious to operation, depending on the design of particular MEMS. The mechanism of the pull-in phenomenon has been studied extensively and many references can be found in the reviews \cite{BatraReview,ZHANG2014}.

The design and operation of MEMS can be aided by placing physical limiters or constraints at locations where contact between the two membranes is more likely \cite{Krylov2010}. These limiters can prevent damage to the device that could occur when the two surfaces meet. In addition, they allow for bistability in the system by creating stable large deflection configurations \cite{Equilibrium2014,Lindsay2016,Qiu2004,Lulinsky2011,Craighead2008}. Therefore, it is important to know at which location(s) in $\Omega$ singularities can form. In the one-dimensional case with $\Omega=(-1,1)$, equation \eqref{eq:intro} can form one singularity at the origin or two singularities located symmetrically about the origin, depending on the particular value of $\eps$ \cite{Capacitor2012}. In the physically relevant two-dimensional scenario, the details of the geometry $\Omega$ and the dependence on the parameter $\eps$ combine to make the set of possible singularity locations much more complex \cite{Lindsay2013,selfsimilar2014}.

Touchdown is a very rapid process in which energy is rapidly focussed in small spatial regions of $\Omega$. This process is manifested in the governing equations \eqref{eq:intro} by a finite time \emph{quenching} singularity. The term quenching refers to the fact that $u(\bx,t)$ is finite at the point of singularity  while $u_t(\bx,t)$ diverges as $t\to T$. Theoretical results on the quenching behavior of fourth-order parabolic equations such as \eqref{eq:intro} have established conditions under which quenching may occur \cite{Capacitor2012,Lin2007}, studied the local form of the profile near singularity \cite{Lindsay2013,Capacitor2012,BGW,GALAK} and given upper and lower estimates of the singularity time \cite{Friedman88,Capacitor2012,Philippin15}. For reviews on the extensive literature on blow-up/quenching for parabolic PDEs, see \cite{Galaktionov2002} and references therein.

\begin{figure}[htbp]
\centering
\subfigure[$\eps= 0.02$.]{\includegraphics[width = 0.3\textwidth]{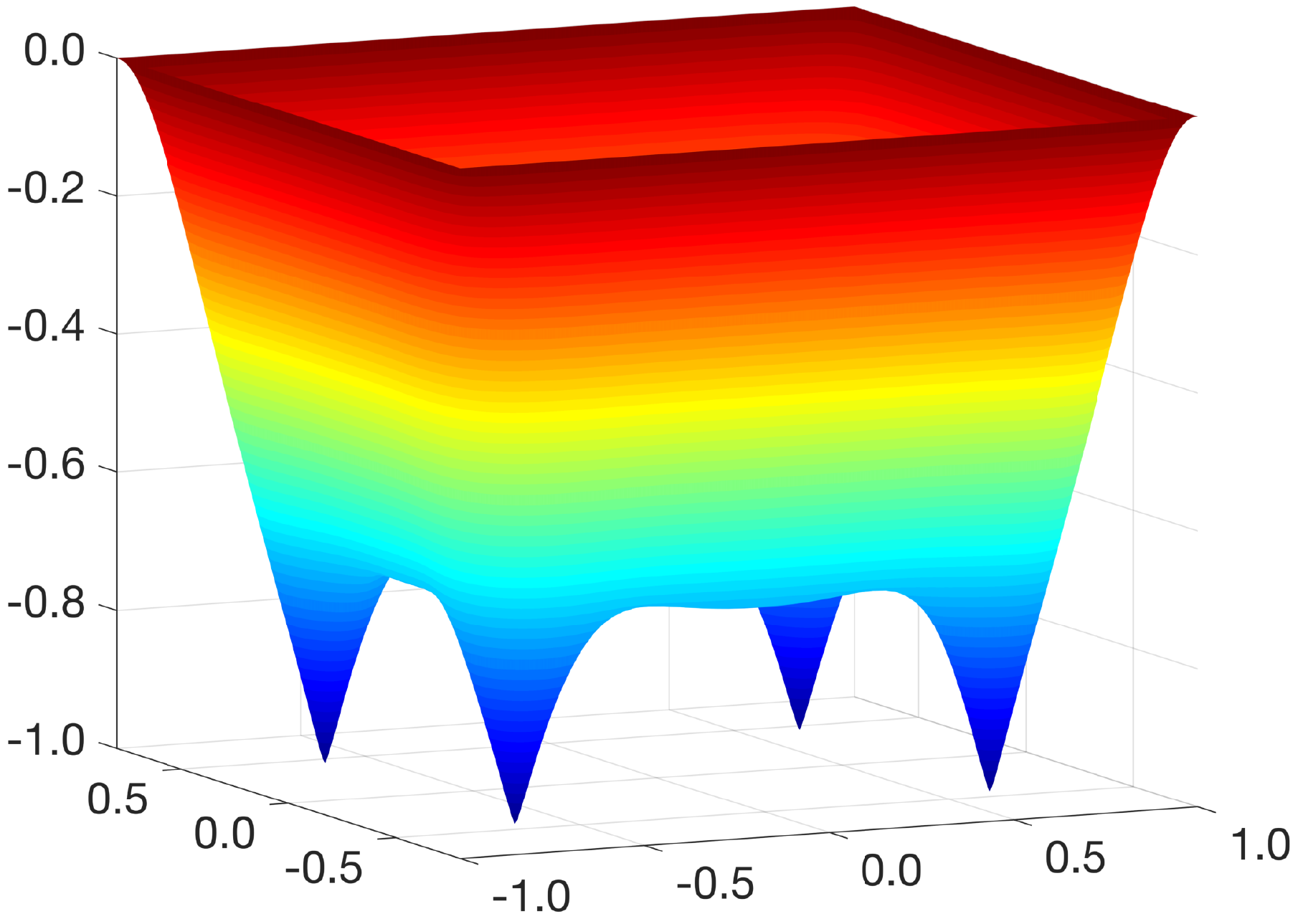} \label{fig:introsol_a}}\quad
\subfigure[$\eps= 0.068$.]{\includegraphics[width = 0.3\textwidth]{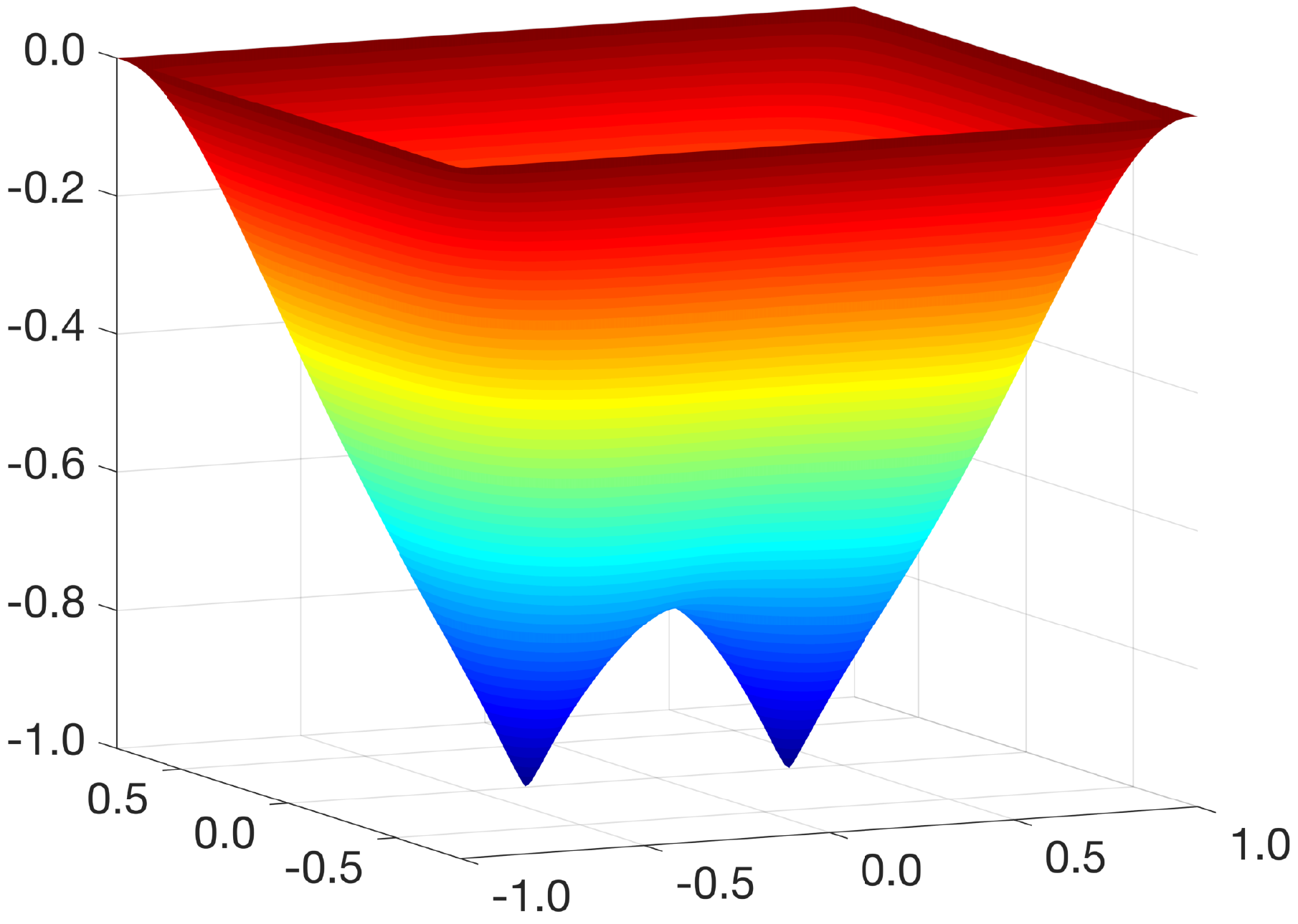} \label{fig:introsol_b}}\quad
\subfigure[$\eps= 0.1$.]{\includegraphics[width = 0.3\textwidth]{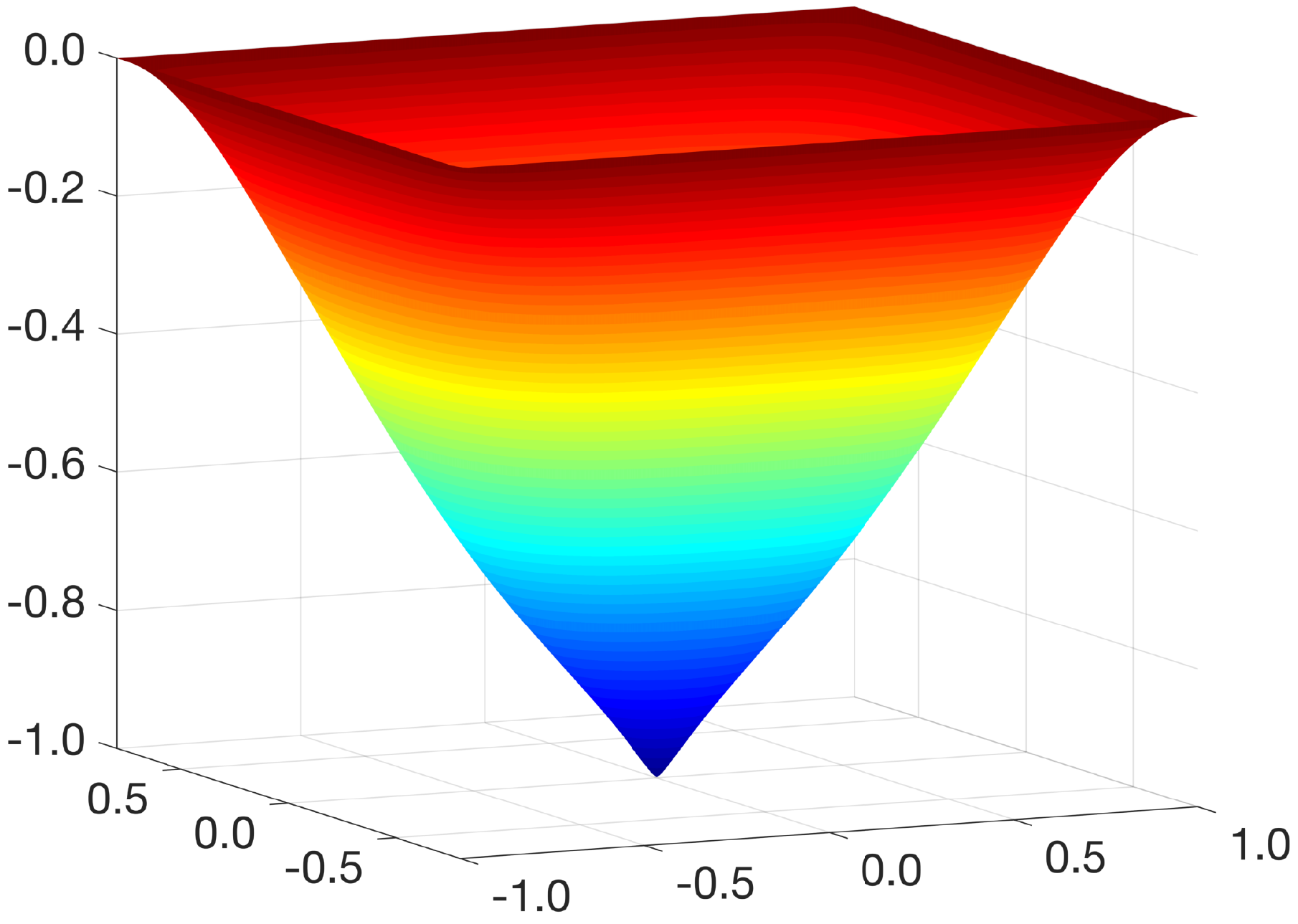} \label{fig:introsol_c}}
\caption{Solutions $u(\bx,t)$ of \eqref{eq:intro} at touchdown for $\eps = 0.02,\, 0.068,\, 0.1$ in the rectangle $\Omega = (-1,1)\times(-0.8,0.8)$. \label{fig:introsol}}
\end{figure}

The aim of this paper is to explore, through numerical simulations and asymptotic analysis of \eqref{eq:intro} as $\eps\rightarrow 0$, the potential set of locations at which touchdown may occur. In particular, we consider how the complex geometry and topology of MEMS devices, as seen in Fig.~\ref{fig:intro_diagram_B}, influences the possible set of contact locations. We present an adaptive moving mesh strategy \cite{HR11} for the solution of \eqref{eq:intro} which dynamically relocates the mesh points to provide resolution in the vicinity of forming singularities. An example of our method for the rectangular domain $\Omega = (-1,1)\times(-0.8,0.8)$ is shown in Fig.~\ref{fig:introsol} which shows either $4$, $2$ or $1$ contact points for different values of $\eps$. The sensitive dependence of the contact set on the domain $\Omega$ and parameter $\eps$ will be explored with a geometric skeleton theory (Sec.~\ref{sec:skeleton}) and adaptive numerical simulations (Sec.~\ref{SEC:fem}). 

In previous computational studies of singularity formation in second-order PDEs, moving mesh methods based on parabolic Monge-Amp\'{e}re (PMA) discretization have been successfully employed in one-dimensional \cite{Budd2010} or rectangular two-dimensional domains \cite{Budd2006,Ceniceros2001,PMA2009}. The PMA moving mesh approach has recently \cite{LD2017} been extended to the fourth-order PDE problem \eqref{eq:intro} by constructing a high regularity mapping between the computational and physical domains. The study \cite{LD2017} was based on a finite difference discretization of the PMA equation that restricted computations to rectangular domains. The main contribution of the present work is a robust adaptive numerical method that can resolve the singularities of \eqref{eq:intro} in general non-simply connected two-dimensional regions such as those utilized in real MEMS devices (cf.~Fig.~\ref{fig:intro_diagram_B}).

The outline of the paper is as follows. In Sec.~\ref{sec:skeleton} we outline a geometric theory for predicting the location of singularities based on a singular perturbation analysis of \eqref{eq:intro} as $\eps\to0$. We find that singularities are more likely to form on a set $\skel\subset\Omega$ known as the \emph{skeleton}. The skeleton of $\Omega$ is defined roughly as the collection of points $\bx\in\Omega$ at which inward facing normal vectors meet at points equidistant to $\partial\Omega$. This geometric construction is a unique minimal  representation of the domain $\Omega$ and is widely used in computer vision to store two- or three-dimensional objects \cite{Malandain1998}.

In Sec.~\ref{SEC:fem} we describe a precision numerical tool for exploring the touchdown set of equation \eqref{eq:intro} in general two-dimensional geometries. The adaptive strategy underpinning this method is a moving mesh partial differential equation (MMPDE) which dynamically relocates the mesh points to increase the density in regions where the solution has fine scale detail requiring additional spatial resolution. A notable strength of this method is the ability to resolve solutions very close to singularity in non-convex and non-simply connected geometries. A finite element method is used to discretize the system \eqref{eq:intro} on the dynamic mesh.

In Sec.~\ref{SEC:numericalresults} we employ the numerical method to investigate the singularity set of \eqref{eq:intro} for a variety of symmetric, asymmetric, and non-simply connected domains. In the case of the rectangular domain with two point symmetries, we find the singularity set to be well described by the skeleton $\skel$. In particular, we observe that the symmetries of the domain allow for touchdown to occur at multiple locations simultaneously. In domains without symmetries, we observe that the set of possible touchdown locations is more complex and that single point touchdown is the expectation for fixed $\eps$ values. In such cases we find the analytical $\skel$ to provide a good qualitative prediction of the set of potential contact locations. Finally in Sec.~\ref{sec:discussion} we summarize the results of the paper and highlight areas for future work.

\section{A geometric theory for singularity set prediction}\label{sec:skeleton}

In this section, we use asymptotic analysis in the limit as $\eps\to0$ to establish a prediction of the singularity set of \eqref{eq:intro}. This theory explains the sensitivity of the contact set and the multiplicity of touchdown points on the parameter $\eps$ and the geometry $\Omega$.

\subsection{Asymptotic analysis}

In the leading order analysis of \eqref{eq:intro} as $\eps\to0$, it is assumed that the term $-\eps^2 \Delta^2 u$ is negligible almost everywhere, except in the vicinity of $\partial\Omega$. This suggests that the solution is largely spatially uniform satisfying $u(\bx,t) \sim u_0(t)$, where $u_0(t)$ is the solution of the initial value problem
\bsub\label{eq:outer1}
\begin{equation}\label{eq:outer1_a}
\frac{d u_0}{d t} = - \frac{1}{(1+u_0)^2}, \quad t \in(0,T_0), \qquad u_0(0) = 0.
\end{equation}
The solution of \eqref{eq:outer1_a} is
\begin{equation}\label{eq:outer1_b}
u_0(t) = -1 + (1 - 3t)^{\frac13}, \qquad t \in\Big(0,\frac13\Big).
\end{equation}
\esub
This gives a leading order approximation of the singularity time as $T_0 = \frac13$. We remark upon the quenching phenomenon whereby $u_0$ is finite as $t\to T_0^{-}$ while $u_{0t}$ diverges. Clearly \eqref{eq:outer1_b} does not satisfy the boundary condition $u = 0$ on $\partial\Omega$ which must be enforced in a boundary layer. To analyze this layer for a general geometry, we introduce an orthogonal coordinate system ($\rho,s)$ where $\rho=\mbox{dist}(\bx,\partial \Omega) > 0$, while $s$, for $\bx\in\partial\Omega$, denotes the arc-length along $\partial\Omega$. In this coordinate system, \eqref{eq:intro} becomes
\bsub\label{eq:outer2}
\begin{align}
\label{eq:outer2a} u_t &= -\eps^2 \left[ \partial_{\rho\rho} - \frac{\kappa}{1- \kappa \rho}  \partial_{\rho} + \frac{1}{1- \kappa \rho} \partial_s \left(\frac{1}{1-\kappa \rho}\partial_s \right)\right]^2 u - \frac{1}{(1+u)^2}, \qquad \rho >0;\\
\label{eq:outer2b} u &=  \left[\partial_{\rho\rho} - \frac{\kappa}{1- \kappa \rho}  \partial_{\rho} + \frac{1}{1- \kappa \rho} \partial_s \left(\frac{1}{1-\kappa \rho}\partial_s \right)\right]u = 0, \qquad \rho = 0,
\end{align}
\esub
where $\kappa = \kappa(s)$ is the curvature of $\partial\Omega$. To analyze the boundary layer in this new coordinate system, we introduce the stretched variables
\begin{equation}\label{eq:outer3}
u = f(t)\, w(z), \qquad z = \frac{\rho}{\phi(t;\eps)}, \qquad \phi(t;\eps) = \eps^{\frac12}f(t)^{\frac14}, \qquad f(t) = -u_0(t) = 1 - (1-3t)^{\frac13} .
\end{equation}
The variables \eqref{eq:outer3} are substituted into \eqref{eq:outer2} and the resulting system is expanded in the form
\begin{equation}\label{eq:outer4}
w(z) = w_0(z) + \phi\, w_1(z) + \cdots.
\end{equation}
Collecting terms at each order gives a sequence of problems for $\{w_0, w_1,\ldots\}$. At leading order we have
\bsub\label{eq:outer5}
\begin{gather}
\label{eq:outer5a} w_{0zzzz} - \frac{z}{4} w_{0z} + w_0 = -1, \quad z>0;\\[5pt]
\label{eq:outer5b} w_0 = w_{0zz} = 0, \quad z = 0; \qquad w_0 \sim -1, \quad z \to \infty.
\end{gather}

\begin{figure}[htbp]
\centering
\includegraphics[width=0.45\textwidth]{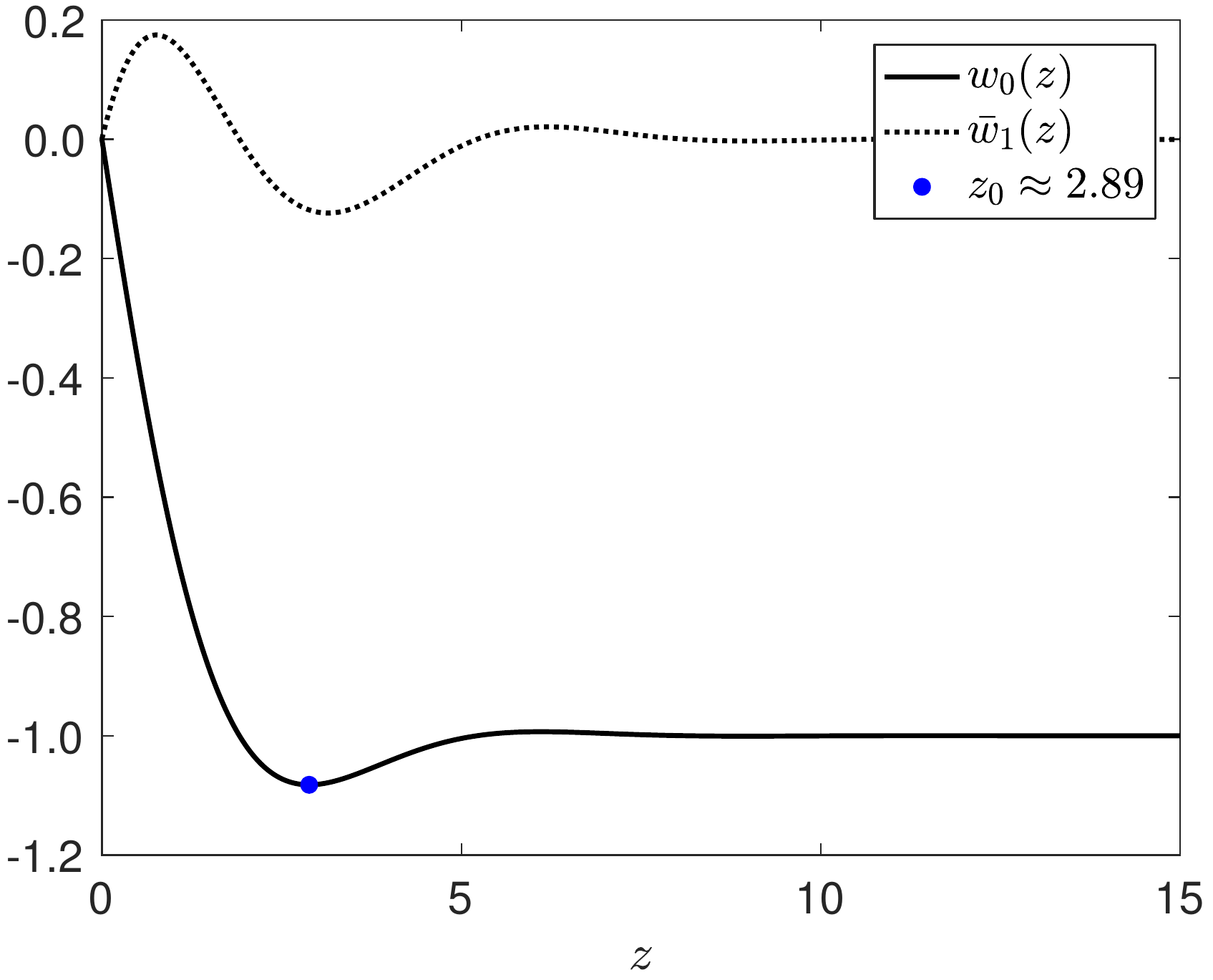}
\caption{The two profiles $w_0(z)$ and $\bar{w}_1(z)$ satisfying \eqref{eq:outer5} and \eqref{eq:outer6}.  \label{fig:twoprofiles} }
\end{figure}

While the solution of \eqref{eq:outer5} can be developed in terms of hypergeometric functions, the resulting expression is quite cumbersome and not particularly useful. The most important property is the behavior of $w_0(z)$ as $z\to\infty$ which can be derived from a WKB analysis \cite{selfsimilar2014}. In particular,
\begin{equation}\label{eq:outer5c}
w_0(z) \sim -1 + Ae^{-\omega z^{\frac43}} \sin [\sqrt{3}\, \omega z^{\frac43} + \psi] + \cdots, \qquad z\to\infty,
\end{equation}
\esub
where $\omega = 3\cdot 2^{-\frac{11}{3}}$ and $A,\psi$ are constants. The crucial observation from the limiting behavior \eqref{eq:outer5c} is the oscillatory decay for $z$ large. This phenomenon is a manifestation of the lack of a maximum principle for fourth-order equations. In particular, $w_0(z)$ attains its global minimum at a finite value, which can be approximated numerically as $z_0\approx 2.89$ - see Fig.~\ref{fig:twoprofiles}.

At the next order, we apply the decomposition $w_1(z,s)  = \kappa(s) \bar{w}_1(z)$ and find that $\bar{w}_1(z)$ satisfies
\bsub\label{eq:outer6}
\begin{gather}
\label{eq:outer6a} \bar{w}_{1zzzz} - \frac{z}{4} \bar{w}_{1z} + \frac{5}{4}\bar{w}_1 = 2w_{0zzz}, \quad z>0;\\[5pt]
\label{eq:outer6b} \bar{w}_1 = \bar{w}_{1zz} - w_{0z} = 0, \quad z = 0; \qquad \bar{w}_1 \sim 0, \quad z \to \infty.
\end{gather}
\esub
The two profiles $w_0(z)$ and $\bar{w}_1(z)$ are shown in Fig.~\ref{fig:twoprofiles}.

The lack of monotonicity in the profiles of the stretching boundary layer lowers the value of the solution at certain points. By superimposing the boundary layer solution with the flat solution $u_0(t)$, and subtracting the overlap term, the following global solution at a particular $\bx\in\Omega$ is
\begin{equation}\label{eqn:mainexp}
u(\bx,t) = u_0(t) - u_0(t)\sum_{j=1}^{N} \left( 1+ w_0\left(\frac{|\bx - \by_j|}{\phi} \right) + \phi \,\kappa(\by_j) \bar{w}_1\left(\frac{|\bx - \by_j|}{\phi} \right) +\cdots \right)
\end{equation}
where $\phi = \eps^{\frac12}|u_0(t)|^{\frac14}$. For each $\bx\in\Omega$, the boundary points $\{\by_1,\ldots,\by_N\}\in\partial\Omega$ are those with inward facing normal vectors that pass through $\bx$, i.e., points such that the straight line between $\bx$ and $\by_j$ is contained in $\Omega$ and meets $\partial\Omega$ orthogonally. The quantity $\kappa(\by_j)$ is the local boundary curvature at the point $\by_j$. The asymptotic expansion \eqref{eqn:mainexp} is valid for $\phi\ll1$ which corresponds to short times $t\ll1$. In this regime, the solution is composed of a flat central region coupled to propagating boundary interfaces.

\subsection{The skeleton of the domain}\label{sec:domain}

\begin{figure}[htbp]
\centering
\subfigure[$\omega(t)$]{\includegraphics[width=0.35\textwidth]{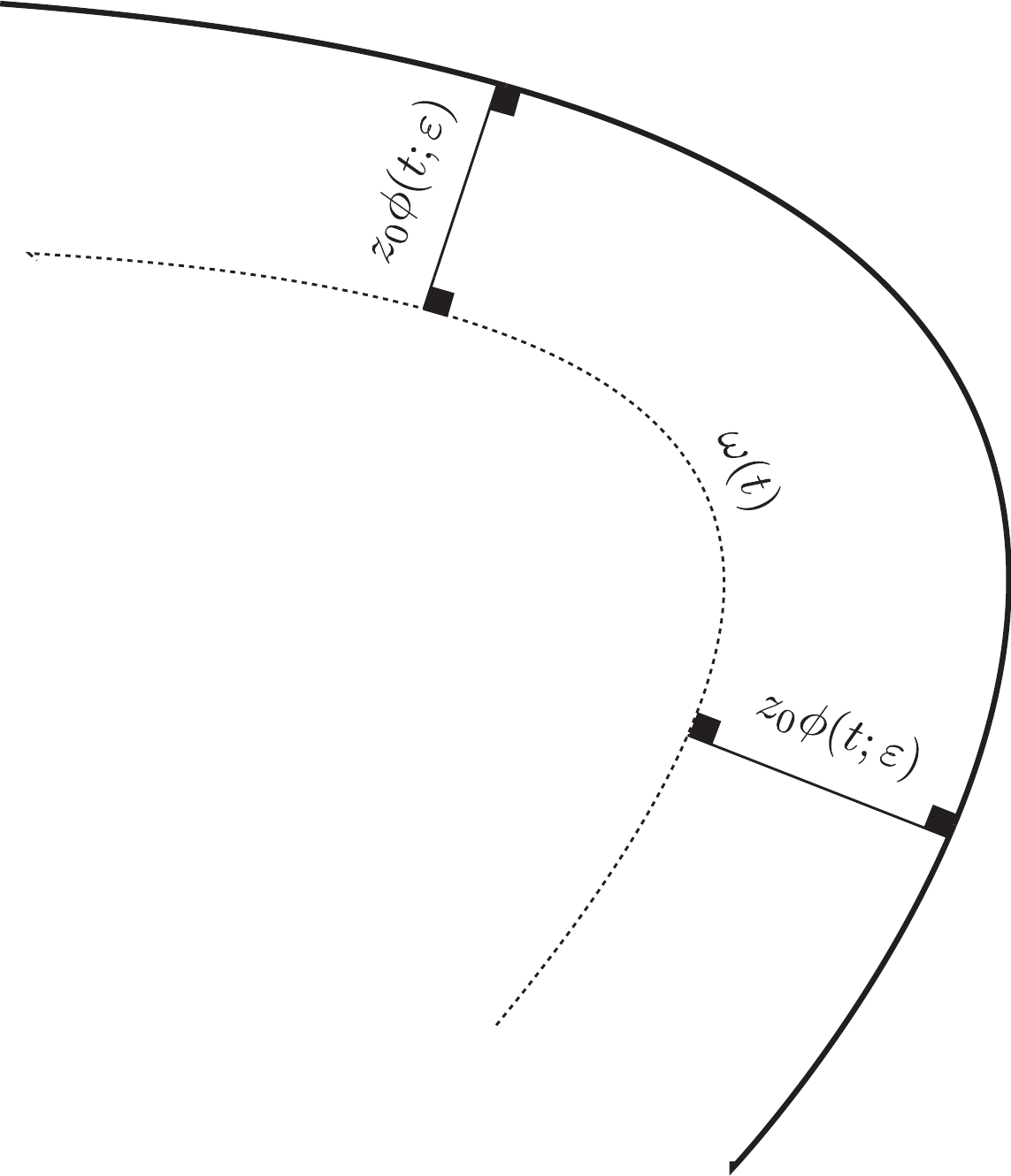} \label{fig:omega_t}}\hspace{0.45in}
\subfigure[$\skel$]{\includegraphics[width=0.35\textwidth]{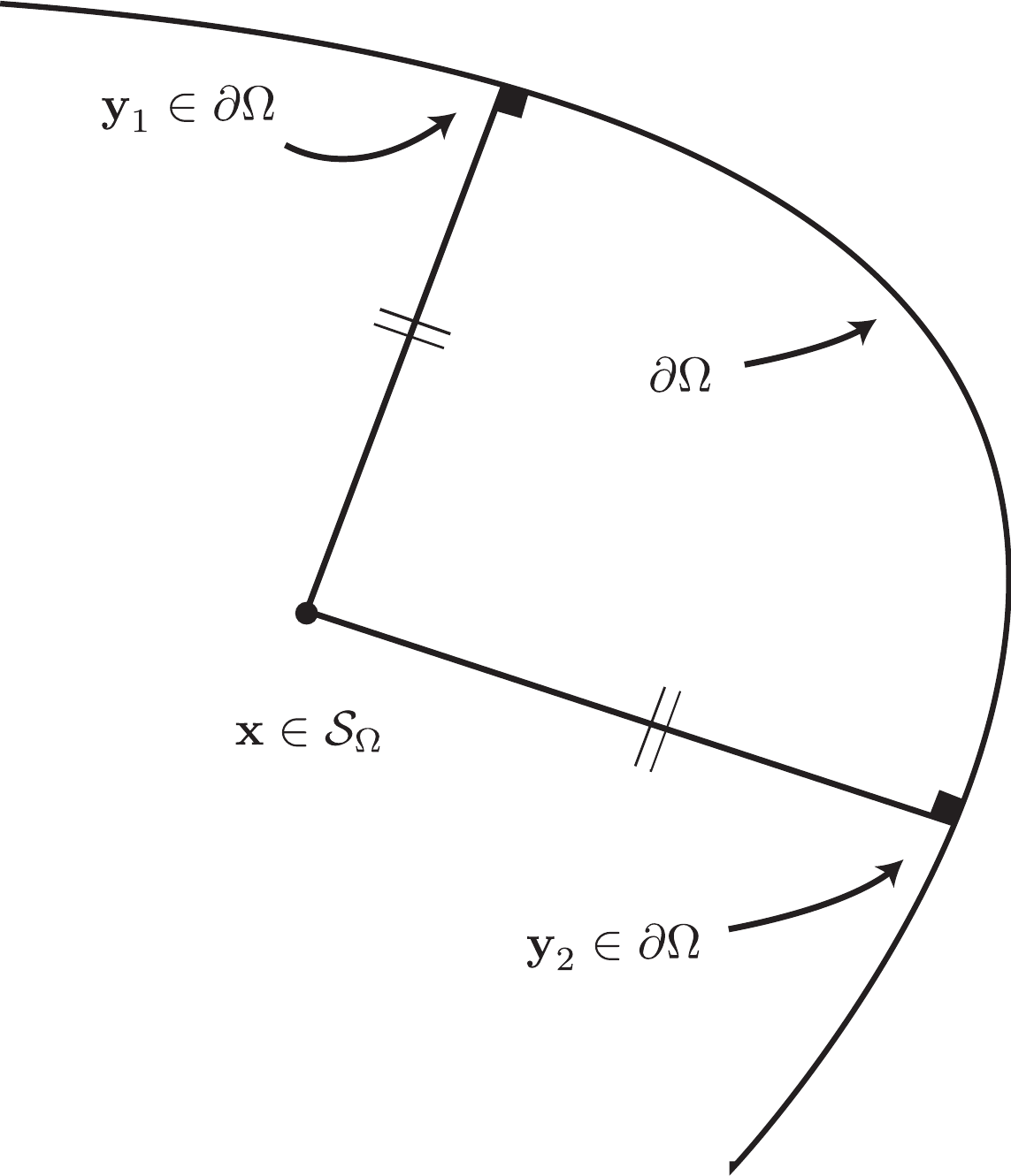}  \label{fig:skel}}
\caption{The two sets $\omega(t)$ and $\skel$ on which touchdown may occur. \label{fig:om_skel}}
\end{figure}

We now use the asymptotic solution \eqref{eqn:mainexp} to develop a predictive theory for how the geometry $\Omega$ and $\eps$ combine to select possible singularity locations. The key reason for the multiple singularity phenomenon (cf.~Fig.~\ref{fig:introsol}) is the non-monotonicity of the profile $w_0(z)$ which lowers the value of the solution at certain points in the domain and promotes faster touchdown there. As shown in Fig.~\ref{fig:twoprofiles}, the profile $w_0(z)$ has a unique global minimum at $z=z_0$ whose value can be estimated numerically as $z_0\approx2.89$. In light of this, and the arguments of the solution \eqref{eqn:mainexp}, the set of points
\begin{equation}\label{eq:omega_t}
\omega(t) = \{ \bx\in\Omega \ \mid \ \mbox{dist} (\bx,\partial\Omega) = z_0\phi(t;\eps) \}
\end{equation}
are particularly important. The condition \eqref{eq:omega_t} describes a curve of points (cf.~Fig.~\ref{fig:omega_t}) that extend inwards from $\partial\Omega$ a distance $z_0\phi(t;\eps)$ and along which the solution of \eqref{eq:intro} is, to first order in $\phi$, at a local minimum. In computer vision literature \cite{Malandain1998}, this set is known as the \emph{firefront}.
In a radially symmetric scenario for which the boundary is of uniform curvature, the singularities may form simultaneously along a ring of points \cite{Witelski00,selfsimilar2014}. For domains whose boundaries have non-uniform curvature, the effect is to promote touchdown at certain points rather than along entire curves. This can be deduced from the asymptotic solution \eqref{eqn:mainexp} by seeking a regular expansion solution $z_{\text{min}} = z_0 + \phi z_1 + \cdots$ of the equation $\nabla u = 0$. This reveals that
\bsub\label{eq:min2}
\begin{equation}\label{eq:min2a}
z_{\text{min}} = z_0 - \alpha  \phi \frac{1}{N}\sum_{j=1}^{N} \kappa(\by_j)  + \bigoh(\phi^2), \qquad \alpha = -\frac{\bar{w}_{1z}(z_0)}{w_{0zz}(z_0)} \approx 0.3533.
\end{equation}
The corresponding asymptotic prediction of the minimum is found from  \eqref{eqn:mainexp} to be
\begin{equation}\label{eq:min2b}
u(\bx,t)\mid_{z=z_{\text{min}}} = u_0(t) - u_0(t) \sum_{j=1}^{N} \Big( 1+ w_0(z_0) + \phi \,\kappa(\by_j) \bar{w}_1 (z_0) + \cdots \Big),
\end{equation}
where numerically we determine the values
\begin{equation}\label{eq:min3b}
w_0(z_0) = -1.0822 , \qquad \bar{w}_1(z_0) = -0.1186.
\end{equation}
\esub
Since $\bar{w}_1(z_0) <0$, we conclude that the solution will take lower values at points of $\omega(t)$ whose contributing boundary points correspond to maxima of the boundary curvature $\kappa(s)$.

As $t$ increases and the curve $\omega(t)$ propagates, it may self-intersect at some time $t= \sat$. If this occurs, the solution minimum \eqref{eq:min2b} goes through a distinct change since the number of contributing boundary points, $N$, increases. For example, in the scenario displayed in Fig.~\ref{fig:om_skel}, the set $\omega(t)$ eventually intersects the point $\bx\in\Omega$ which receives boundary contributions from the two points $\{ \by_1, \by_2\}$ and the number of boundary contributions increases from $N=1$ to $N=2$. These points are important as multiple boundary contributions arrive simultaneously and superimpose to lower the value of \eqref{eq:min2b} further - this set of points is called the \emph{skeleton} of the domain and denoted $\skel$. The time $\sat$ is then known as \emph{the skeleton arrival time} and can be defined explicitly as
\begin{equation}\label{eq:sat}
\sat = \inf_{t\geq0} \left\{ t \ \mid \ \exists \bx\in \skel, \  \mbox{dist}(\bx,\partial\Omega) = z_0 \phi(t;\eps) \right\}.
\end{equation}

As indicted in Fig.~\ref{fig:skel}, a point $\bx\in\skel$ if it has at least two closest boundary points, i.e., $\mbox{dist}(\bx,\partial\Omega) = \mbox{dist}(\bx,\by_j)$ for $\{\by_1,\ldots,\by_N \}\in\partial\Omega$ and $N\geq2$. The skeleton $\skel\subset\Omega$ is a minimal representation of the domain $\Omega$. They are homotopic to one another so that each $\Omega$ defines a unique $\skel$ and vice-versa \cite{Malandain1998}. A particular point $\bx\in\skel$ can potentially be associated with multiple boundary distances in which case the most pertinent one is the shortest as that is where the trough associated with $w_0(z)$ will reach first.

In summary the asymptotic analysis makes the following prediction for the touchdown set. Let $\skel$ be the skeleton of the domain, $\sat\geq0$ be the skeleton arrival time and define $T$ to be maximum global existence time of \eqref{eq:intro}. Then, the leading order asymptotic analysis predicts the following the dichotomy of possibilities:
\begin{enumerate}
\item If $T <\sat$, the singularities form on $\omega(t)$ at point(s) corresponding to maximum boundary curvature.
\item If $T \geq \sat$, the singularities form on $\skel$ at points $\bx\in\skel$ satisfying $\mbox{dist}(\bx,\partial\Omega) = \phi(T;\eps)$.
\end{enumerate}
% section
In the following sections, we explore the predictive ability of this asymptotic framework for describing the potential touchdown sets of \eqref{eq:intro}. In the following section, we describe a new moving mesh finite element method for resolving solutions of \eqref{eq:intro} very close to singularity. In Sec.~\ref{SEC:numericalresults}, we demonstrate the asymptotic and numerical methodologies on a variety of two-dimensional regions $\Omega$.

\section{The adaptive moving mesh finite element method}\label{SEC:fem}

In this section we describe an adaptive MMPDE finite element method for solving
the MEMS system. We apply the MMPDE method to dynamically concentrate the mesh nodes
around the places where the solution is touching down.

It has been proven in \cite{HK2015} that the mesh governed by the MMPDE
stays non-singular if it is non-singular initially. This result holds at semi- and fully-discrete levels and
for convex or concave domains. Numerical examples will be presented in Section~\ref{SEC:numericalresults}.

It should be pointed out that a number of other moving mesh methods have been developed in the past
and there is a vast literature in the area. The interested reader is referred to the books or review articles \cite{Bai94a,Baines-2011,BHR09,HR11,Tan05} and references therein.

% section
\subsection{Finite element discretization}

We now describe the finite element approximation of the PDE system \eqref{eq:intro} up to a finite time instant $T$.
We introduce the auxiliary variable $v$ defined as
\[
v = \Delta u.
\]
Then, the fourth-order PDE \eqref{eq:intro} can be written as the second-order system
\begin{equation}
\begin{cases}
u_t = -\eps^2 \Delta v - \ds\frac{1}{(1+u)^2}, & (\bx, t) \in \Omega \times (0,T) \\[5pt]
v = \Delta u, & (\bx, t) \in \Omega \times (0,T)\\[5pt]
u = v = 0, & (\bx, t) \in \partial \Omega \times (0,T) \\[5pt]
u(\bx,0) = v(\bx,0) = 0, & \bx \in \Omega .
\end{cases}
\label{MEMsv}
\end{equation}
The computation alternates between the integration of the PDE and the mesh equation.
Assume that we are given time instants
\[
0 = t_0 < t_1 < \ldots < t_{n_f} = T
\]
and the physical mesh $ \mathcal{T}_h^n$ and the numerical solution $u_h^n(\cdot,t)$
and $v_h^n(\cdot,t)$ defined thereon at $t_n$.
The new physical mesh $ \mathcal{T}_h^{n+1}$ is first generated by an MMPDE-based strategy
(to be described in Section~\ref{SEC:MMPDE}) and then the physical PDEs are integrated
from $t_n$ to $t_{n+1}$.  The procedure is repeated until $T$ is reached.
The number of the mesh elements and the mesh connectivity are fixed throughout the computation.

Denote the coordinates of the vertices of $ \mathcal{T}_h^{n}$ and $ \mathcal{T}_h^{n+1}$ by
$\bx_j^n$ and $\bx_j^{n+1}$, $j = 1, 2, \ldots, N_v$, respectively. We define the coordinates of the vertices
between $t_{n}$ and $t_{n+1}$ as
\[
\bx_j(t) = \frac{t - t_n}{t_{n+1} - t_n} \bx_j^{n+1}+ \frac{t_{n+1} - t}{t_{n+1} - t_n}\bx_j^n,\quad j = 1,\ldots,N_v,
\quad t \in [t_n, t_{n+1}].
\]
The corresponding mesh is denoted by $\mathcal{T}_h(t)$.
Then, a linear finite element approximation for (\ref{MEMsv}) is to find $u_h(\cdot,t)$, $v_h(\cdot,t) \in V_h^0(t)$,
for $t \in (t_0, T]$, such that
\begin{equation}
\begin{cases}
\ds\int_{\Omega} \frac{\partial u_h}{\partial t} \psi\, d\bx - \eps^2 \ds\int_{\Omega} \nabla v_h \cdot \nabla \psi\, d\bx + \ds\int_{\Omega} \ds\frac{\psi}{(1+u_h)^2} d\bx = 0, & \forall \psi \in V_h^0(t) \\[10pt]
\ds\int_{\Omega} v_h \psi\, d\bx + \int_{\Omega} \nabla u_h \cdot \nabla \psi\, d\bx = 0, & \forall \psi \in V_h^{0}(t)
\end{cases}
\label{mem-fem}
\end{equation}
where $V_h^0(t)$ is the span of the linear basis functions that are compactly supported on $\mathcal{T}_h(t)$
at $t$. Notice that linear basis functions and the linear finite element function space are time dependent.
For simplicity, we assume that the first $N_{v_i}$ out of $N_v$ vertices are interior vertices.
Denoting the linear basis function associated with the vertex $\bx_j$ by $\psi_{j}(\cdot,t)$, $V_h^0(t)$
can be expressed as
\[
V_h^0(t) = \text{span}\{\psi_1(\cdot,t),\, \ldots,\,\psi_{N_{v_i}}(\cdot,t)\}.
\]

With the linear basis functions being time dependent, the main difference between the integration of (\ref{mem-fem})
from that on a fixed mesh lies in the term $\frac{\partial u_h}{\partial t}$. To see this, expressing $u_h$ as
\begin{align}
u_h(\bx, t) = \sum_{i = 1}^{N_{v_i}} u_i(t) \psi_i(\bx, t)
\label{udis}
\end{align}
and differentiating it with respect to time, we get
\[
\frac{\partial u_h(\bx, t)}{\partial t} = \sum_{ i = 1}^{N_{v_i}} \frac{du_i}{dt}\psi_i(\bx, t)
+ \sum_{i = 1}^{N_{v_i}}u_i(t) \frac{\partial \psi_i(\bx,t)}{\partial t}.
\]
It has been proven (e.g., see \cite{HR11}) that
\[
\frac{\partial \psi_i}{\partial t} = -\nabla \psi_i \cdot \dot{\bX}, \quad\text{a.e. in } \Omega
\]
where the mesh velocity $\dot{\bX}$ is defined as
\[
\dot{\bX} = \sum_{i = 1}^{N_v} \dot{\bx}_i \psi_i(\bx,t)
\]
and the term $\dot{\bx}_i$ denotes the nodal mesh speed. Combining the results above, we get
\[
\frac{\partial u_h}{\partial t} = \sum_{i = 1}^{N_{v_i}} \frac{du_i}{dt} \psi_i - \nabla u_h \cdot \dot{\bX}.
\]

Inserting these into (\ref{mem-fem}) and taking $\psi = \psi_j$ successively,  we can rewrite
(\ref{mem-fem}) into a system of differential-algebraic equations in the form as
\begin{equation}
\begin{cases}
M(\bX) \dot{\bU} = \eps^2 B(\bX)\bV + F(\bX,\dot{\bX},\bU),\\[5pt]
\textbf{0} = M(\bX) \bV + B(\bX)\bU ,
\end{cases}
\label{simode}
\end{equation}
where $\bX$ is a vector representing the location of the vertices, $M(\bX)$ is the mass matrix, $B(\bX)$ is the stiffness matrix and $\bU$, $\bV$ are vectors of the unknown nodal values. This system for $\bU$ and $\bV$ is integrated from
$t_n$ to $t_{n+1}$ using the fifth-order Radau IIA method (e.g., see Hairer and Wanner \cite{HW96}),
and the time step is chosen by a standard selection procedure \cite{HW96} with a two-step error
estimate of Gonzalez-Pinto et al. \cite{Montijano2004}.

% section
\subsection{The MMPDE moving mesh strategy}
\label{SEC:MMPDE}

We now describe the MMPDE moving mesh strategy to generate $\mathcal{T}_h^{n+1}$,
given the physical mesh $\mathcal{T}_h^n$ and the computed solutions $u_h(t_n)$
and $v_h(t_n)$ at time step $t = t_n$. We use here a discrete approach of \cite{HK2014}
for the MMPDE method.

We shall use three meshes, the physical mesh $\mathcal{T}_h = \{ \bx_1, \ldots, \bx_{N_v} \}$,
the computational mesh $\mathcal{T}_{c,h} = \{ \bxi_1,\ldots,$ $\bxi_{N_v} \}$, and
the reference computational mesh $\hat{ \mathcal{T}}_{c,h} = \{ \hat{\bxi}_1,\ldots, \hat{\bxi}_{N_v} \}$,
with all of them having the same number of elements and the same connectivity.
Typically, $\hat{\mathcal{T}}_{c,h}$ is chosen to be a mesh as uniform as possible (under the Euclidean metric)
and kept fixed throughout the computation. We may take $\mathcal{T}_h = \mathcal{T}_h^n$ or
$\mathcal{T}_{c,h} = \hat{ \mathcal{T}}_{c,h}$, depending on which formulation, the $\bxi$- or $\bx$-formulation,
we use for generating $\mathcal{T}_h^{n+1}$. (See the description below.)
Notice that for each $K \in \mathcal{T}_h$, there exists a unique corresponding element $K_c \in \mathcal{T}_{c,h}$.
The affine mapping between the elements is denoted by $F_K: K_c \rightarrow K$ and its Jacobian matrix
by $F_K'$.

The MMPDE method employs a metric tensor $\M = \M(\bx)$ to specify the size, shape, and orientation
of the mesh elements throughout the domain. Here we always assume that $\M$ is symmetric and
uniformly positive definite on $\Omega$. We define $\M$ as a piecewise constant function as
\begin{equation}
\M_K = \det(I + \alpha_h^{-1} |H_K|)^{-\frac{1}{d+4}}(I + \alpha_h^{-1} |H_K|), \quad \forall K \in \mathcal{T}_h
\label{M-1}
\end{equation}
where $H_K$ is an approximate Hessian of $u_h$ on element $K$ that is obtained using
a least-squares Hessian recovery technique,
$|H_K| = Q \text{diag}\,\,(|\lambda_1|,\ldots,|\lambda_d|)Q^T$, with the eigen-decomposition of $H_K$ being
$Q \text{diag}(\lambda_1,\ldots,\lambda_d)$ $Q^T$, and $\alpha_h$ is chosen such that
\[
\sum \limits_{K\in \mathcal{T}_h} |K| \sqrt{\det(\M_K)} = 2 |\Omega| .
\]
The choice (\ref{M-1}) of $\M$ is known to be optimal
with respect to the $L^2$ norm of the linear interpolation error \cite{Hua05b},
with the expectation that the mesh points will be concentrated around the regions where the recovered Hessian
of $u_h$ has a large determinant.

The main idea of the MMPDE method is viewing any adaptive mesh $\mathcal{T}_h$ as a uniform one
under the metric $\M$ and in reference to the mesh $\mathcal{T}_{c,h}$.
Geometrically, this is equivalent to the requirements that the volume of each $K \in \mathcal{T}_h$
be of the same multiple of the volume of $K_c\in \mathcal{T}_{c,h}$ and that $K$ be similar to $K_c$,
all in the metric $\M$. Recall that the distance between two points $\bx$ and $\bx + d\bx$
under the metric tensor $\M$ is defined as
\[
\sqrt{d\bx^T \M(\bx) d\bx }.
\]
Then the above requirements can be approximately expressed as the equidistribution and alignment conditions
(e.g., see \cite{HR11})
\begin{align}
& |K|\sqrt{\det(\M_K)}  = \frac{\sigma_h |K_c|}{|\Omega_c|}, \qquad \forall K \in \mathcal{T}_h
\label{equi}
\\
& \frac{1}{d}\text{tr}\left ( (F_K')^{-1}\M_K^{-1}(F_K')^{-T}\right ) =
\det \left ( (F_K')^{-1}\M_K^{-1}(F_K')^{-T}\right )^{\frac{1}{d}},
\qquad \forall K \in \mathcal{T}_h
\label{align}
\end{align}
where $|K|$ and $|K_c|$ denote the volume of $K$ and $K_c$, respectively, $d$ is the dimension of $\Omega$
($d = 2$ for the current situation),
and $|\Omega_c|$ and $\sigma_h$ denote the total volume of the computational
mesh (in the Euclidean metric) and the physical mesh (in the metric $\M$), respectively,
\[
|\Omega_c| = \sum_{K_c \in \mathcal{T}_{c,h}} |K_c| ,\qquad
\sigma_h = \sum \limits_{K\in \mathcal{T}_h} |K| \sqrt{\det(\M_K)}.
\]
An energy function associated with these conditions \cite{Hua01b} is given by
\begin{align}
I_h(\mathcal{T}_h, \mathcal{T}_{c,h})  = & \; \frac{1}{3} \sum_{K \in \mathcal{T}_h} |K| \sqrt{\det(\M_K)}
\left (\mathbf{\text{tr}}( (F_k')^{-1}\M_K^{-1}(F_K')^{-T}) \right )^{\frac{3 d}{4}}
\notag \\
+ & \; \frac{d^{\frac{3 d}{4}}}{3}    \sum_{K \in \mathcal{T}_h} |K| \sqrt{\det(\M_K)}
\left (\frac{|K_c|}{|K|\sqrt{\det(\M_K)}}\right )^{\frac{3}{2}} .
\label{ih}
\end{align}
Minimization of this energy function tends to produce a mesh satisfying (\ref{equi}) and (\ref{align}).

We integrate the gradient system of the energy function for minimization. In the $\bx$-formulation,
we take $\mathcal{T}_{c,h} = \hat{\mathcal{T}}_{c,h}$ and $I_h$ is a function of the coordinates
of the physical vertices. The MMPDE mesh equation is defined as
\begin{equation}
\frac{\partial \bx_i}{\partial t} = - \frac{P_i}{\tau} \left (\frac{\partial I_h}{\partial \bx_i}\right )^T, \quad i = 1,\ldots, N_v,
\quad t \in (t_n, t_{n+1}],
\label{xmethod}
\end{equation}
where $\frac{\partial I_h}{\partial \bx_i}$ is a row vector, $P_i$ is a positive function chosen
as $P_i = \det(\M(\bx_i))^{\frac{1}{4}}$ (to make (\ref{xmethod}) to be invariant under the scaling transformation of $\M$),
and $\tau > 0$ is a positive parameter used to adjust the response time of mesh movement to the changes in $\M$.
It has been proven in \cite{HK2015} that the mesh governed by (\ref{xmethod}) stays non-singular if it is
non-singular initially. This result holds for any convex or concave domain in any dimension
and for the semi-discrete form (\ref{xmethod}) or a fully-discrete form of (\ref{xmethod}).
(In the latter case, the time step is required to be sufficiently small but not diminishing.)
The drawback of this formulation is that $\M$, a function of $\bx$, needs to be constantly updated during
the integration, which can be costly especially in higher dimensions.

To avoid this disadvantage, we now consider the $\bxi$-formulation with which $\mathcal{T}_h$ is taken
as $\mathcal{T}_h^{n}$ in (\ref{ih}) and $I_h$ is minimized with respect to the coordinates of the computational
vertices. Then the MMPDE mesh equation reads as
\begin{equation}
\frac{\partial \bxi_i}{\partial t} = - \frac{P_i}{\tau}\left (\frac{\partial I_h}{\partial \bxi_i}\right )^T,
\qquad i = 1, \ldots, N_v,  \quad t \in (t_n, t_{n+1}] .
\label{ximethod}
\end{equation}
This equation, with proper modifications for the boundary vertices (to keep them on the boundary),
is integrated from the initial mesh $\mathcal{T}_{c,h}(t_n) = \hat{\mathcal{T}}_{c,h}$.
The Matlab\textsuperscript \textregistered\, function {\tt ode15s},
a Numerical Differentiation Formula based integrator,  is used for this purpose in our computation.
Since $\mathcal{T}_h^n$ is fixed during the integration, there is no need of constantly reassigning
the metric tensor $\M$. The new computational mesh obtained in this way is denoted by $\mathcal{T}_{c,h}^{n+1}$.
Notice that $\mathcal{T}_h^n$ and $\mathcal{T}_{c,h}^{n+1}$ form a correspondence, i.e.,
$\mathcal{T}_h^n = \Psi_h(\mathcal{T}_{c,h}^{n+1})$.  Then, the new physical mesh at $t^{n+1}$
is defined as
\[
\mathcal{T}_h^{n+1} =  \Psi_h( \hat{\mathcal{T}}_c),
\]
which can be approximated readily by linear interpolation.

The derivative $\partial I_h/\partial \bxi_i$ in (\ref{ximethod})
can be found analytically using scalar-by-matrix differentiation and has a
relatively simple matrix form \cite{HK2014}. In this way, we can rewrite \eqref{ximethod} as
\begin{equation}
\frac{\partial \bxi_i}{\partial t} = \frac{P_i}{\tau} \sum_{K \in \omega_i} |K| v_{i_K}^K, \qquad i = 1,\ldots, N_v,
\label{ximethod-2}
\end{equation}
where $\omega_i$ is the set of all the elements having $\bx_i$ as a vertex and $v_{i_K}^K$ is the local velocity
contributed by the element $K$ to vertex $\bx_i$, with $i_K$ denoting the local index of $\bx_i$ in $K$.
The local velocities on element $K$ are given by
\begin{equation}
\begin{bmatrix}
(v_1^K)^T \\
\vdots \\
(v_d^K)^T
\end{bmatrix} = -E_K^{-1} \frac{\partial G}{\partial \det(\J)}
- \frac{\partial G}{\partial \det(\J)} \frac{\det(\hat{E}_K)}{\det(E_K)} {\hat{E}_K}^{-1},
\qquad v_0^{K} = - \sum_{ i = 1}^d v_d^K,
\label{velocity}
\end{equation}
where $E_K = [\bx_1^K-\bx_0^K, \ldots, \bx_d^K-\bx_0^K]$ and $\hat{E}_K = [\bxi_1^K-\bxi_0^K, \ldots, \bxi_d^K-\bxi_0^K]$
denote the edge matrices of $K$ and $K_c$, respectively,  $\J = (F_K)^{-1} = \hat{E}_K E_K^{-1}$, and
$G = G(\J, \det(\J), \M_K)$ is the function associated with the energy functional \eqref{ih}, i.e.,
\[
 G = \frac{1}{3} \sqrt{\det(\M_K)} (\text{tr}(\J\M_K^{-1}\J^T))^{\frac{3 d}{4}}
+ \frac{d^{\frac{3 d}{4}}}{3} \sqrt{\det(\M_K)}\left (\frac{\det(\J)}{\sqrt{\det(\M_K)}}\right )^{\frac{3}{2}}.
\]
The partial derivatives  $ \partial G/\partial \J$ (a matrix-valued function)
and $\partial G/\partial \det(\J)$ (a scalar function) can be found as
\begin{align*}
 \frac{\partial G}{\partial \J} &= \frac{d}{2} \sqrt{\det(\M_K)}(\text{tr}(\J\M_K^{-1}\J^T))^{\frac{3d}{4} - 1}\M_K^{-1}\J^T, \\[5pt]
 \frac{\partial G}{\partial \det(\J)} &= \frac{d^{\frac{3d}{4}}}{2} \det(\M_K)^{-\frac{1}{4}} \det(\J)^{\frac{1}{2}} .
\end{align*}

We note that the MMPDE equation \eqref{ximethod-2} is already discrete in space (and thus no further spatial
discretization is needed). Moreover, its computation mainly involves the calculation of the edge matrices
and matrix inversion and multiplications.

%% section
\section{Numerical results}
\label{SEC:numericalresults}

In this section, we demonstrate the efficacy of adaptive numerical methodology and the asymptotic predictions on a variety of examples. For each of the domains $\Omega$ considered, we first calculate the skeleton $\skel$ of the region defined in Sec.~\ref{sec:domain}. The numerical integration of the PDE system \eqref{mem-fem} is performed until $\min_{\bx \in \Omega} u_h =-0.99$.

\vspace{10pt}

%% example 1
 \begin{exam}[Rectangular domain]
 {\em We first consider the rectangular domain $\Omega = (-1,1)\times (-0.8,0.8)$.}
 \end{exam}

\begin{figure}[htbp]
\centering
\subfigure[Skeleton with mesh size $N = 6240\; (40 \times 39)$.]{\includegraphics[width = 0.4\textwidth]{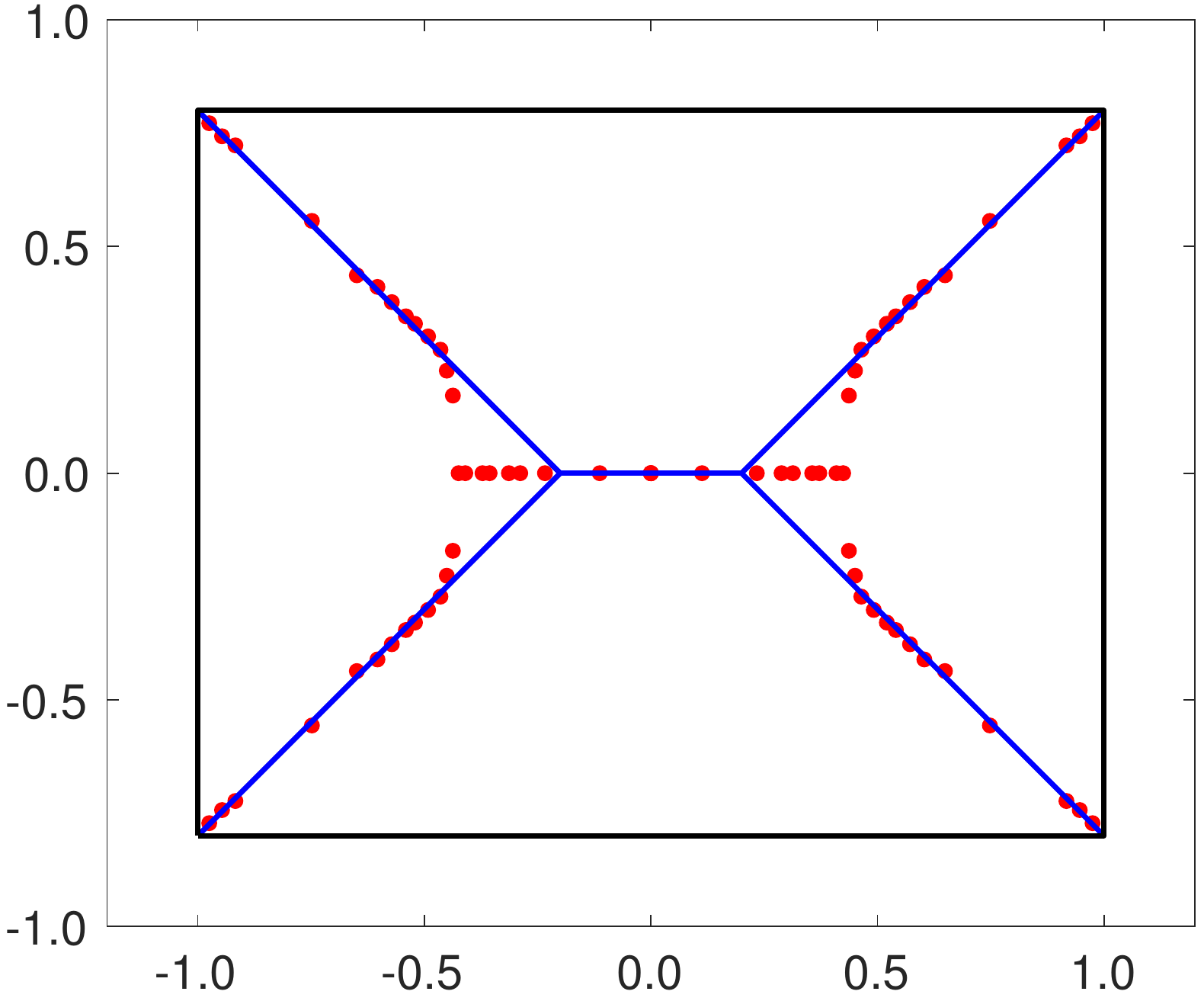}\label{noholesk8}}\quad\quad\quad
\subfigure[Skeleton with mesh size $N = 15680\; (70 \times 56)$.]{\includegraphics[width = 0.4\textwidth]{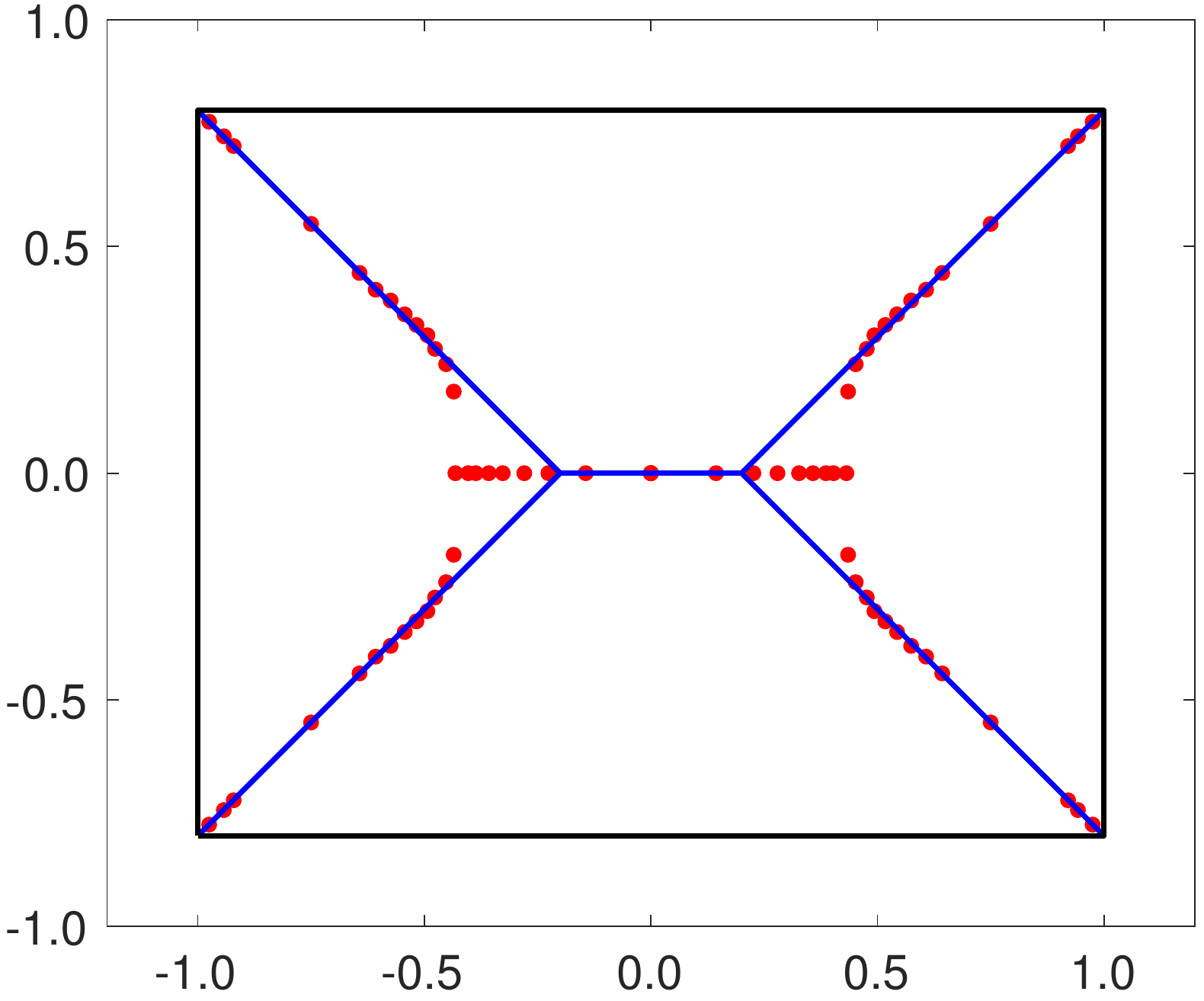}\label{noholesk14}}\quad
%\subfigure[]{\includegraphics[width = 0.35\textwidth]{rectangleSkeleton.pdf}\label{fig:rectangleSkel}}\qquad
%\subfigure[]{\includegraphics[width = 0.35\textwidth]{noholesk.pdf}\label{fig:rectangleTD}}
%\includegraphics[width = 0.4\textwidth]{noholesk.pdf}
\caption{Skeleton for rectangular domain (solid blue) with numerically obtained touchdown locations (red dots). Figs.~\ref{noholesk8} and \ref{noholesk14} show results obtained with mesh sizes $N = 6240\; (40 \times 39)$ and $N = 15680\; (70 \times 56)$ respectively.\label{noholesk}}
\end{figure}

In Fig.~\ref{noholesk} we show $\Omega$ and the skeleton $\skel$ together with numerically obtained touchdown points for the parameter range $\eps\in(10^{-4},10^{-1})$.  In Fig.~\ref{noholesk8} and Fig.~\ref{noholesk14} results are shown for meshes of size $N = 6240\; (40 \times 39)$ and  $N = 15680\; (70 \times 56)$, respectively. We observe the location of touchdown is robust as the mesh size increases. The set $\skel$ meets the boundary $\partial\Omega$ and so the skeleton arrival time satisfies $\sat=0$.

%
%\[
%\eps \in \{ 10^{-4}, \quad 10^{-3}, \quad 5\times 10^{-3}, \quad 10^{-2}, \quad linspace(0.02, 0.1, 21) \}.
%\]
%
 \begin{figure}[htbp]
\centering
\subfigure[Mesh at touchdown, $\eps= 0.02$.]{\includegraphics[width = 0.3\textwidth]{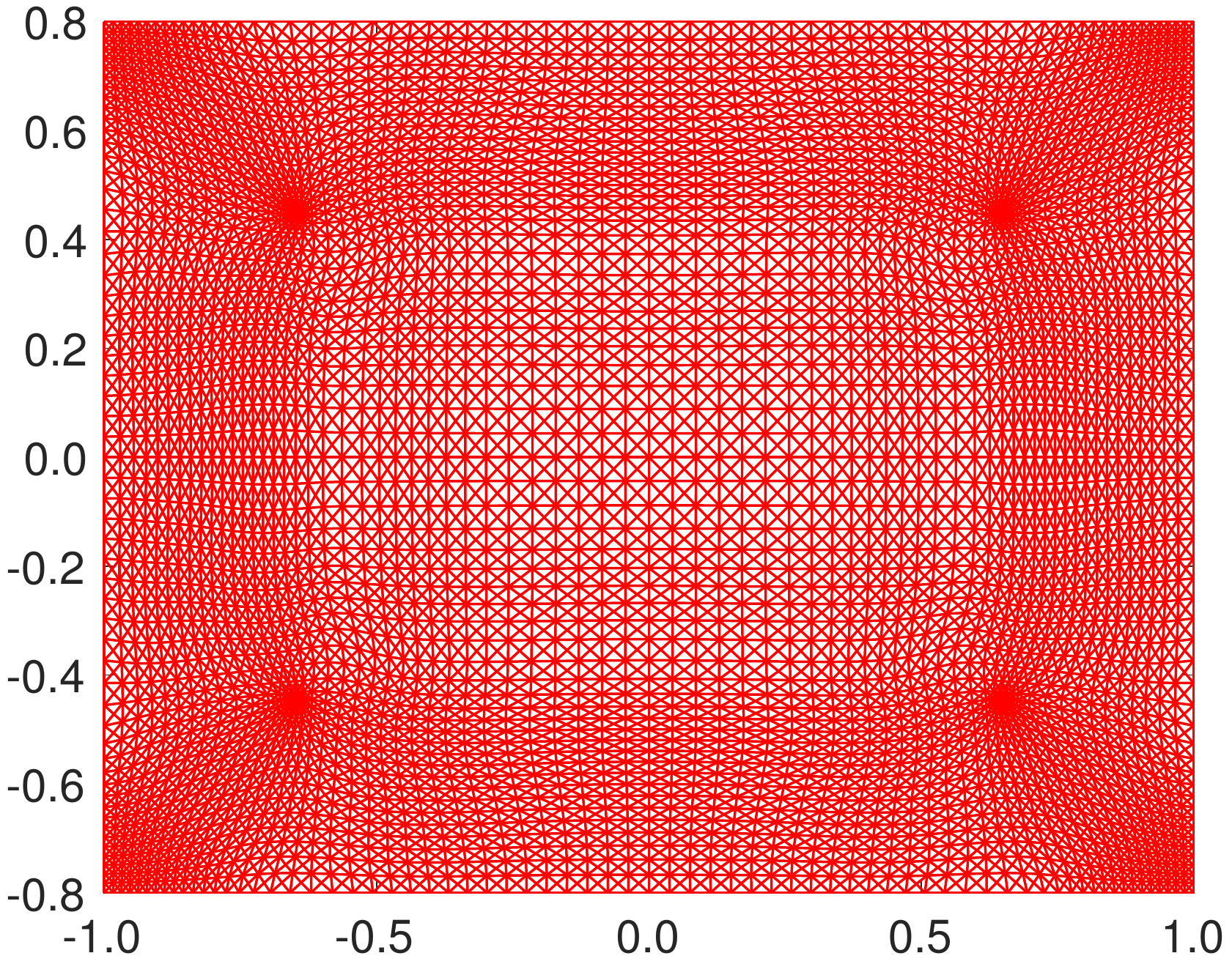} \label{nh1a}}\quad
\subfigure[Mesh at touchdown, $\eps= 0.068$.]{\includegraphics[width = 0.3\textwidth]{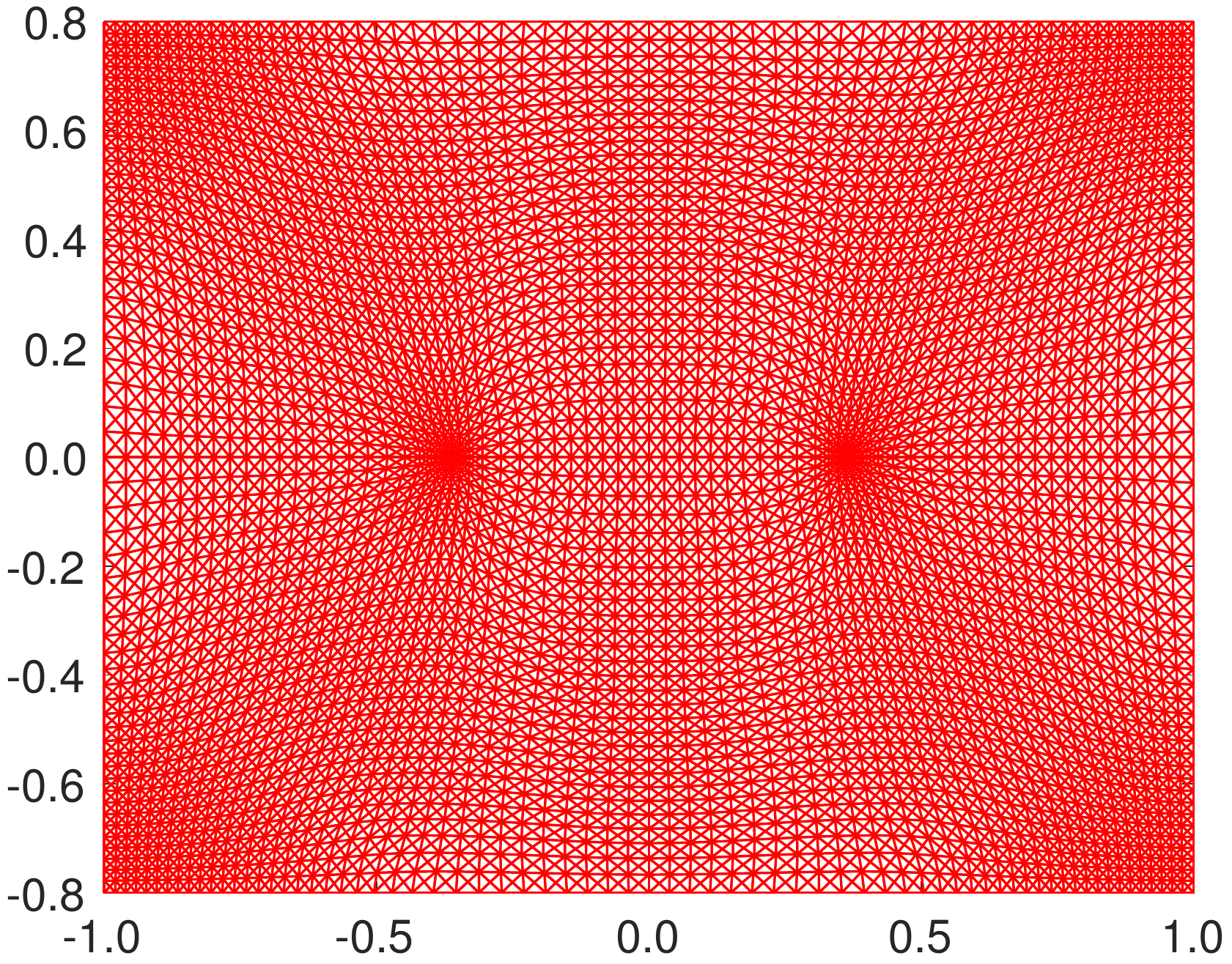} \label{nh1b}}\quad
\subfigure[Mesh at touchdown, $\eps= 0.1$.]{\includegraphics[width = 0.3\textwidth]{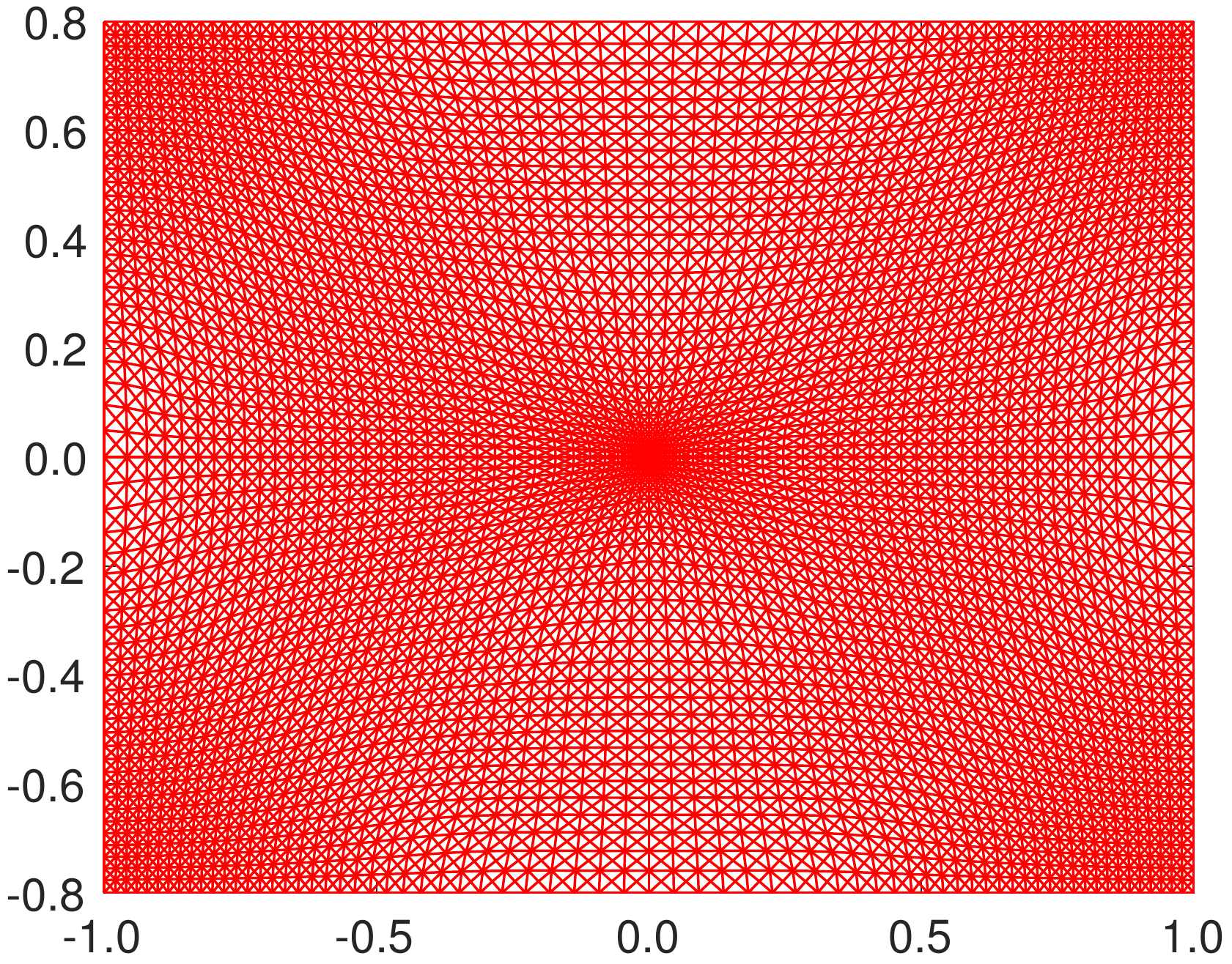} \label{nh1c}}\\
\subfigure[Solution at touchdown, $\eps= 0.02$.]{\includegraphics[width = 0.3\textwidth]{rectu1.pdf} \label{nh1d}}\quad
\subfigure[Solution at touchdown, $\eps= 0.068$.]{\includegraphics[width = 0.3\textwidth]{rectu2.pdf} \label{nh1e}}\quad
\subfigure[Solution at touchdown, $\eps= 0.1$.]{\includegraphics[width = 0.3\textwidth]{rectu3.pdf} \label{nh1f}}
\caption{The profiles $u(\bx,t)$ of \eqref{eq:intro} and associated meshes very close to singularity for $\eps = 0.02$, $\eps = 0.068$ and $\eps = 0.1$ in the rectangle $(-1,1)\times(-0.8,0.8)$. The mesh size is $N = 15680\;  (70 \times 56)$. \label{nh1}}
\end{figure}

 \begin{figure}[htbp]
\centering
\subfigure[Mesh at touchdown, $\eps= 0.02$. Mesh size $N = 6240\; (40 \times 39)$.]{\includegraphics[width = 0.3\textwidth]{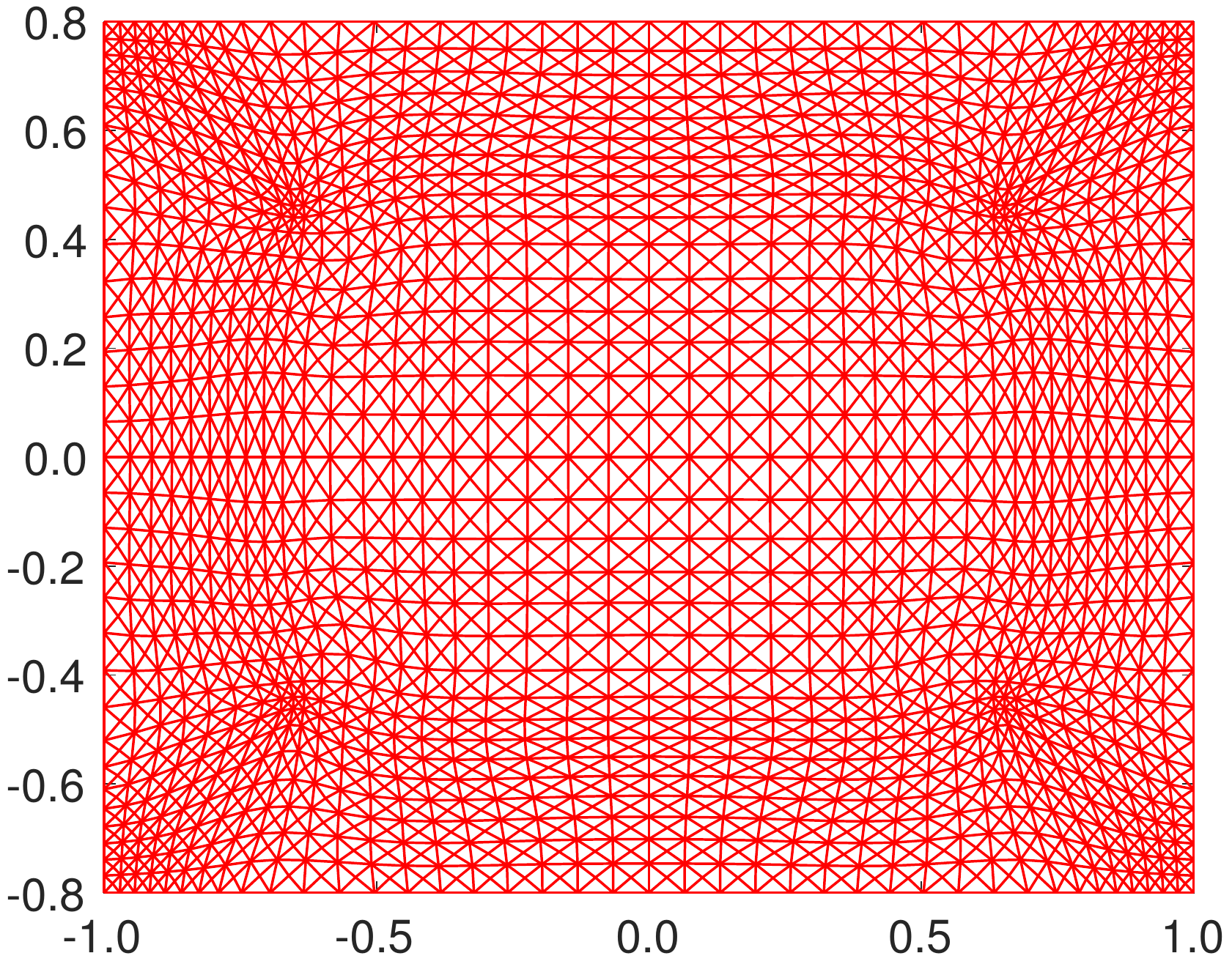} \label{nhcom1a}}\quad
\subfigure[Mesh at touchdown, $\eps= 0.068$. Mesh size $N = 6240\; (40 \times 39)$.]{\includegraphics[width = 0.3\textwidth]{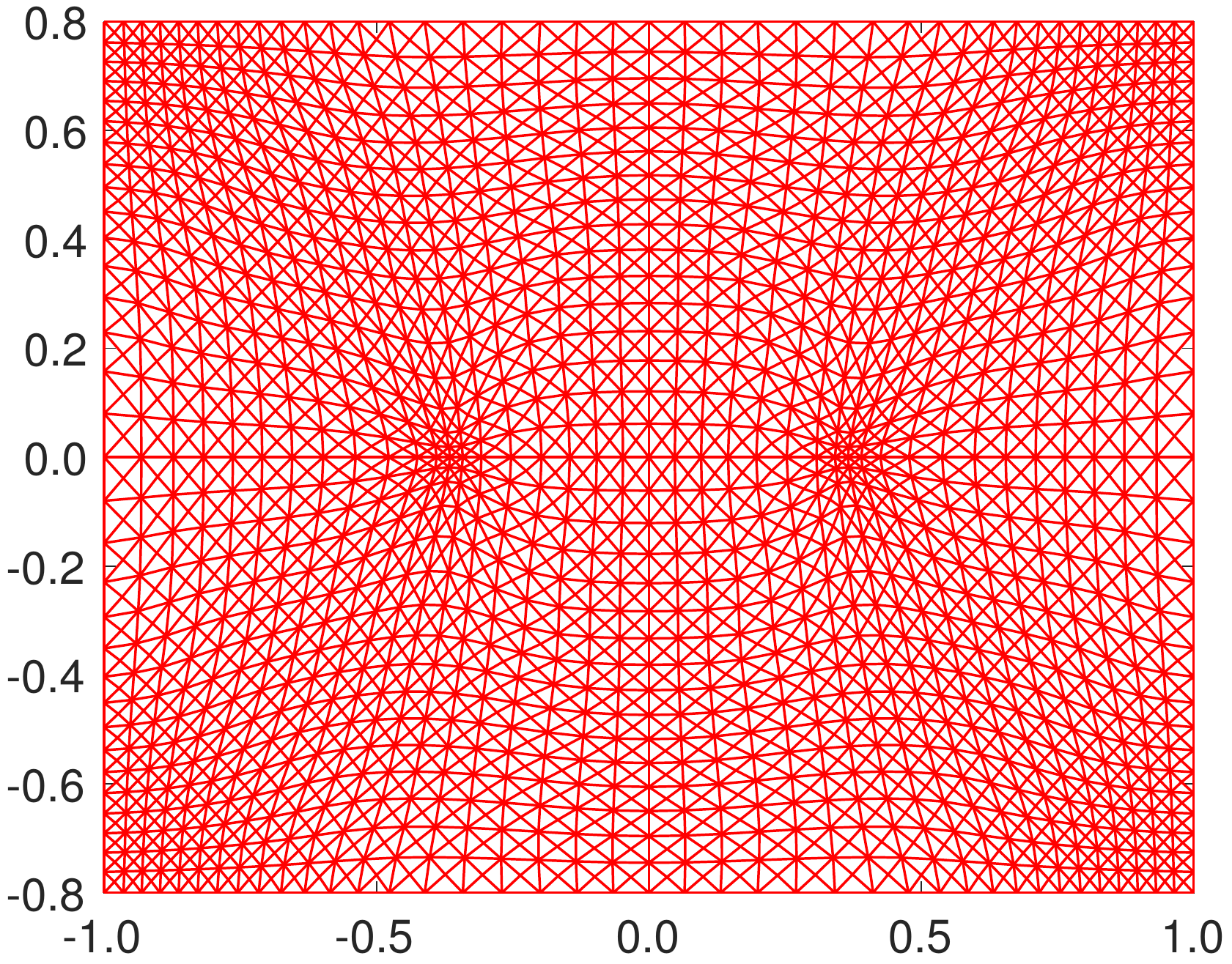} \label{nhcom1b}}\quad
\subfigure[Mesh at touchdown, $\eps= 0.1$. Mesh size $N = 6240\; (40 \times 39)$.]{\includegraphics[width = 0.3\textwidth]{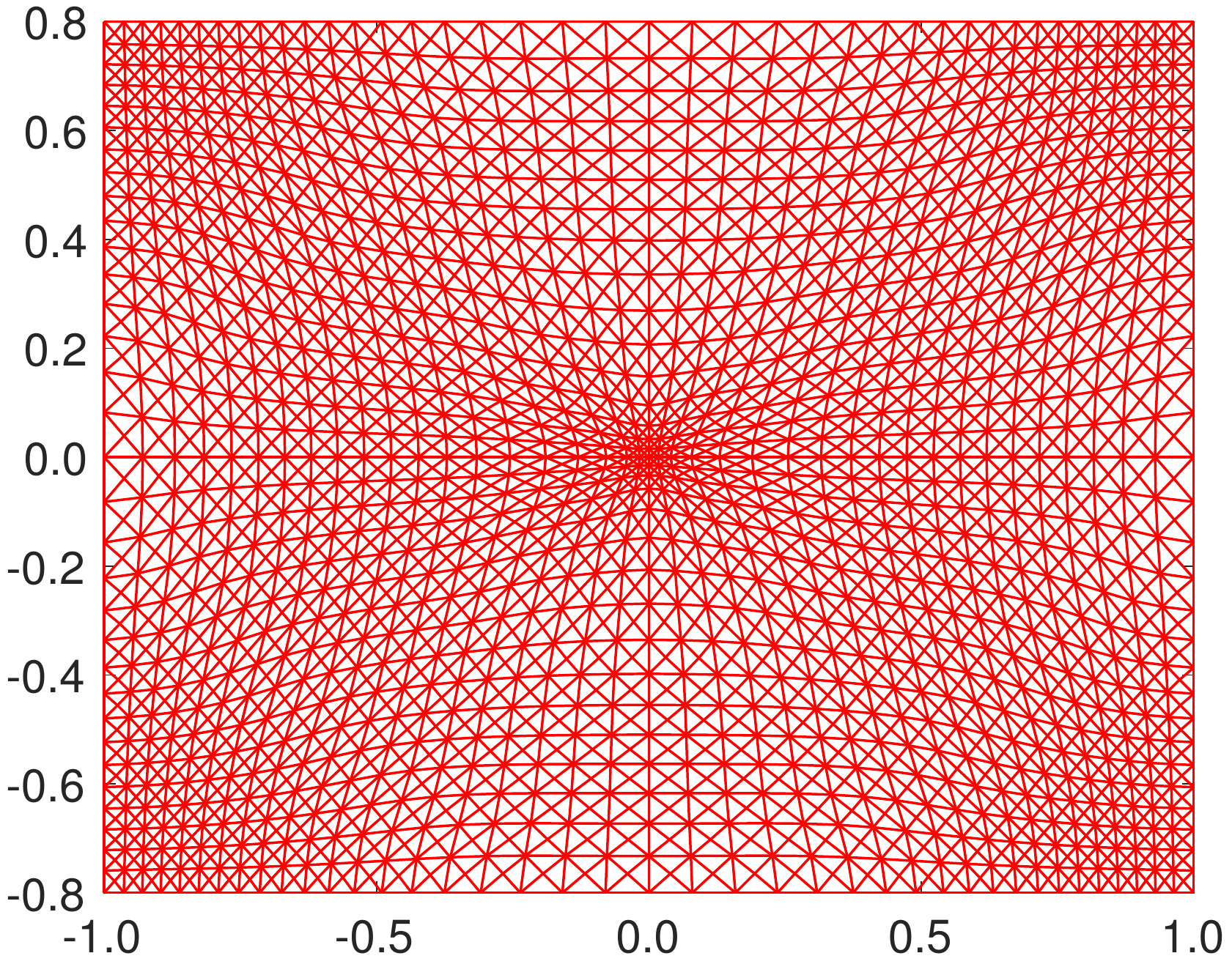} \label{nhcom1c}}\\
\subfigure[Mesh at touchdown, $\eps= 0.02$. Mesh size $N = 15680\; (70 \times 56)$.]{\includegraphics[width = 0.3\textwidth]{rectMesh1.pdf} \label{nhcom1a14}}\quad
\subfigure[Mesh at touchdown, $\eps= 0.068$. Mesh size $N = 15680\; (70 \times 56)$.]{\includegraphics[width = 0.3\textwidth]{rectMesh2.pdf} \label{nhcom1b14}}\quad
\subfigure[Mesh at touchdown, $\eps= 0.1$. Mesh size $N = 15680\; (70 \times 56)$.]{\includegraphics[width = 0.3\textwidth]{rectMesh3.pdf} \label{nhcom1c14}}\\
\caption{The profiles $u(\bx,t)$ of \eqref{eq:intro} and associated meshes very close to singularity for $\eps = 0.02$, $\eps = 0.068$ and $\eps = 0.1$ in the rectangle $(-1,1)\times(-0.8,0.8)$. The top row is obtained by mesh size $N = 6240\; (40 \times 39)$, and the bottom row is obtained by mesh size $N = 15680\; (70 \times 56)$. \label{nhcom1}}
\end{figure}

 As predicted by the analytical skeleton, there are four singularities close to each corner for small $\eps$.
 As $\eps$ increases, the four singularities move inwards along $\skel$ merging first into two singularities and, as $\eps$ increases, eventually into one.
 The final mesh and solution for values $\eps = 0.02$, $0.068$ and $0.1$
 are shown in Fig.~\ref{nh1}. In each case shown in Figs.~\ref{nh1a}-\ref{nh1c}, the numerical algorithm correctly locates the position and multiplicity of the forming singularities and increases local mesh density in their vicinity to accurately resolve the solution. To demonstrate that the solution is robust with respect to grid refinement, we present the final mesh for $\eps = 0.02$, $0.068$ and $0.1$ obtained with mesh size $N = 6240\; (40 \times 39)$ and $N =15680\; (70 \times 56)$ in Fig.~\ref{nhcom1}.

In comparing $\skel$ with the numerical touchdown points, we see that at smaller values of $\eps$  (for which the singularities are confined to the corners) the set $\skel$ accurately predicts the touchdown set. At larger values of $\eps$, in particular those at which the four singularities merge into two, we observe that $\skel$ has reduced accuracy in predicting the contact set. This reduction in the quality of the prediction is not surprising considering the asymptotic formulation relies on the peaks being well separated at touchdown which is not valid for larger $\eps$ values. Nevertheless, the skeleton theory gives a qualitatively accurate description of the possible touchdown locations and multiplicities.

\begin{figure}[htbp]
\centering
\subfigure[Mesh at $t = 0.002$.]{\includegraphics[width = 0.295\textwidth]{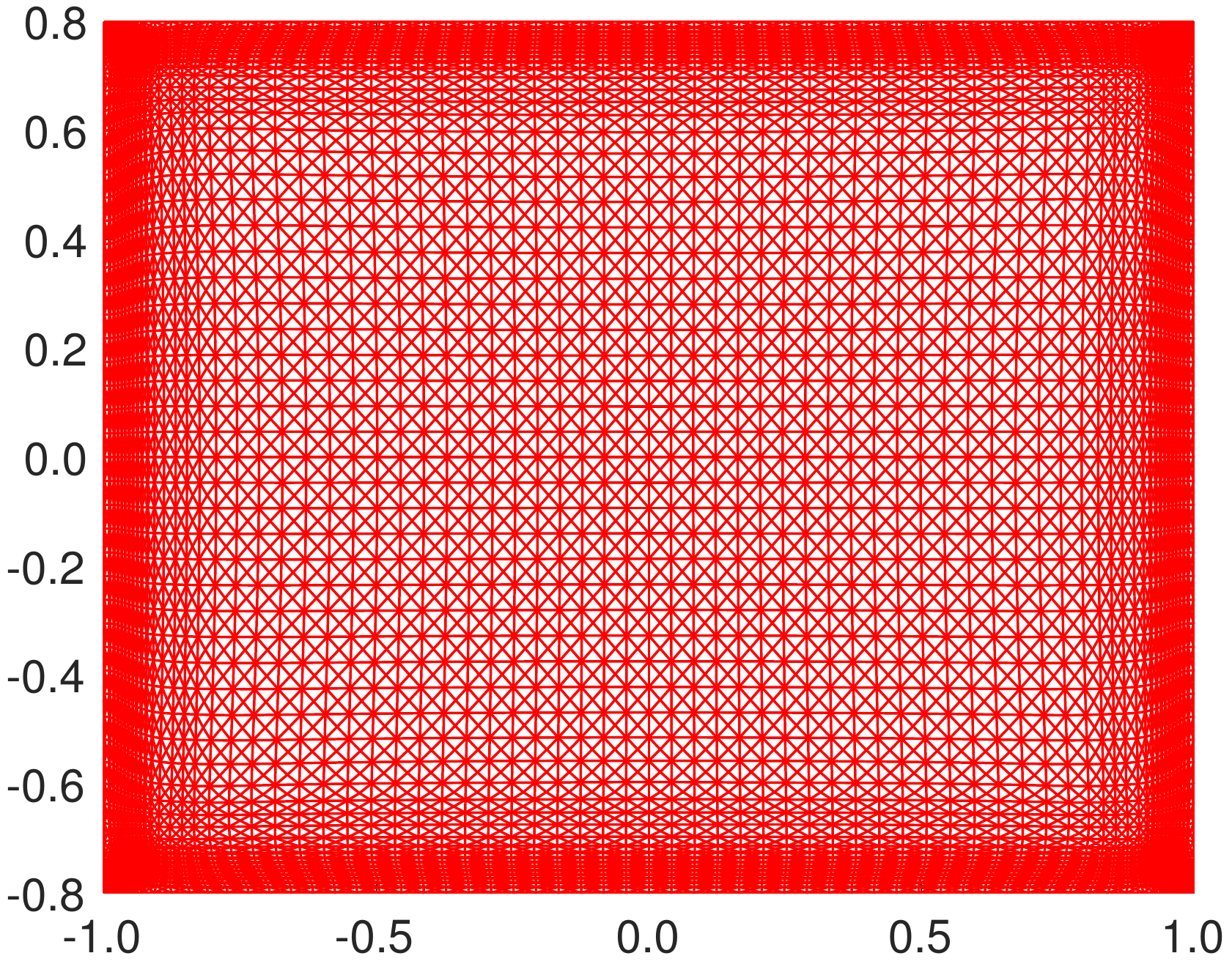}\label{fig:RectEvolution_a} }\qquad
\subfigure[Mesh at $t = 0.163$.]{\includegraphics[width = 0.295\textwidth]{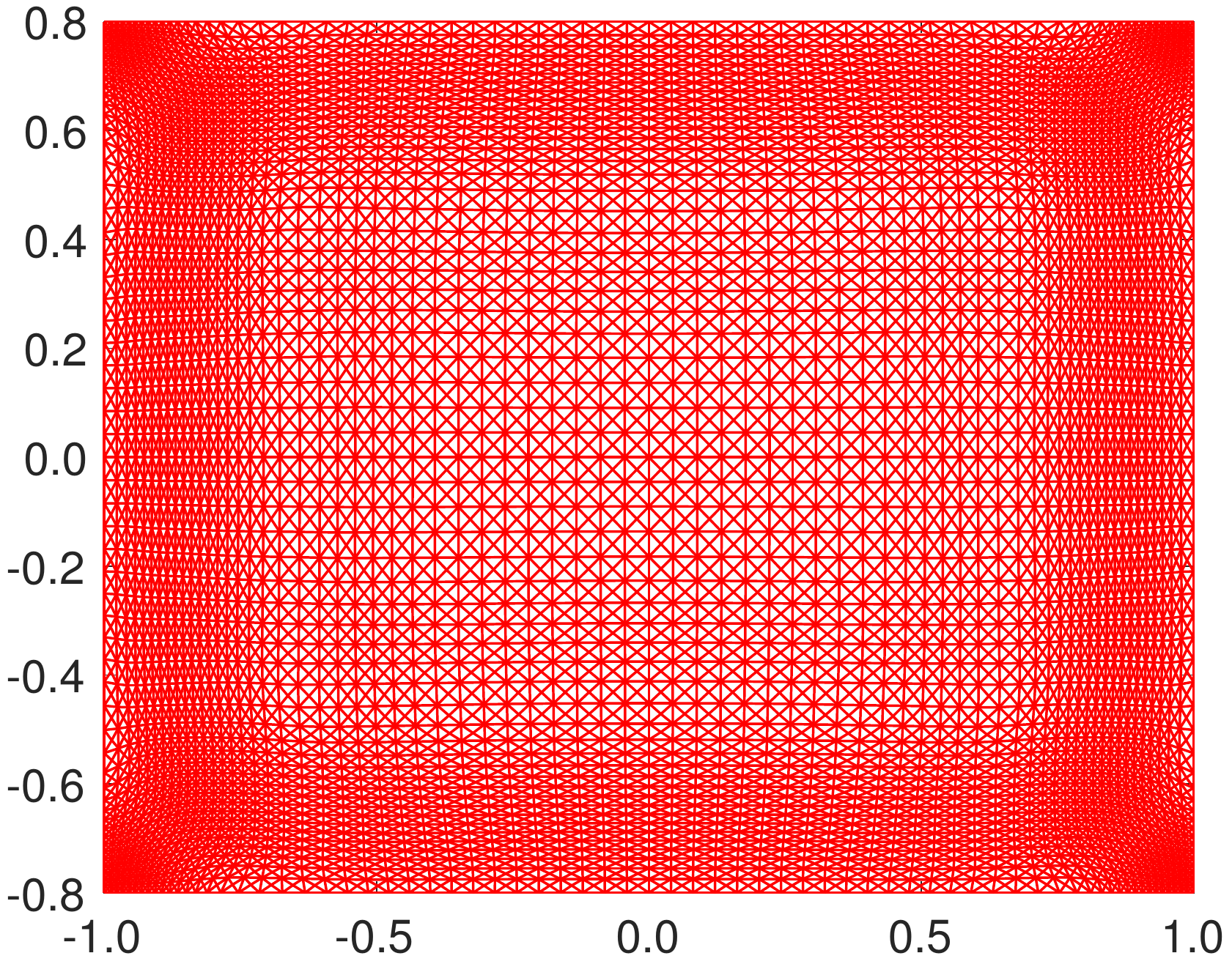}\label{fig:RectEvolution_b} }\qquad
\subfigure[Mesh at $t = 0.312$.]{\includegraphics[width = 0.295\textwidth]{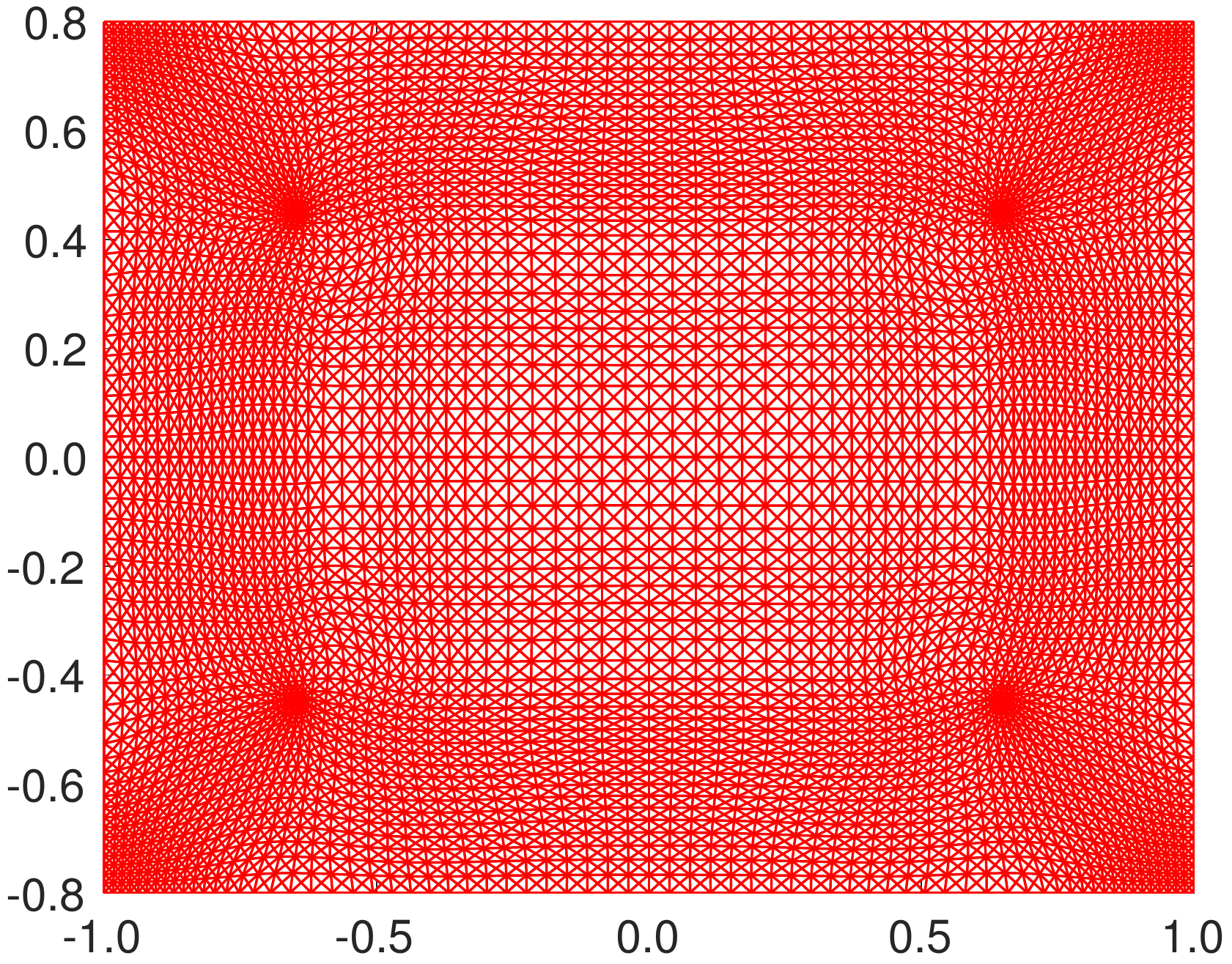}\label{fig:RectEvolution_c} }\\
\subfigure[Solution at $t = 0.003 $.]{\includegraphics[width = 0.295\textwidth]{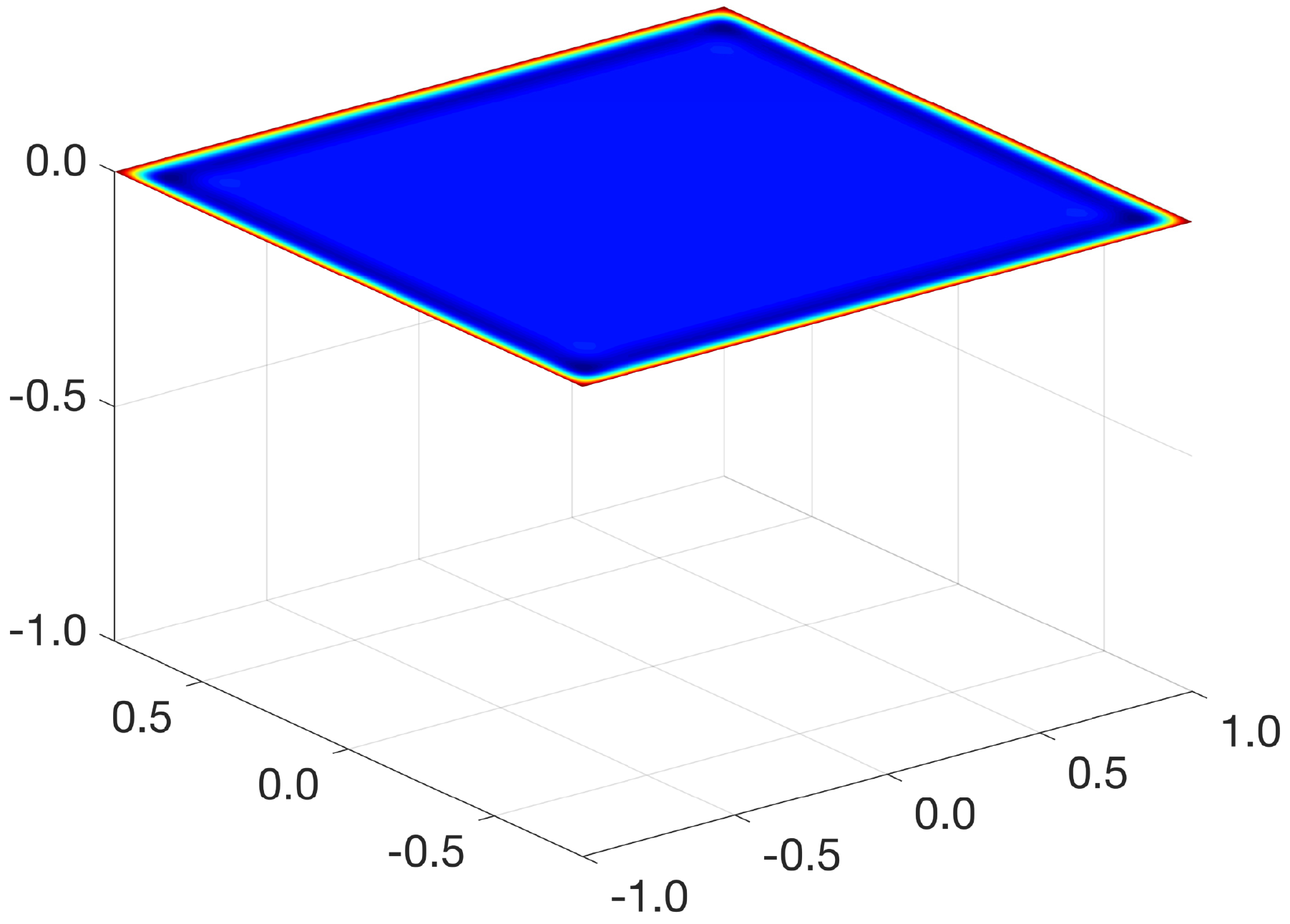}\label{fig:RectEvolution_d} }\qquad
\subfigure[Solution at $t = 0.163$.]{\includegraphics[width = 0.295\textwidth]{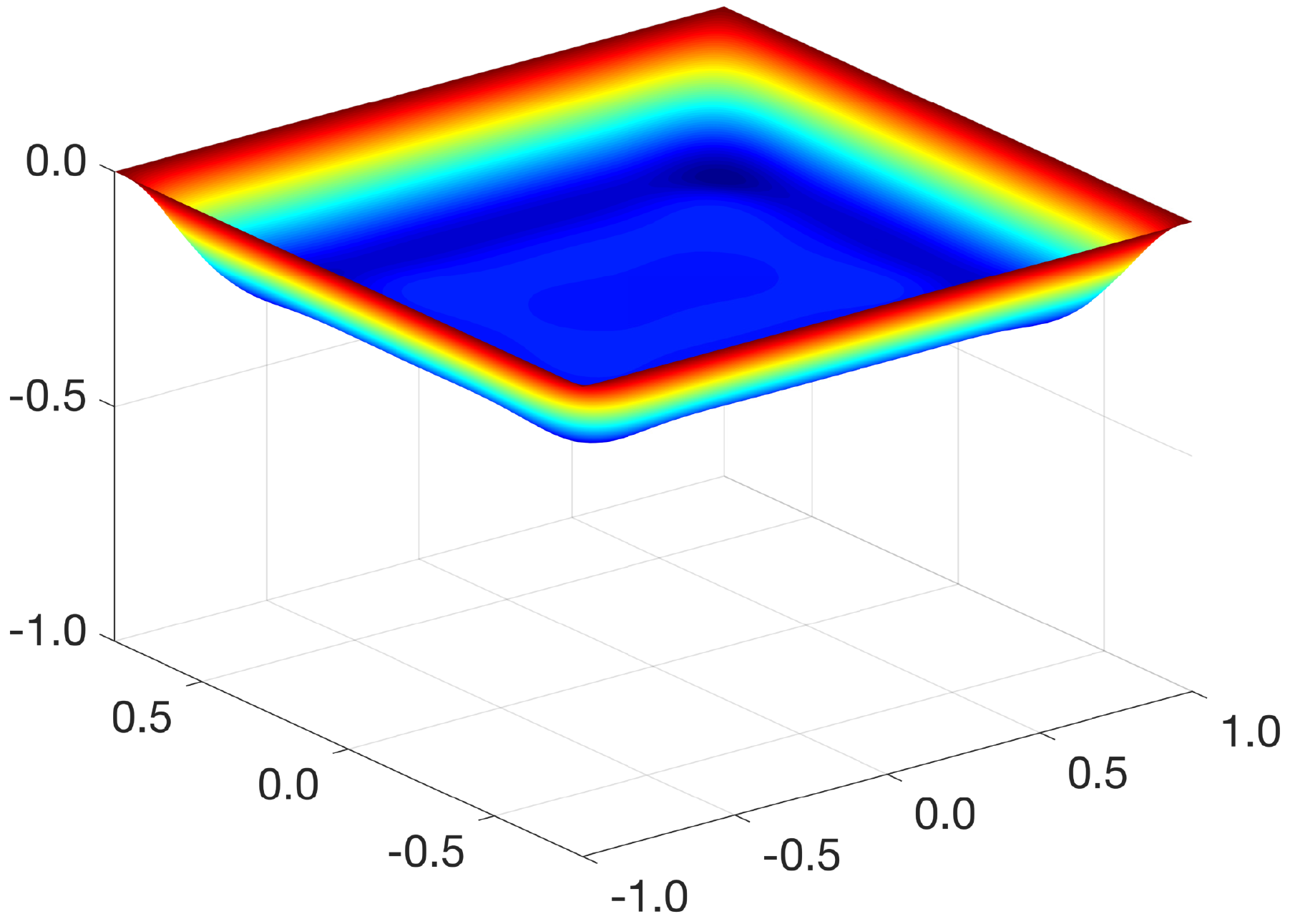}\label{fig:RectEvolution_e} }\qquad
\subfigure[Solution at $t = 0.312$.]{\includegraphics[width = 0.295\textwidth]{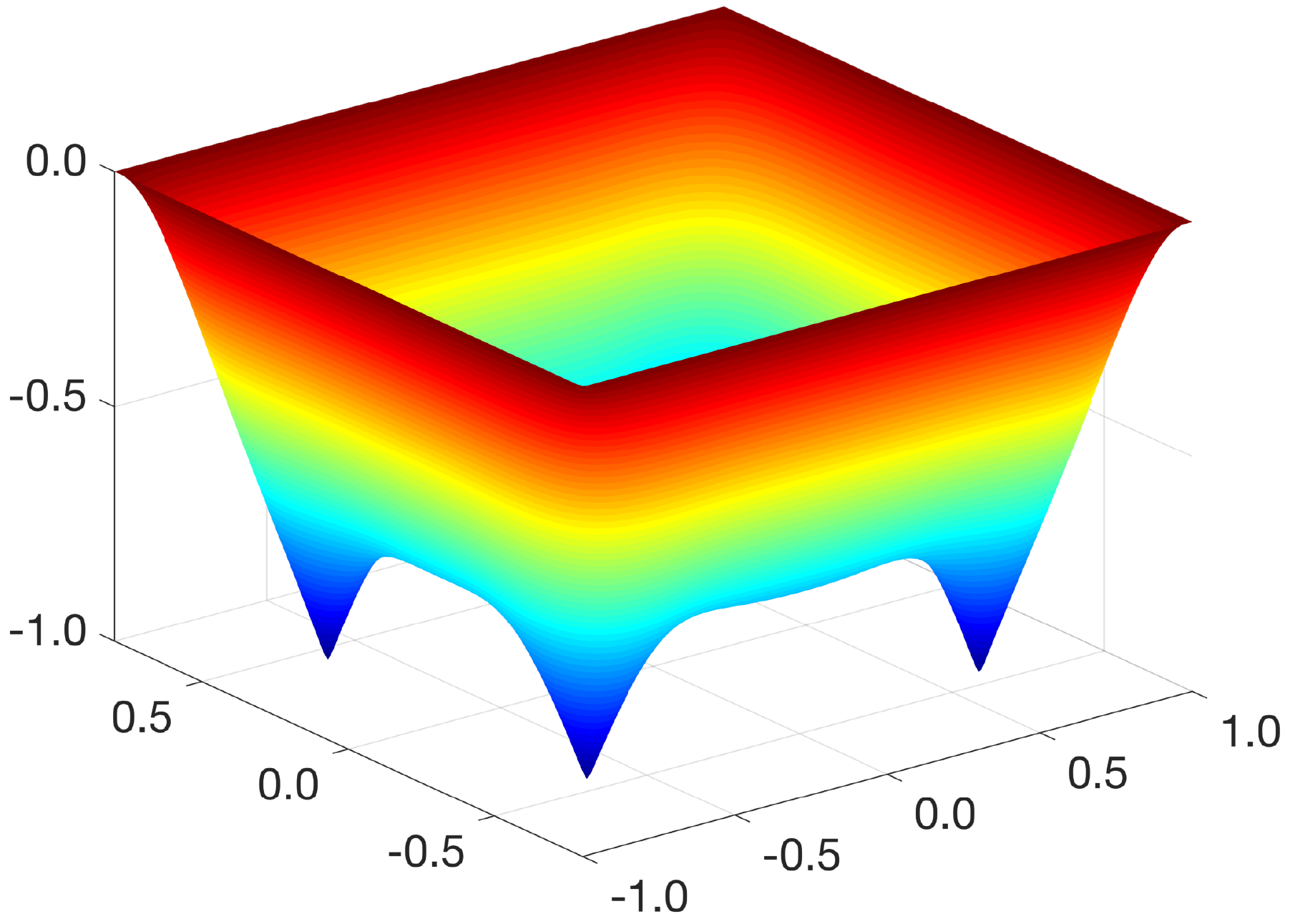}\label{fig:RectEvolution_f} }
\caption{Evolution of the solution of \eqref{eq:intro} and the associated mesh for $\eps = 0.02$ in the rectangular domain for three time instants. \label{fig:RectEvolution}}
\end{figure}

In Fig.~\ref{fig:RectEvolution} we display the evolution of the solution to \eqref{eq:intro} for the fixed value $\eps=0.02$ and several temporal snapshots with the accompanying mesh. For this value of $\eps$, touchdown is observed at four points simultaneously. At short times (Fig.~\ref{fig:RectEvolution_a}), the computational mesh is adapted to the propagating boundary layers emanating around $\partial\Omega$. By the touchdown time (Fig.~\ref{fig:RectEvolution_c}), the mesh generation algorithm allocates resources between each of the four forming singularities and the sharp ridges that join them.

%\begin{figure}[htbp]
%\centering
%\subfigure[]{\includegraphics[width =0.3\textwidth]{nh2u1.eps}\label{nh2a}}\qquad
%\subfigure[]{\includegraphics[width =0.3\textwidth]{nh2mesh1.eps}\label{nh2b}}\\
%\subfigure[]{\includegraphics[width =0.3\textwidth]{nh2u2.eps}\label{nh2c}}\qquad
%\subfigure[]{\includegraphics[width =0.3\textwidth]{nh2mesh2.eps}\label{nh2d}}\\
%\subfigure[]{\includegraphics[width =0.3\textwidth]{nh2u3.eps}\label{nh2e}}\qquad
%\subfigure[]{\includegraphics[width =0.3\textwidth]{nh2mesh3.eps}\label{nh2f}}
%\caption{The evolution of the solution for $\eps =0.068$. The mesh size is $N = 7644\;  (39 \times 49)$. \label{nh2}}
%\end{figure}
%%
%\begin{figure}[htbp]
%\centering
%\subfigure[]{\includegraphics[width =0.45\textwidth]{nh3u1.eps}\label{nh3a}}\qquad
%\subfigure[]{\includegraphics[width =0.45\textwidth]{nh3mesh1.eps}\label{nh3b}}\\
%\subfigure[]{\includegraphics[width =0.45\textwidth]{nh3u2.eps}\label{nh3c}}\qquad
%\subfigure[]{\includegraphics[width =0.45\textwidth]{nh3mesh2.eps}\label{nh3d}}\\
%\subfigure[]{\includegraphics[width =0.45\textwidth]{nh3u3.eps}\label{nh3e}}\qquad
%\subfigure[]{\includegraphics[width =0.45\textwidth]{nh3mesh3.eps}\label{nh3f}}
%\caption{The evolution of the solution for $\eps =0.1$. The mesh size is $N = 7644\;  (39 \times 49)$.\label{nh3}}
%\end{figure}

\vspace{10pt}

%% example 2
 \begin{exam}[Rectangular domain with a hole]
{\em Here we consider the rectangular domain $\Omega = (-1,1)\times (-0.8,0.8)$, with a circular hole of radius $0.2$, centered at $(0.2,0.3)$. In this example $\Omega$ is non-convex.}
\label{examrec}
\end{exam}

For this example, we have found that it is important to keep a level of mesh concentration
around the hole. To this end, we modify the metric tensor as 
\begin{equation} 
\tilde{\M}_K = \M_K + \beta I,
\label{M-2}
\end{equation}
where $\M_K$ is defined as in (\ref{M-1}) and $\beta$ is chosen as 
\[
\beta =  \left [ e^{4(0.2-\sqrt{(x - 0.2)^2 +(y-0.3)^2}\,)} - 1
+ \frac{2}{\max\limits_{K \in \mathcal{T}_h} \sqrt{\det(\M_K)} } \right ]^{-1}.
\]  
Notice that for $(x,y)$ on the circle, this gives
\[ 
\tilde{\M}_K = \M_K + \frac{1}{2} \max_{K \in \mathcal{T}_h}\sqrt{\det(\M_K)}\; I,
\]
which will give a level of mesh concentration around the circle comparable to that in the regions
with largest $\sqrt{\det(\M_K)}$. The exponential term makes $\beta$ decrease sufficiently fast such that
the mesh elements are not over concentrated near the circle.

\begin{figure}[htbp]
\centering
%\subfigure[]{\includegraphics[width = 0.3\textwidth]{recholesk.pdf}}\qquad
%\subfigure[Skeleton of the Domain]{\includegraphics[width = 0.3\textwidth]{rectangleSkeletonHole.pdf}\label{fig:rectHolSkel}}
\includegraphics[width = 0.4\textwidth]{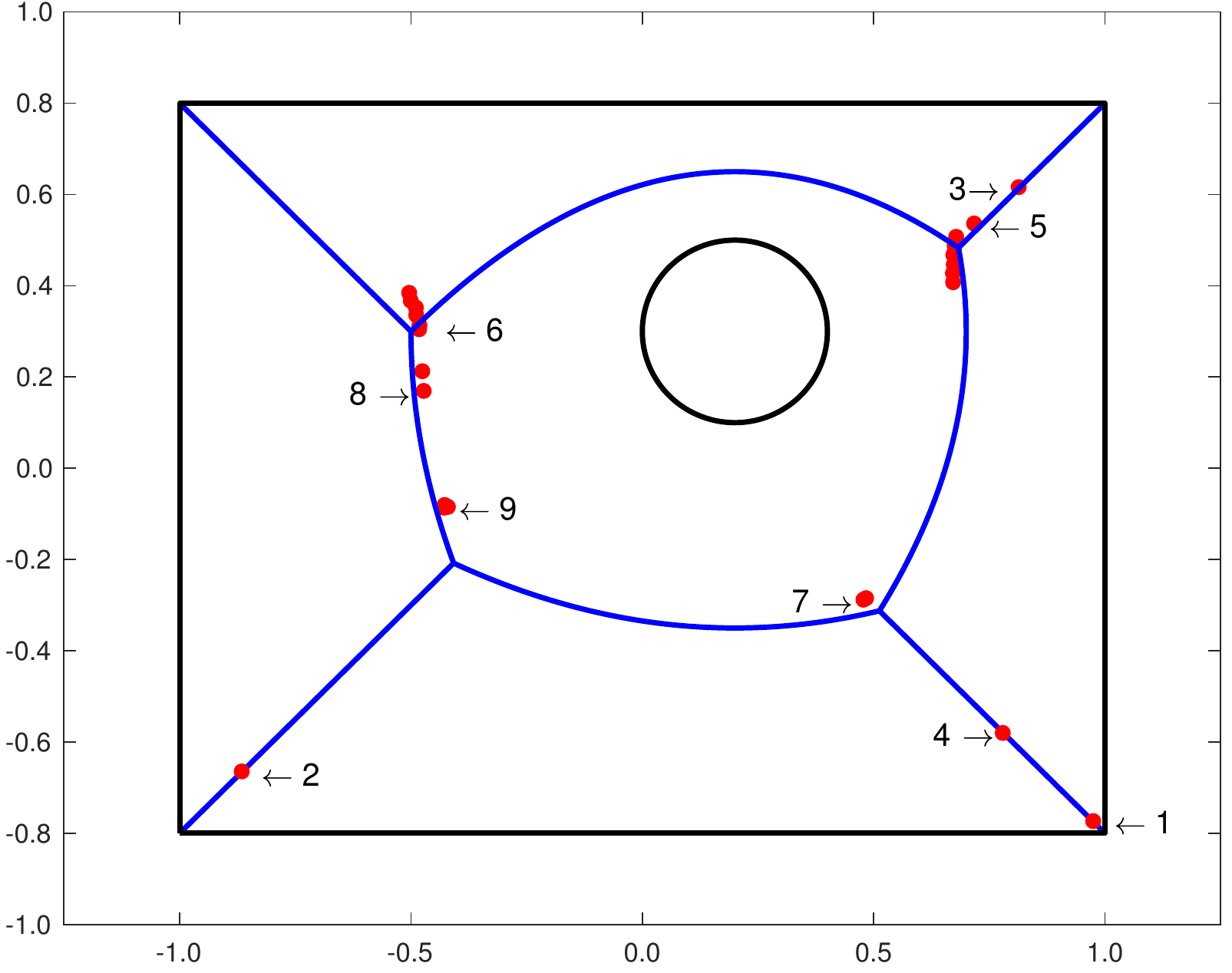}
\caption{Skeleton of rectangular domain hole (blue solid line) and numerically computed touchdown locations (red dots).
The points marked $1-9$ correspond to the first touchdown location for solutions of \eqref{eq:intro} for values $\eps = 10^{-4}$, $2.662 \times 10^{-3}$, $5.2 \times 10^{-3}$, $7.78 \times 10^{-3}$, $0.01$, $0.036$, $0.044$, $0.051$, and $0.06$, respectively.
The solution and mesh for Mark 5 ($\eps = 0.01$), Mark 7 ($\eps = 0.044$), and
Mark 9 ($\eps = 0.06$) are shown in Fig.~\ref{fig:RectHole}.
\label{recsk}}
\end{figure}

The skeleton $\skel$ of the domain which is displayed in Fig.~\ref{recsk}. In this example $\skel$ is formed from straight line segments that originate from each corner and are linked by four curved segments contorted around the hole. The expressions for the parabolic segments of $\skel$ are found analytically by considering the points that are equidistant from the boundary of the outer rectangle and the perturbing hole. In this example, $\skel$ intersects with the boundary hence the skeleton arrival time is $\sat=0$.

\begin{figure}[!ht]
\centering
\subfigure[Mesh at touchdown, $\eps= 0.01$.]{\includegraphics[width = 0.295\textwidth]{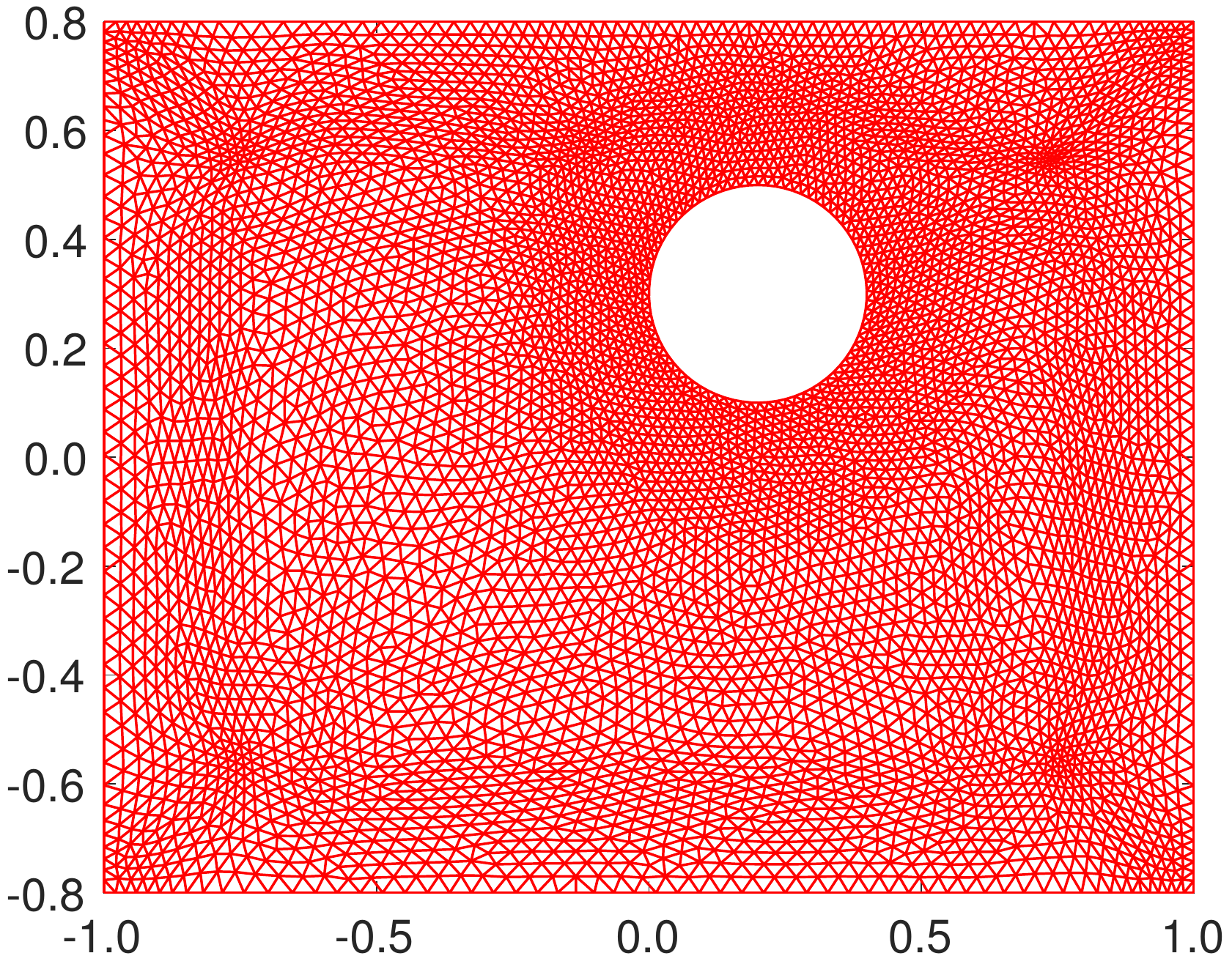}\label{fig:RectHole_a} }\qquad
\subfigure[Mesh at touchdown, $\eps= 0.044$.]{\includegraphics[width = 0.295\textwidth]{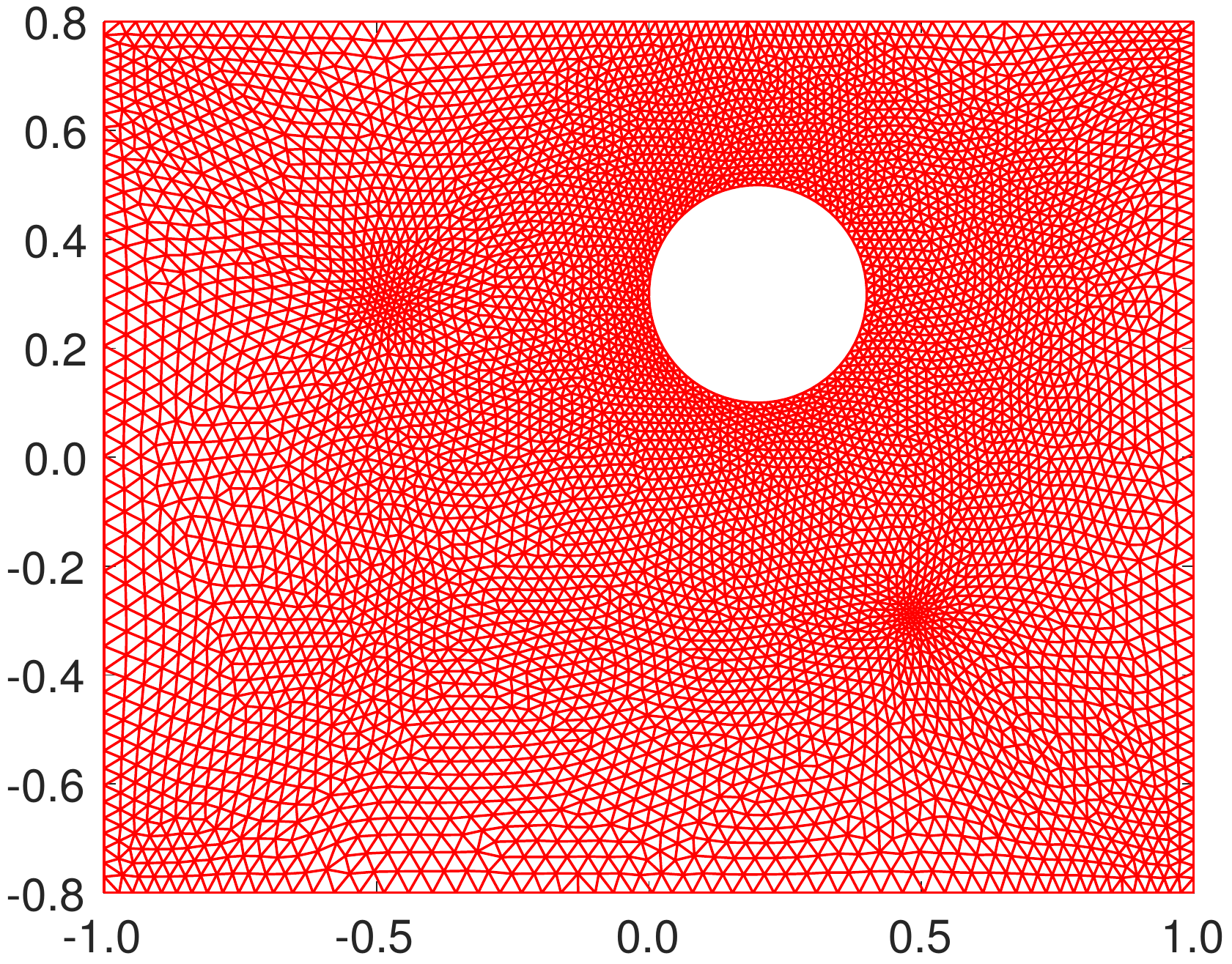}\label{fig:RectHole_b} }\qquad
\subfigure[Mesh at touchdown, $\eps= 0.06$.]{\includegraphics[width = 0.295\textwidth]{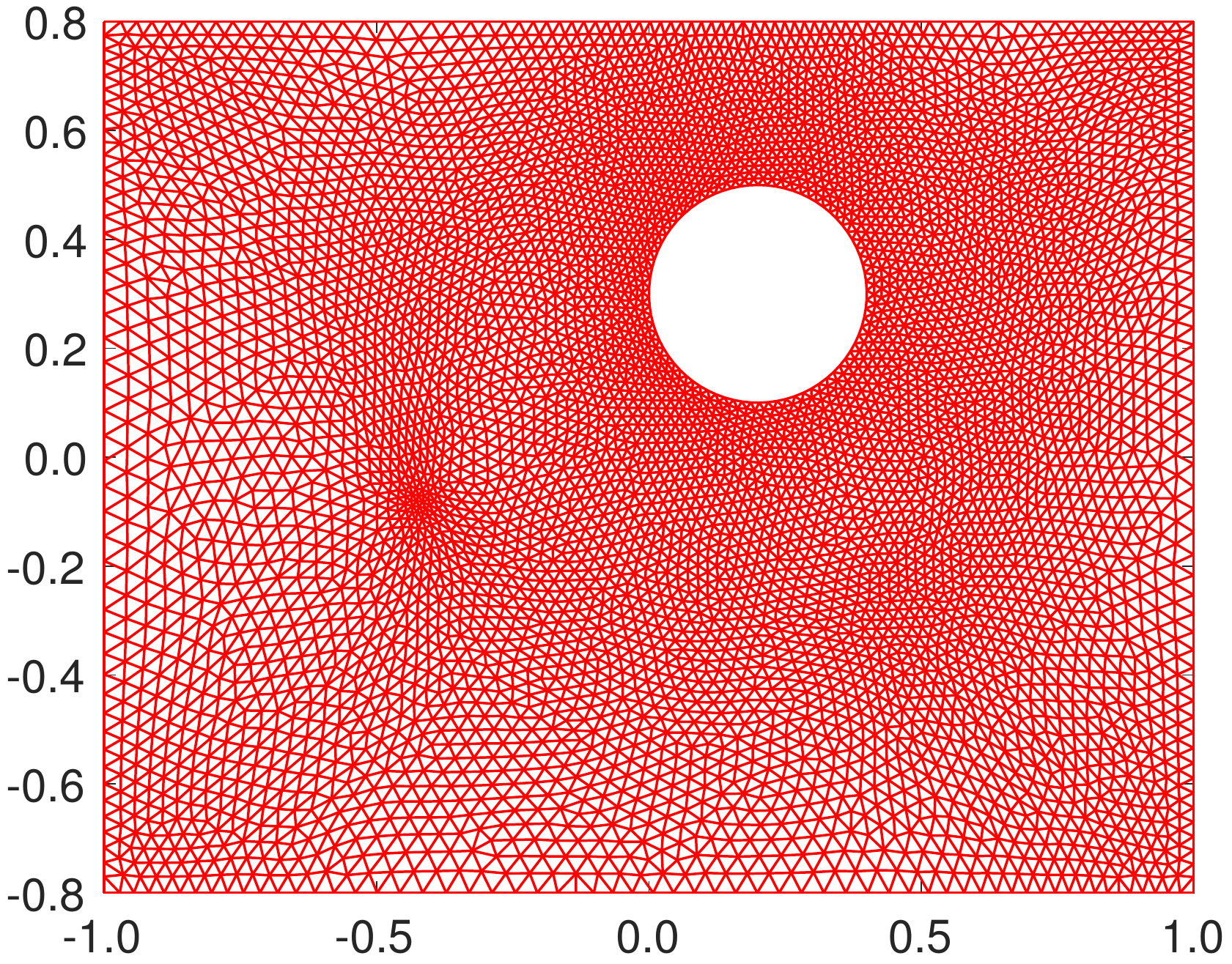}\label{fig:RectHole_c} }\\
\subfigure[Solution at touchdown, $\eps= 0.01$.]{\includegraphics[width = 0.295\textwidth]{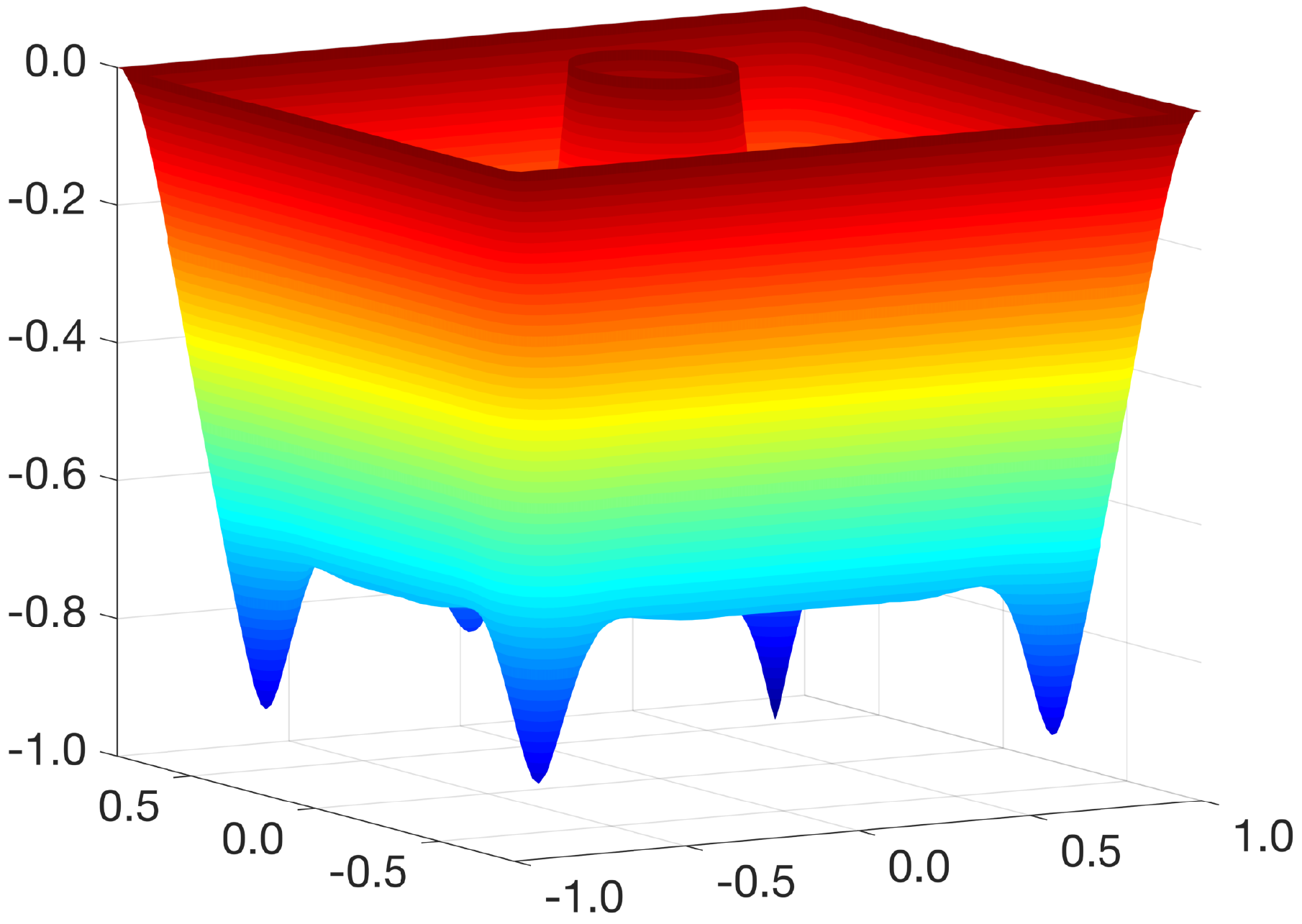}\label{fig:RectHole_d} }\qquad
\subfigure[Solution at touchdown, $\eps= 0.044$.]{\includegraphics[width = 0.295\textwidth]{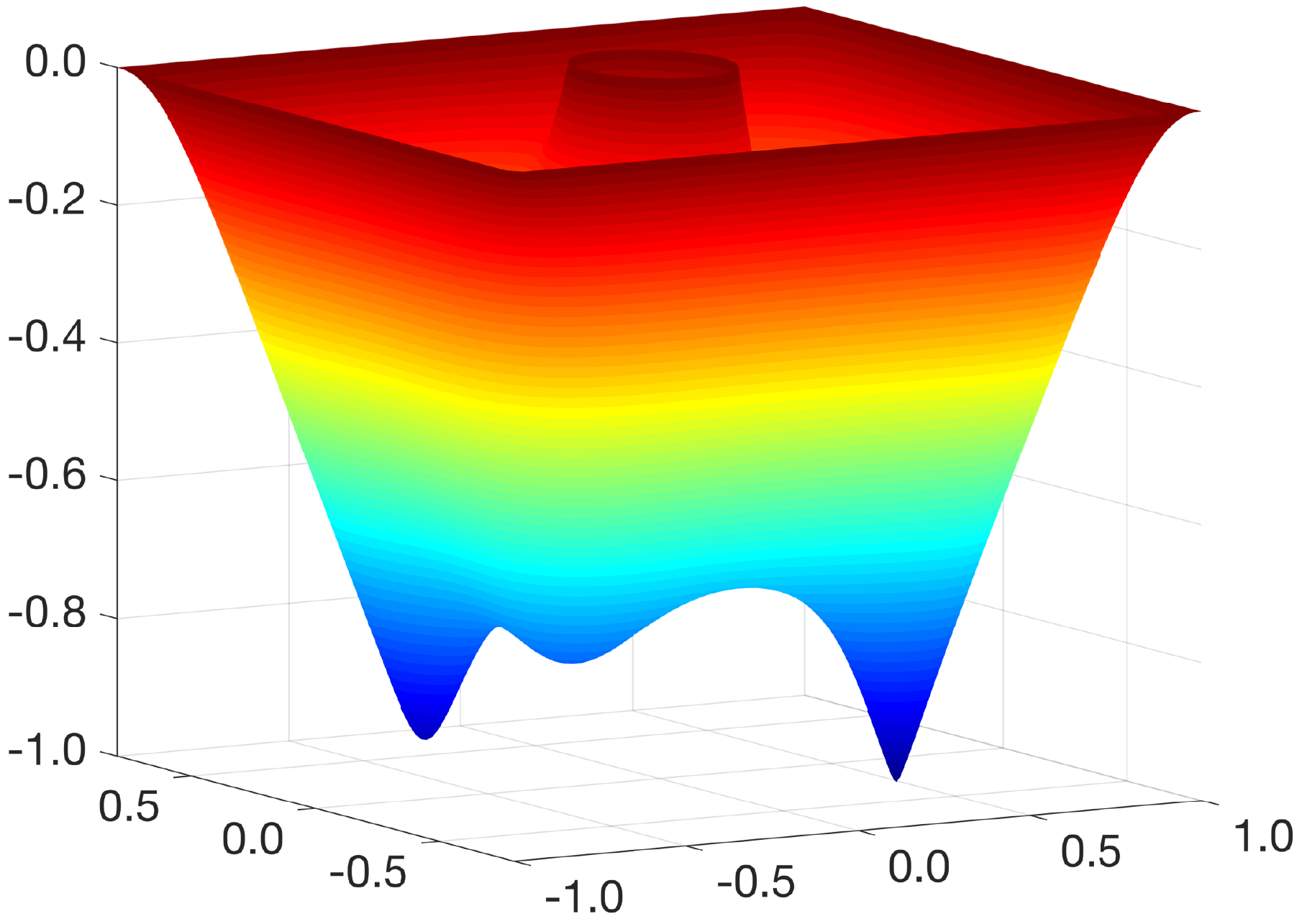}\label{fig:RectHole_e} }\qquad
\subfigure[Solution at touchdown, $\eps= 0.06$.]{\includegraphics[width = 0.295\textwidth]{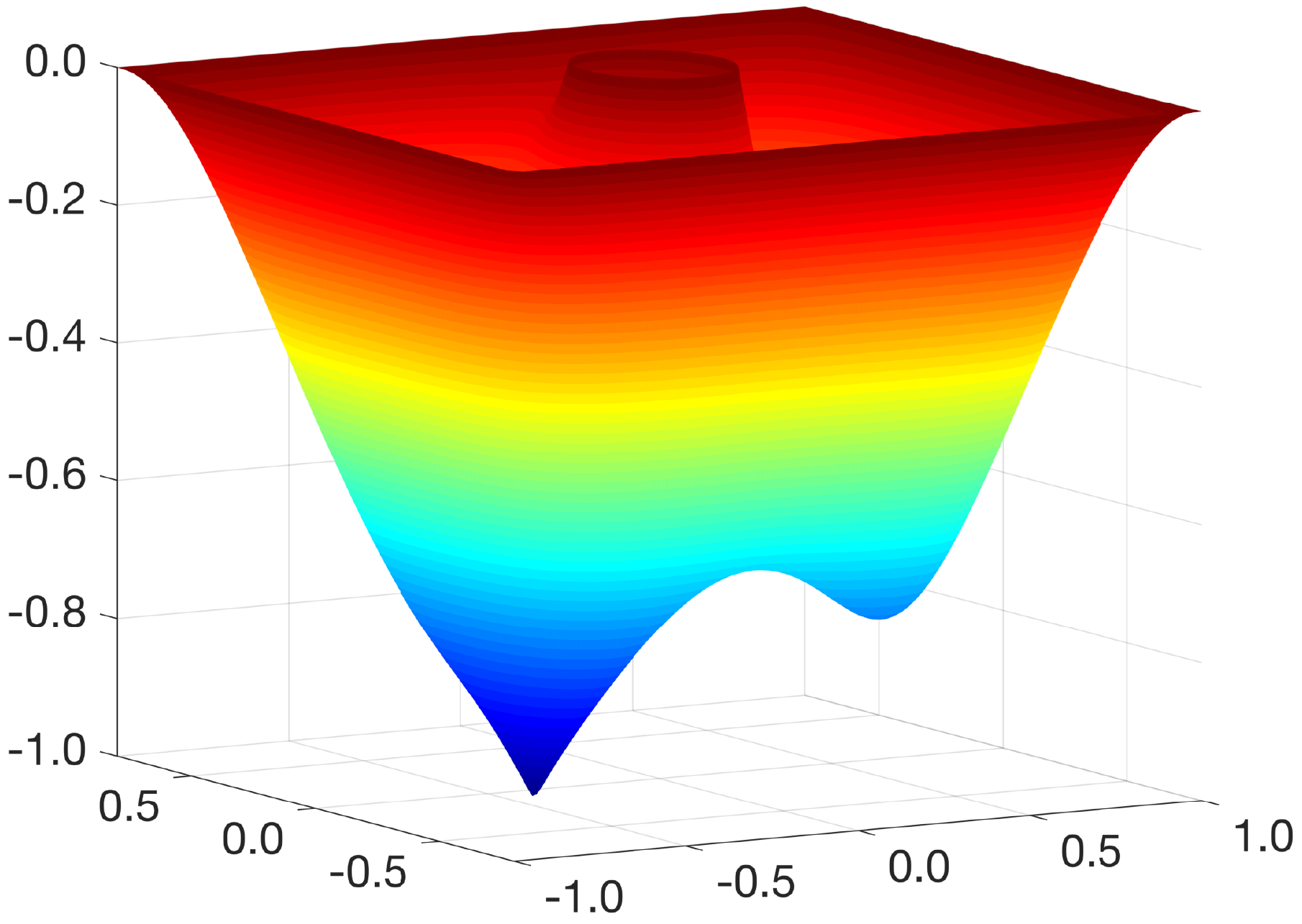}\label{fig:RectHole_f} }
\caption{Solutions of \eqref{eq:intro} and meshes at singularity for values $\eps = 0.01, 0.044, 0.1$ in the rectangular domain with hole. \label{fig:RectHole}}
\end{figure}

The presence of the hole breaks the symmetry of the domain. In simulations, we observe that this precludes touchdown at multiple points simultaneously, except for certain fixed values of $\eps$. Simultaneous touchdown at multiple points, as observed in the previous example of the rectangle with no hole, relies on the symmetric properties of the domain. In the absence of such symmetries, single point touchdown is the expectation for solutions of equation \eqref{eq:intro}. However, as is clear from the solution profile of \eqref{eq:intro} for $\eps=0.044$ shown in Fig.~\ref{fig:RectHole_e}, the solution may be forming multiple troughs.  When there are multiple troughs present, the singularity location will be selected by the trough which has the lowest value as the singularity is approached.  In terms of applications to MEMS, each of these troughs can form contacts between the two surfaces and are important to track. We remark that $\skel$ describes a set of potential points at which the asymptotic solution has a local minimum and therefore considers all potential contact locations.

%In Fig. \ref{recsk} we present the skeleton for $\eps \in linspace(0.02,0.1,21)$ as well as a few other values including
%\[
%\eps = \{ 10^{-4},2\times10^{-4},5 \times 10^{-4},10^{-3} ,3 \times 10^{-3}, 3\times 10^{-3} ,5 \times 10^{-3} , \]
%\[
%0.042 , 0.046, 0.15 , 0.2, 0.25 \}.
%\]
%For all $\eps$ values above, the singularities happen at one point, and as $\eps$ increases,
%the skeleton is tracked by the numbers marked in the figure, with ascending order in $\eps$ value.
%For small $\eps$ values, touchdown behavior occurs at four points which are close to the four corners,
%yet singularity only happens at one point (see Fig. \ref{5e-4}). For $\eps = 10^{-4}$ and $2 \times 10^{-4}$,
%$u$ touches down at Mark $1$. As $\eps$ increases ($10^{-3}$, $3 \times 10^{-3}$, $5 \times 10^{-3}$),
%the singularities travel around at Mark $2$, $3$, $4$, $5$, respectively. Then for $\eps = 5 \times 10^{-3}$ and $0.01$, the singularity point is at Mark $5$ and $6$. Mark $7$ to $11$ are the singularity points for $\eps \in \text{linspace}(0.02,0.1,21)$. For $\eps > 0.1$, the singularities just concentrate around Mark $12$.

%The evolution at Mark $2$ ($\eps = 5 \times10^{-4}$), Mark $6$ ($\eps = 0.01$) and
%Mark $12$ ($\eps = 0.2$) have been presented in Fig. \ref{5e-4}, Fig. \ref{r2e-2}, and Fig. \ref{r1e-1}.

\begin{figure}[!ht]
\centering
\subfigure[Mesh at $t = 0.035$.]{\includegraphics[width = 0.3\textwidth]{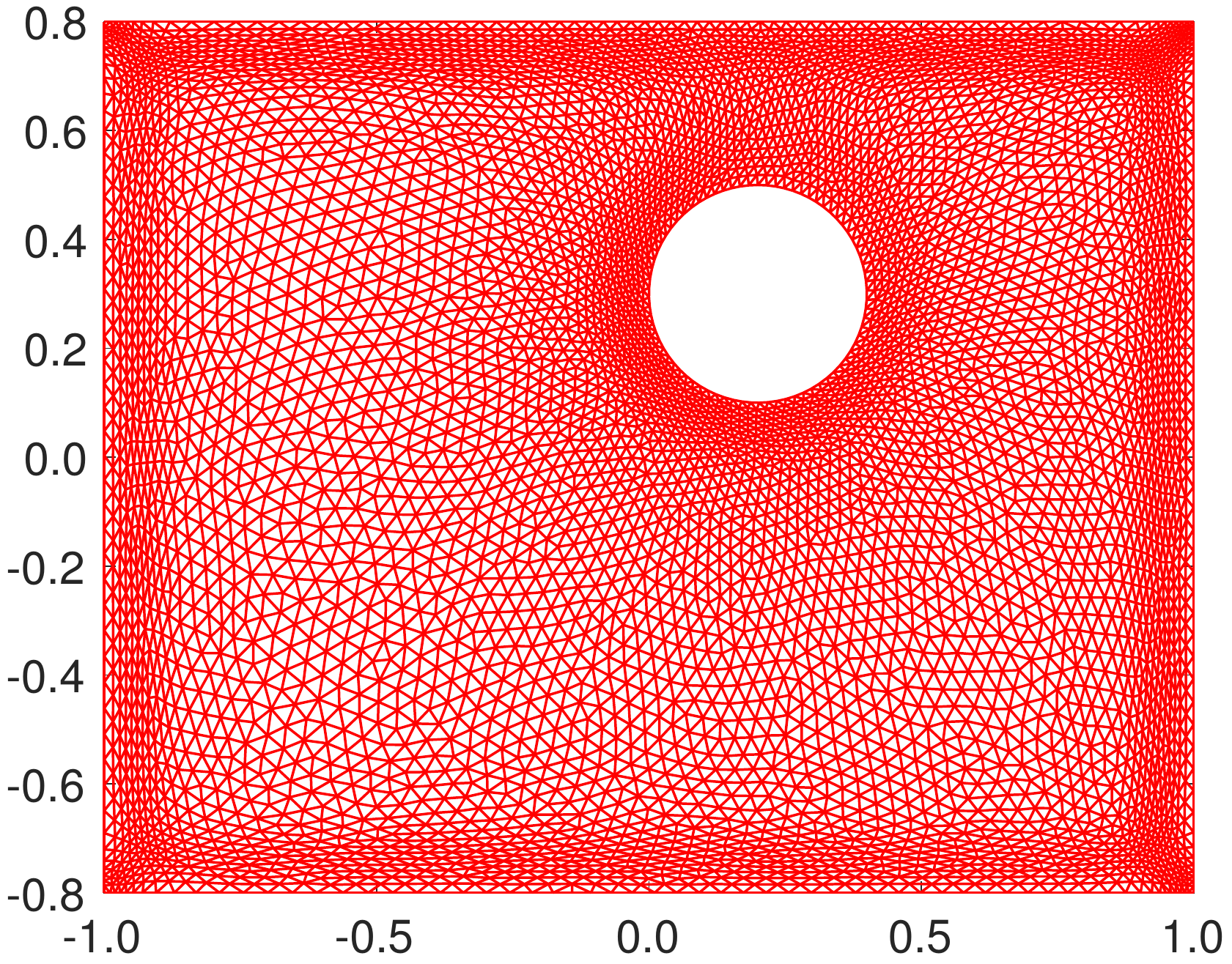}}\qquad
\subfigure[Mesh at $t = 0.275$.]{\includegraphics[width = 0.3\textwidth]{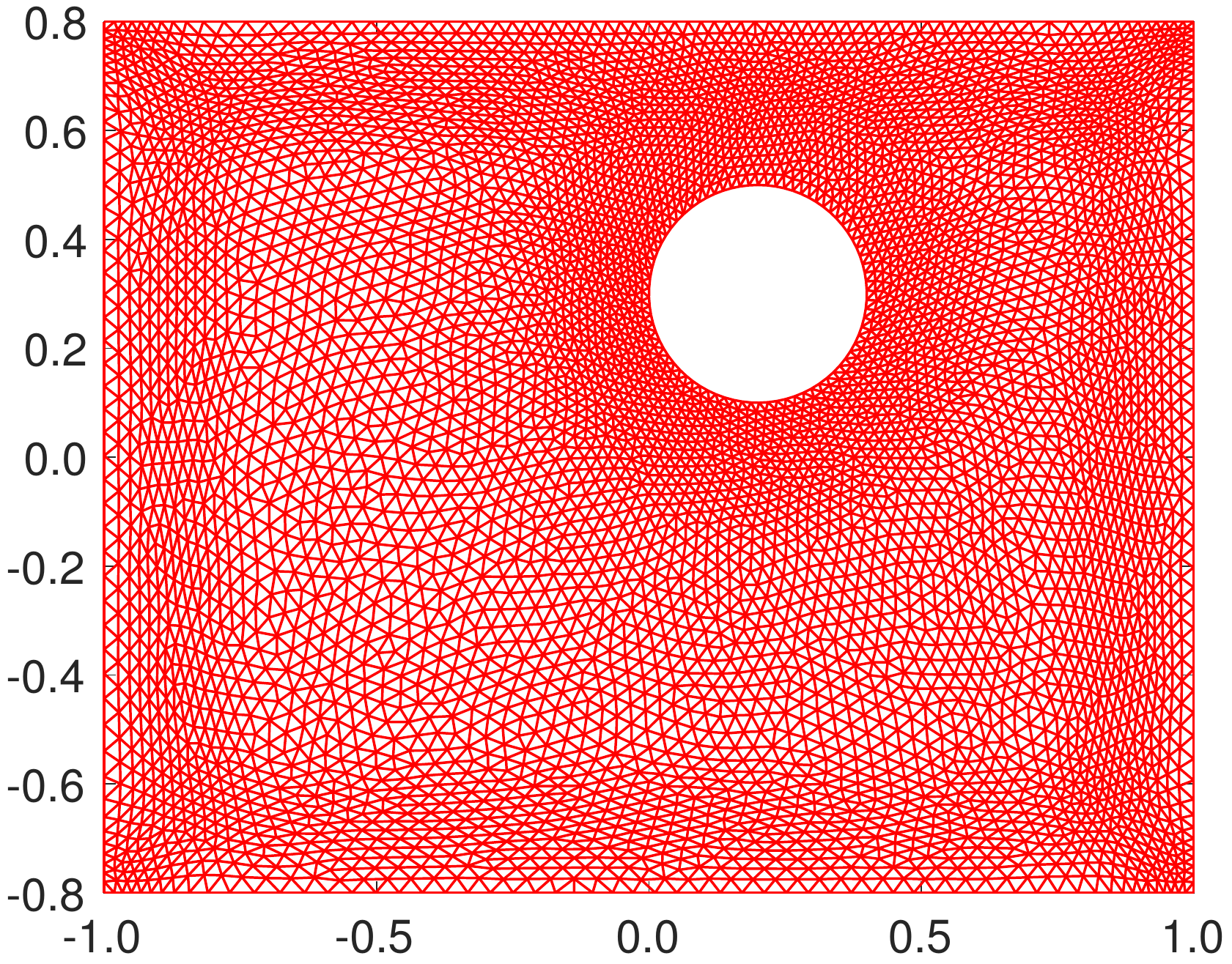}}\qquad
\subfigure[Mesh at $t = 0.311$.]{\includegraphics[width = 0.3\textwidth]{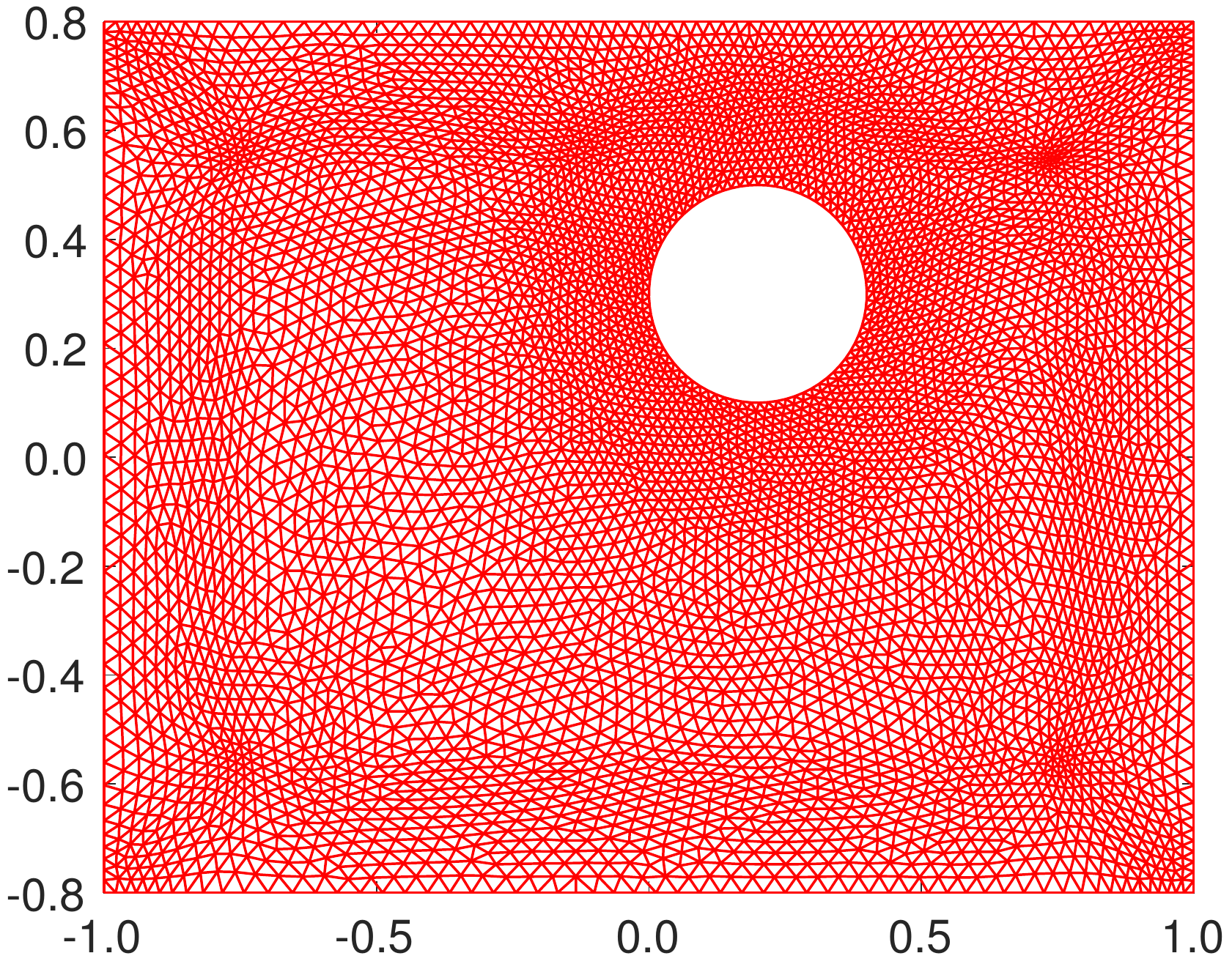}}
\subfigure[Solution at $t = 0.035$.]{\includegraphics[width = 0.3\textwidth]{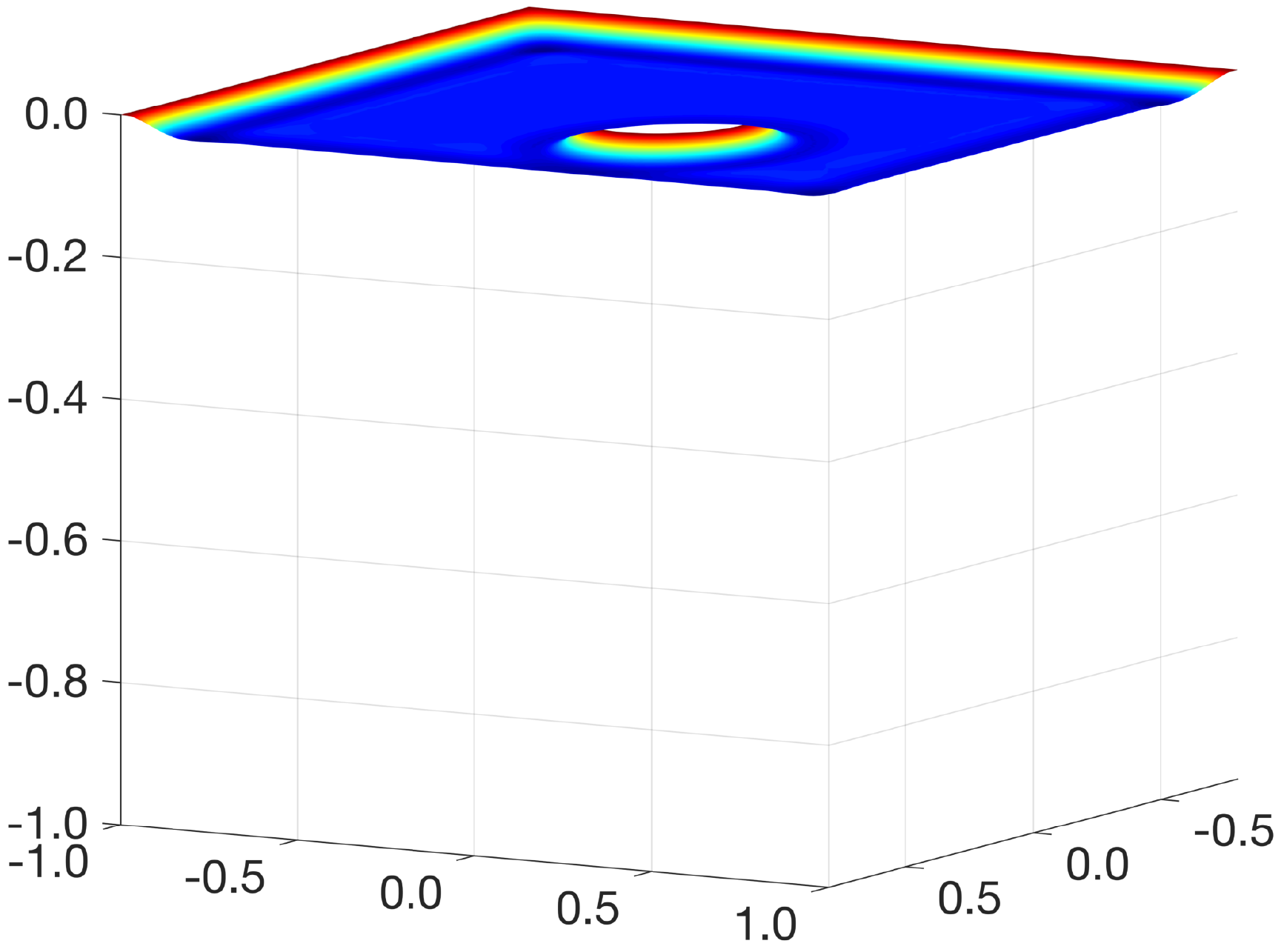}}\qquad
\subfigure[Solution at $t = 0.275$.]{\includegraphics[width = 0.3\textwidth]{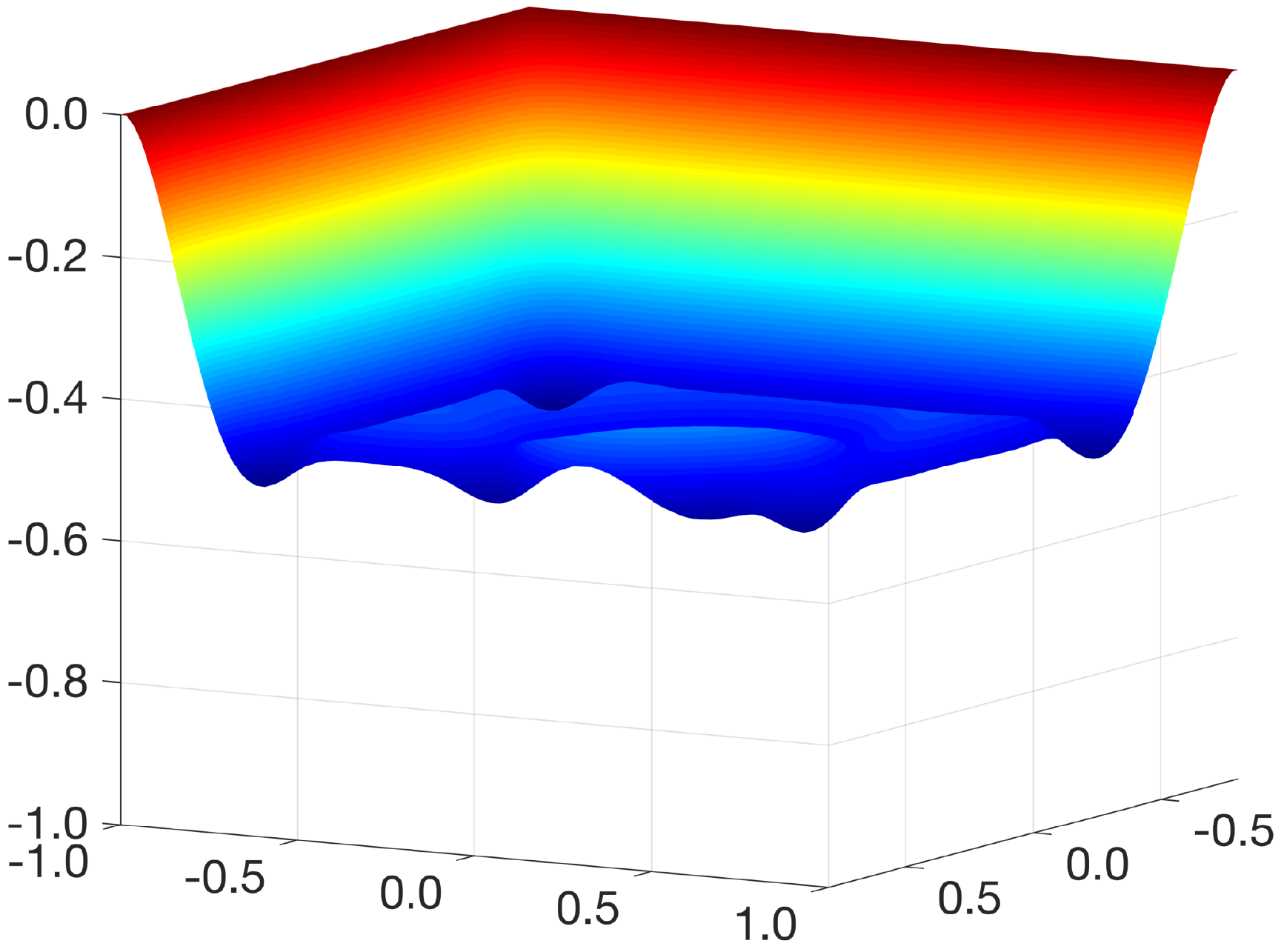}}\qquad
\subfigure[Solution at $t = 0.311$.]{\includegraphics[width = 0.3\textwidth]{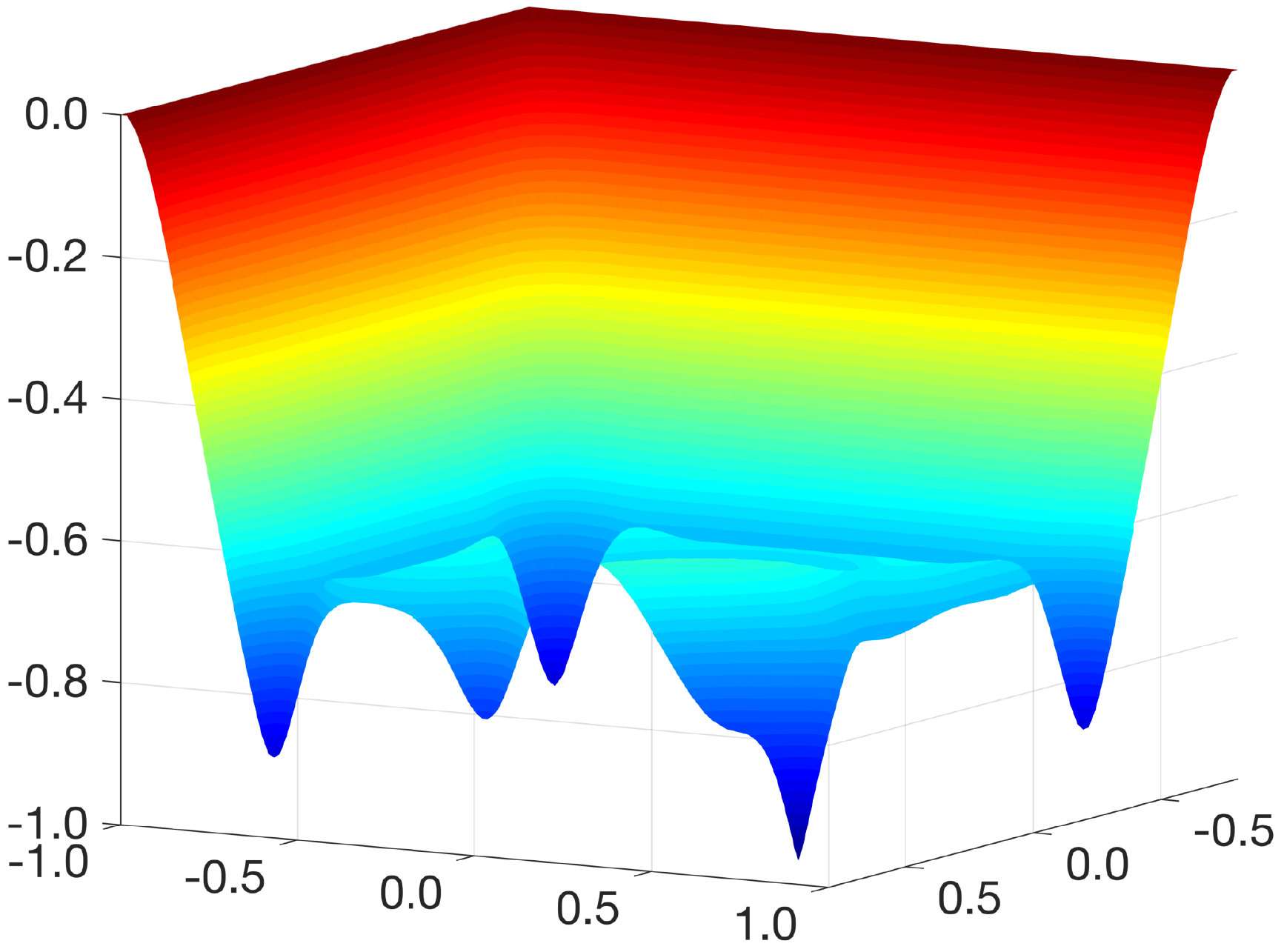}}
\caption{The evolution of the solution for $\eps = 0.01$. The mesh size is $N = 11658$. \label{r2e-2}}
\end{figure}

One interesting observation from the skeleton and singularity points shown in Fig.~\ref{recsk} is that the track of the first touchdown point does not vary continuously with $\eps$. We observe that the first singularity point switches between branches several times suggesting that multiple simultaneous singularities are possible only at fixed values of $\eps$. In Figs.~\ref{fig:RectHole} and \ref{r2e-2}, single point touchdown is observed, however, other troughs in the solution are also very close to singularity which helps explain the sensitivity of the touchdown set on $\eps$.

%For example, $\eps \in \{ 0.42, 0.44, 0.46\}$, the touchdown points concentrate in the fourth quadrant
%while when $\eps$ takes values greater than or equal to $0.48$, the singularities jump to the third quadrant.
%For $\eps = 0.42$, $0.44$, $0.46$, $0.48$, there are actually two points in which the solution decreases faster,
%and one is around Mark 9 in Fig. \ref{recsk}, another one is in region around Mark 10 in Fig. \ref{recsk}.
%Throughout the whole evolution, values of $u$ at these two points compete on touching down faster,
%and can be close at this two points (see Fig. \ref{two}). Thus it is reasonable to speculate that
%for some $\eps \in (0.46,0.48)$, the solution touches down at two points.

%\begin{figure}[htbp]
%\centering
%\subfigure[]{\includegraphics[width = 0.45\textwidth]{mrechole42-2.eps}}\qquad
%\subfigure[]{\includegraphics[width = 0.45\textwidth]{mrechole48-2.eps}}
%\caption{The picture on the left shows singularity for $\eps = 0.042 $ and the picture on the right shows singularity of $\eps = 0.048$, $u$ values are close to $-1$ at the two points in both pictures. \label{two}}
%\end{figure}

%\begin{figure}[htbp]
%\centering
%\subfigure[]{\includegraphics[width = 0.3\textwidth]{rec0mesh1.eps}}\qquad
%\subfigure[]{\includegraphics[width = 0.3\textwidth]{rec0mesh2.eps}}\qquad
%\subfigure[]{\includegraphics[width = 0.3\textwidth]{rec0mesh3.eps}}
%\subfigure[]{\includegraphics[width = 0.3\textwidth]{rec0u1.eps}}\qquad
%\subfigure[]{\includegraphics[width = 0.3\textwidth]{rec0u2.eps}}\qquad
%\subfigure[]{\includegraphics[width = 0.3\textwidth]{rec0u3.eps}}\caption{The evolution of the solution for $\eps = 5 \times 10^{-4}$. The mesh size is $N = 31236$. \label{5e-4}}
%\end{figure}

%\begin{figure}[htbp]
%\centering
%
%\subfigure[]{\includegraphics[width = 0.3\textwidth]{rec2mesh1.eps}}\qquad
%\subfigure[]{\includegraphics[width = 0.3\textwidth]{rec2mesh2.eps}}\qquad
%\subfigure[]{\includegraphics[width = 0.3\textwidth]{rec2mesh3.eps}}\\
%\subfigure[]{\includegraphics[width = 0.3\textwidth]{rec2u1.eps}}\qquad
%\subfigure[]{\includegraphics[width = 0.3\textwidth]{rec2u2.eps}}\qquad
%\subfigure[]{\includegraphics[width = 0.3\textwidth]{rec2u3.eps}}
%\caption{The evolution of the solution for $\eps = 0.2$. The mesh size $N = 15401$. \label{r1e-1}}
%
%\end{figure}

\clearpage

\vspace{10pt}

% Example 3.
\begin{exam}[Rectangular domain with four holes]
{\em We consider an example based on a domain similar to the real MEMS device in Fig.~\ref{fig:intro_diagram_B} by solving equation \eqref{eq:intro} on the rectangular domain $\Omega = (-1,1)\times (-0.8,0.8)$ punctured by four circular holes. The holes have common radius $0.15$ and are arranged symmetrically at the points $\pm(0.5,\pm 0.3)$.
A similar modification to the metric tensor as in \eqref{M-2} has been used for this example
to keep a level of mesh concentration near the holes.
}

\end{exam}

\begin{figure}[htbp]
\centering
\includegraphics[width = 0.4\textwidth]{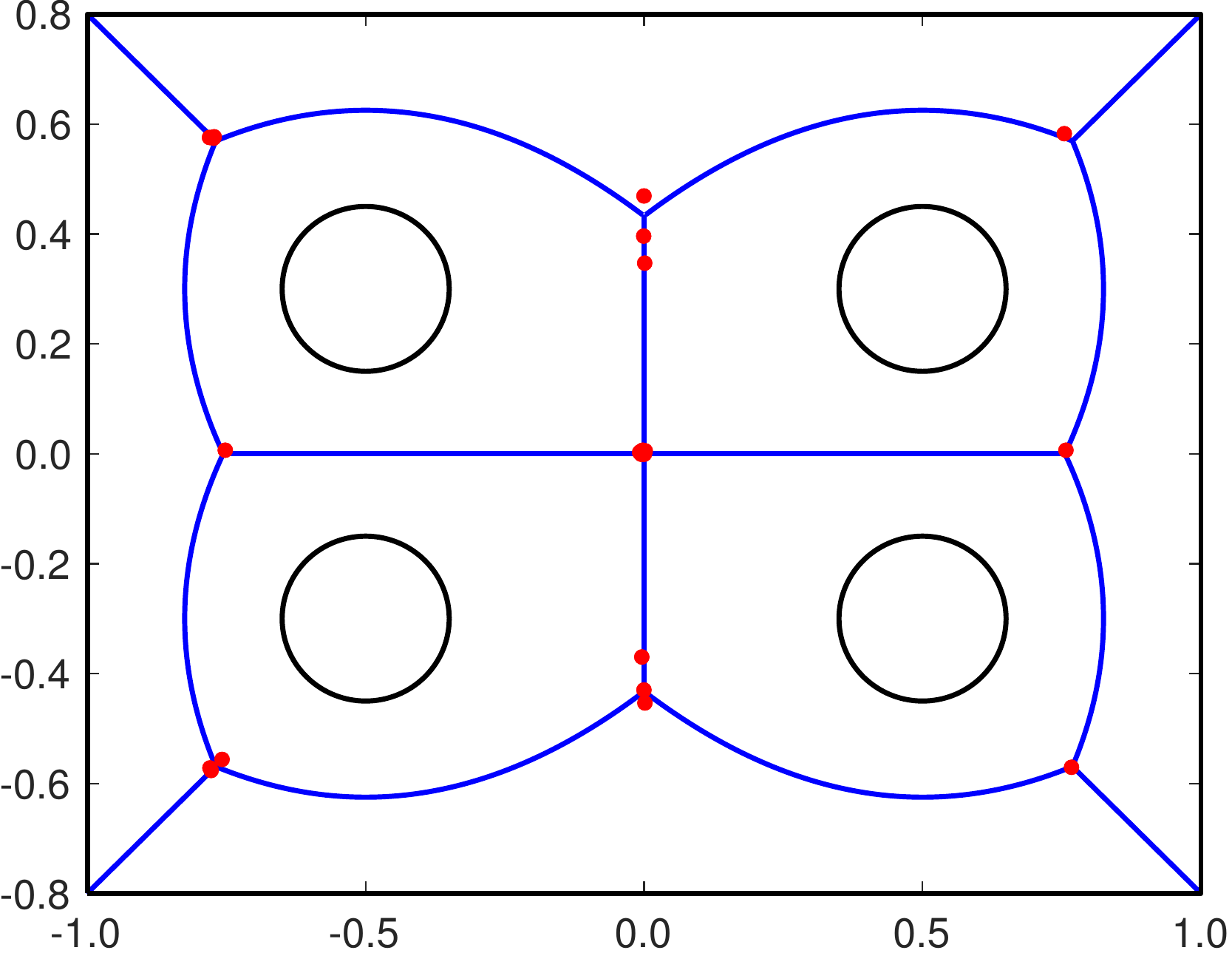}
\caption{Skeleton (blue curve) and touchdown points (red points) for the rectangular domain with four holes. \label{fourholesk}}
\end{figure}

\begin{figure}[!ht]
\centering
\subfigure[Mesh at $t = 0.095$.]{\includegraphics[width = 0.3\textwidth]{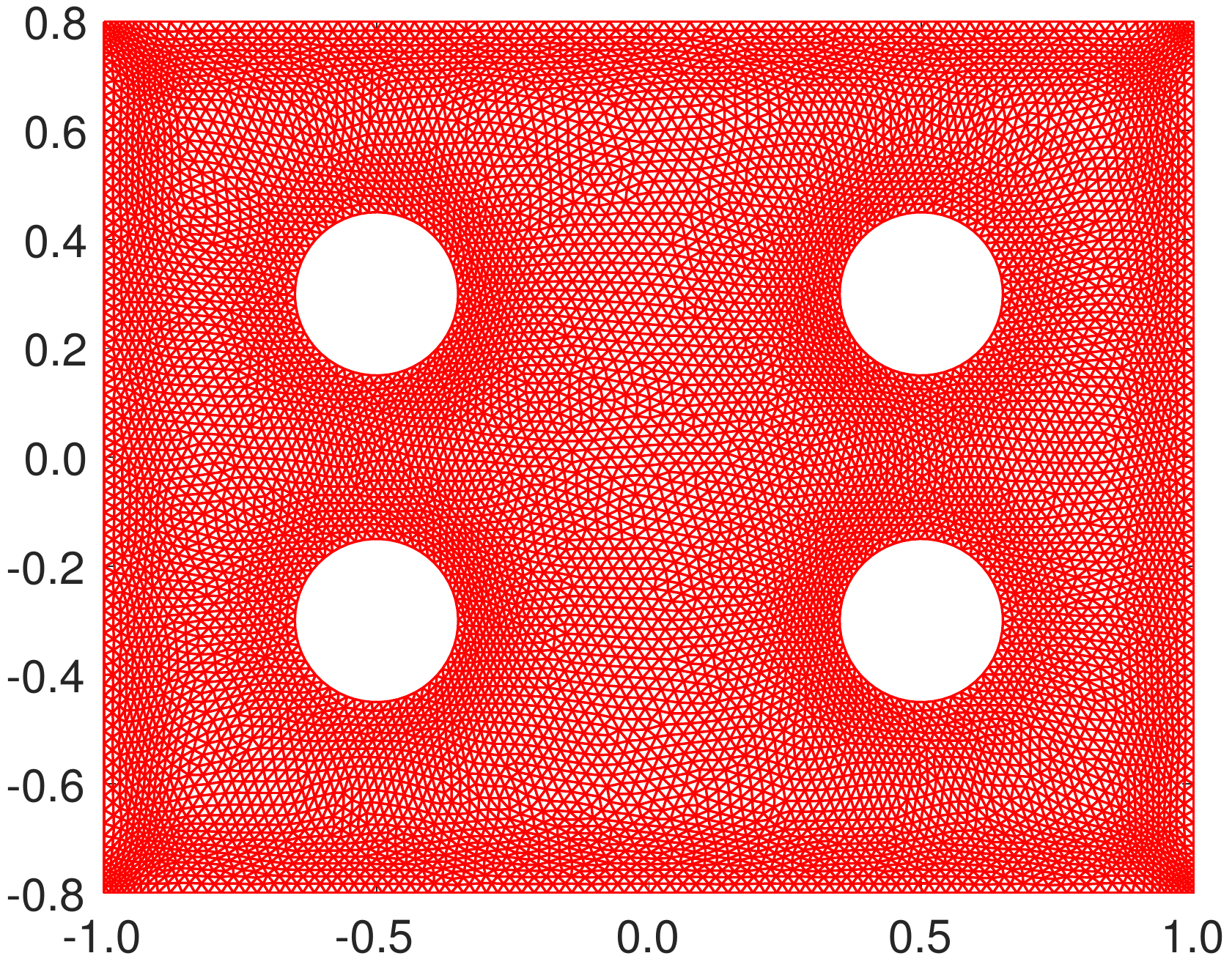}}\qquad
\subfigure[Mesh at $t = 0.275$.]{\includegraphics[width = 0.3\textwidth]{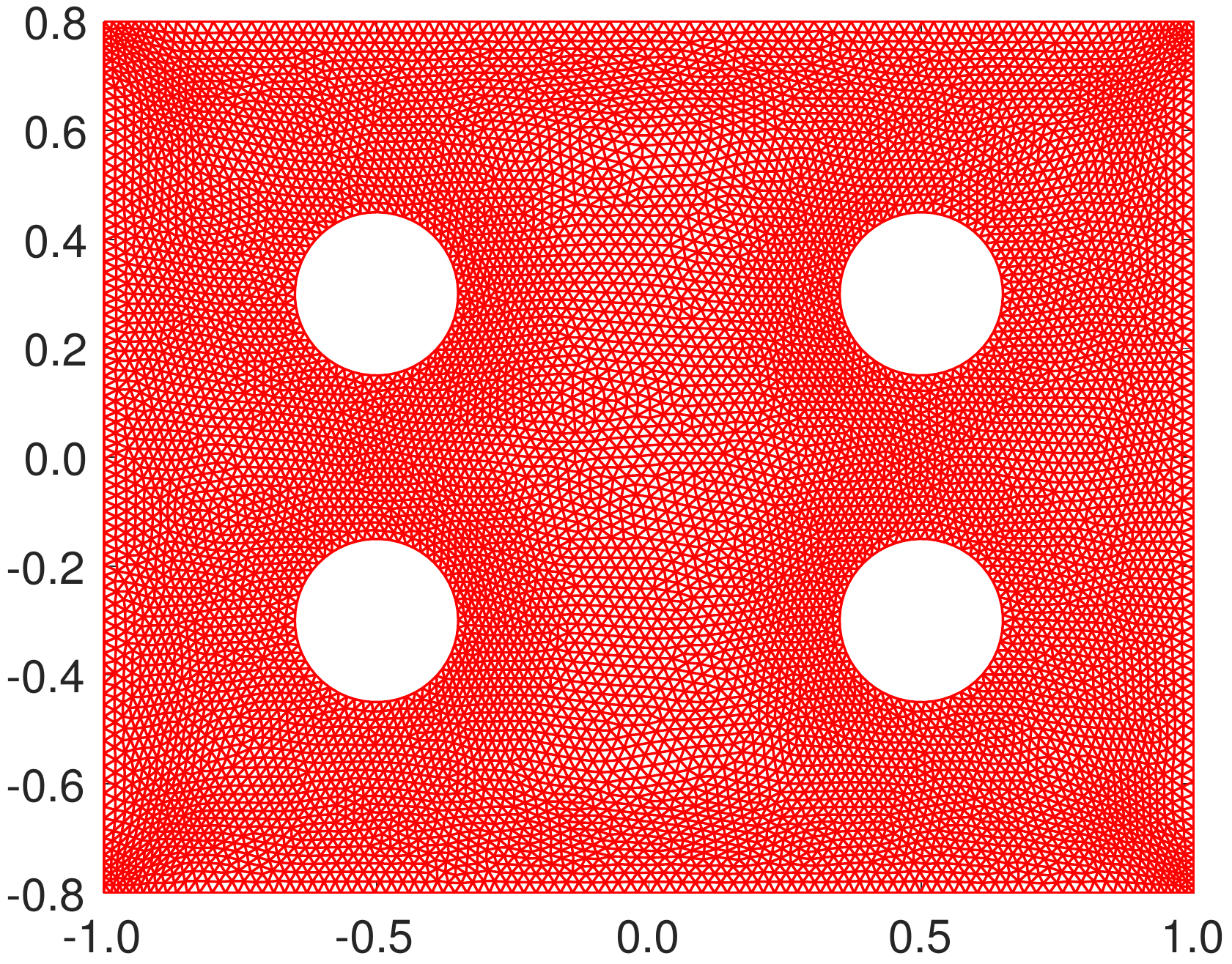}}\qquad
\subfigure[Mesh at $t = 0.307$.]{\includegraphics[width = 0.3\textwidth]{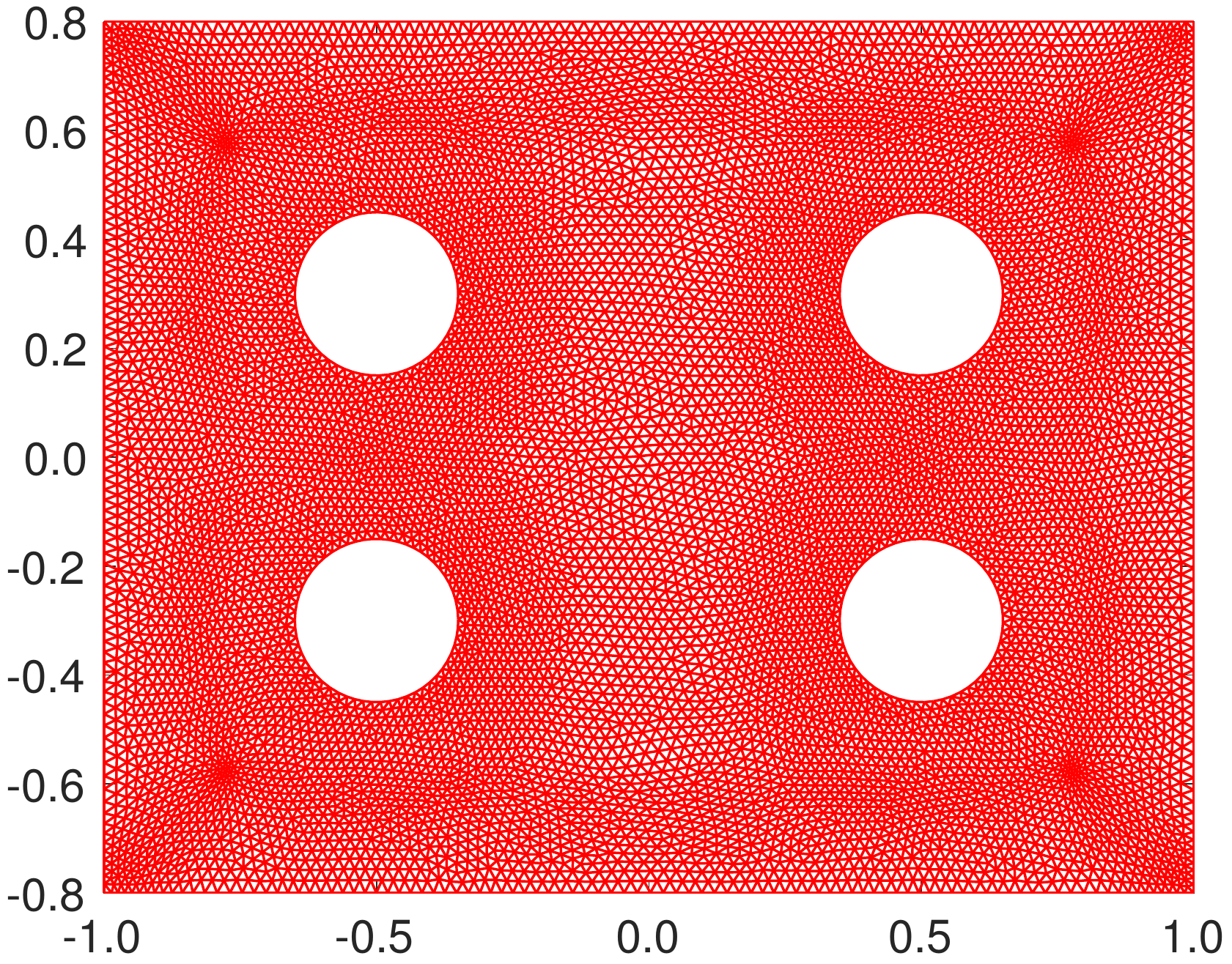}}
\subfigure[Solution at $t = 0.095$.]{\includegraphics[width = 0.3\textwidth]{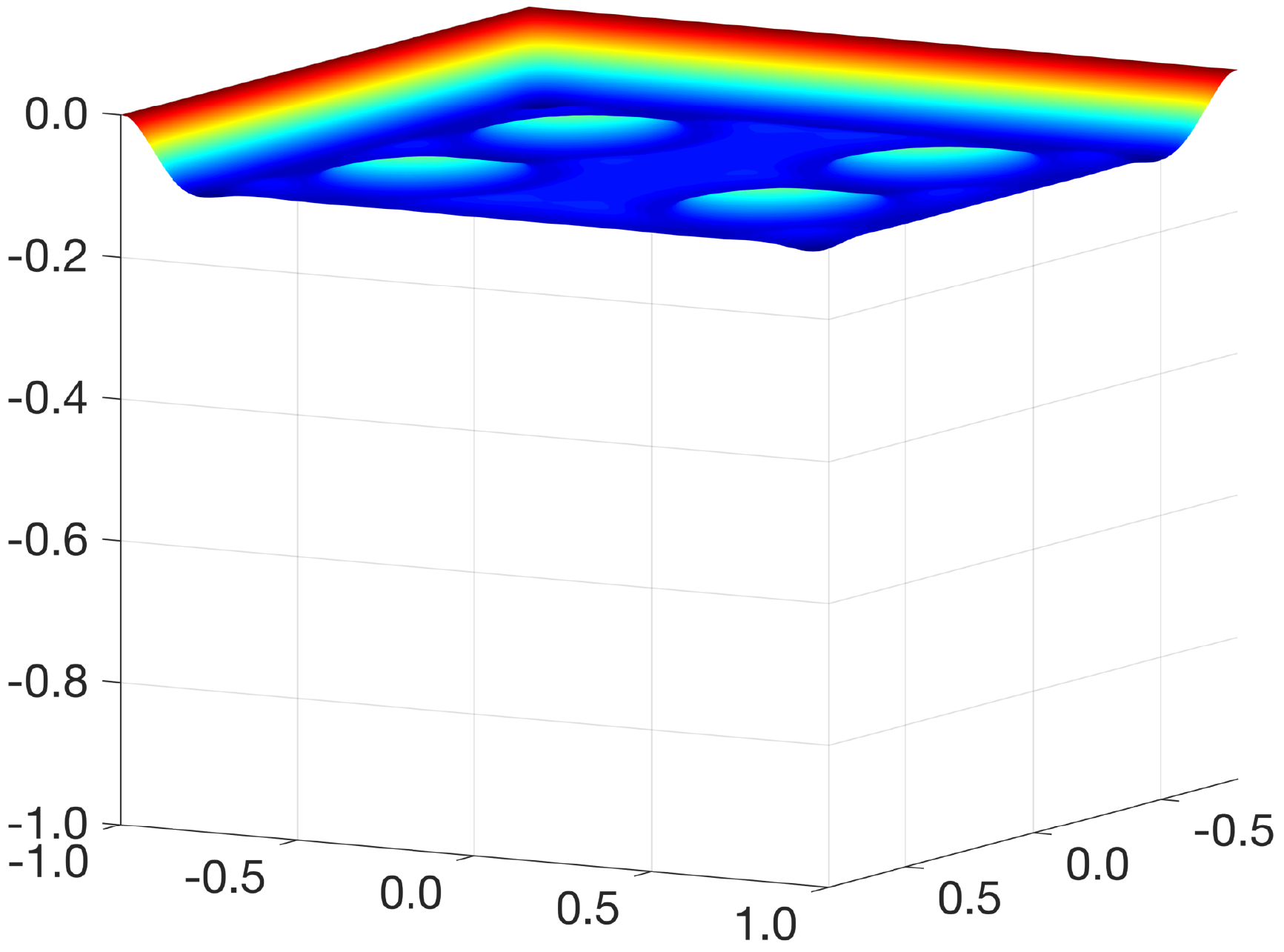}}\qquad
\subfigure[Solution at $t = 0.275$.]{\includegraphics[width = 0.3\textwidth]{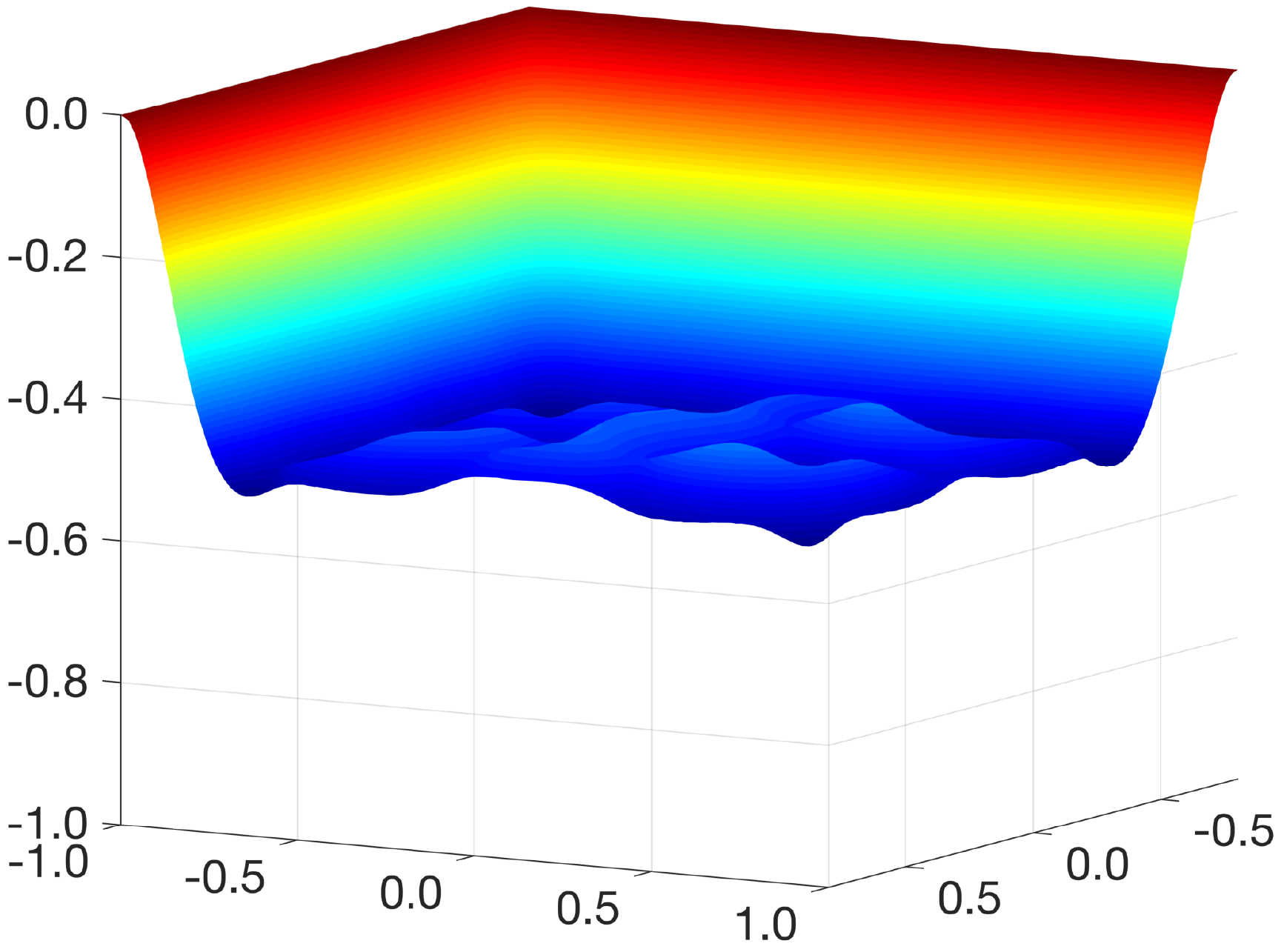}}\qquad
\subfigure[Solution at $t = 0.307$.]{\includegraphics[width = 0.3\textwidth]{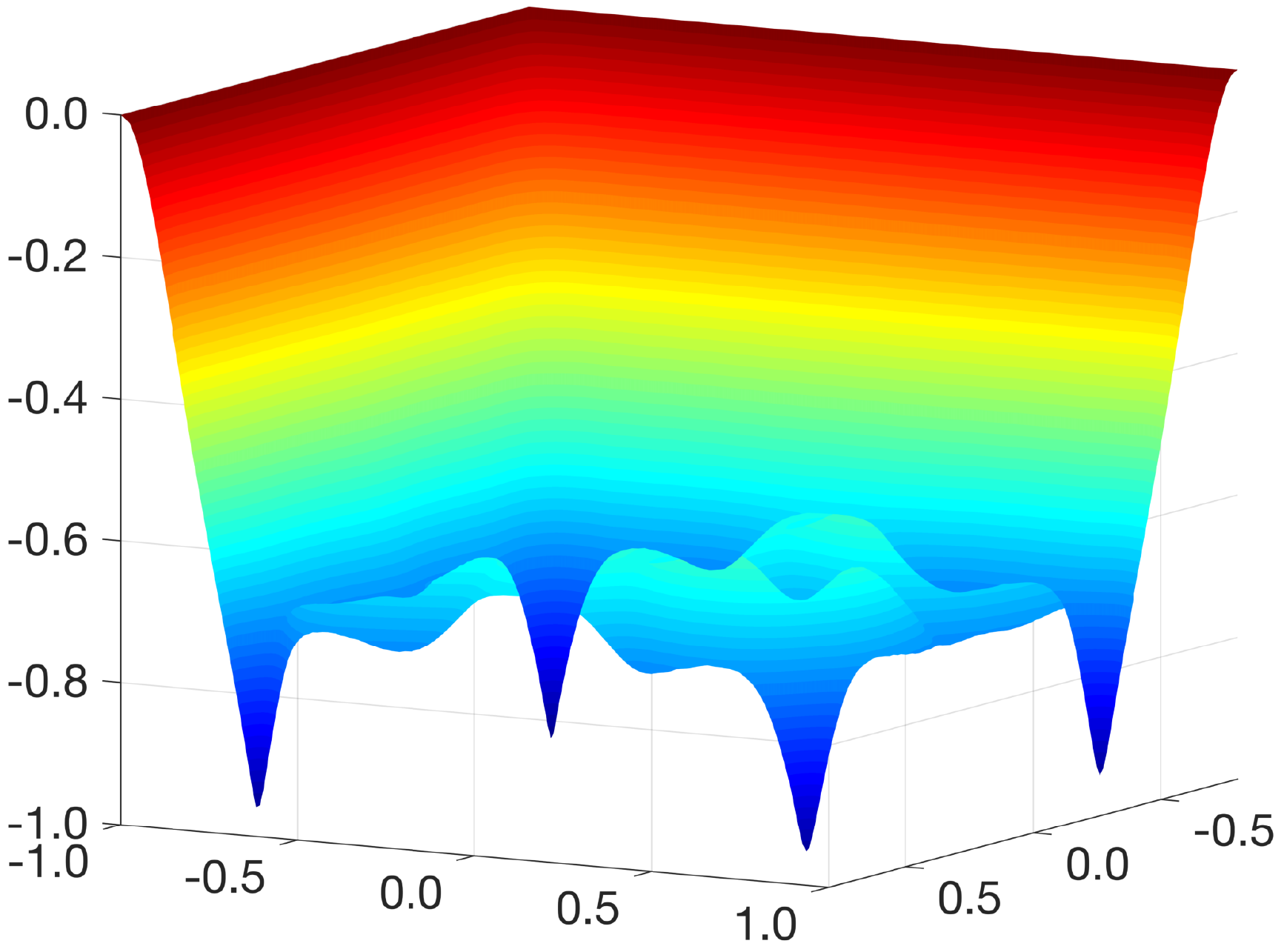}}
\caption{The evolution of the solution for $\eps = 0.007$. The mesh size is $N = 20018$. \label{fh7e-3}}
\end{figure}

\begin{figure}[!ht]
\centering
\subfigure[Mesh at $t = 0.092$.]{\includegraphics[width = 0.3\textwidth]{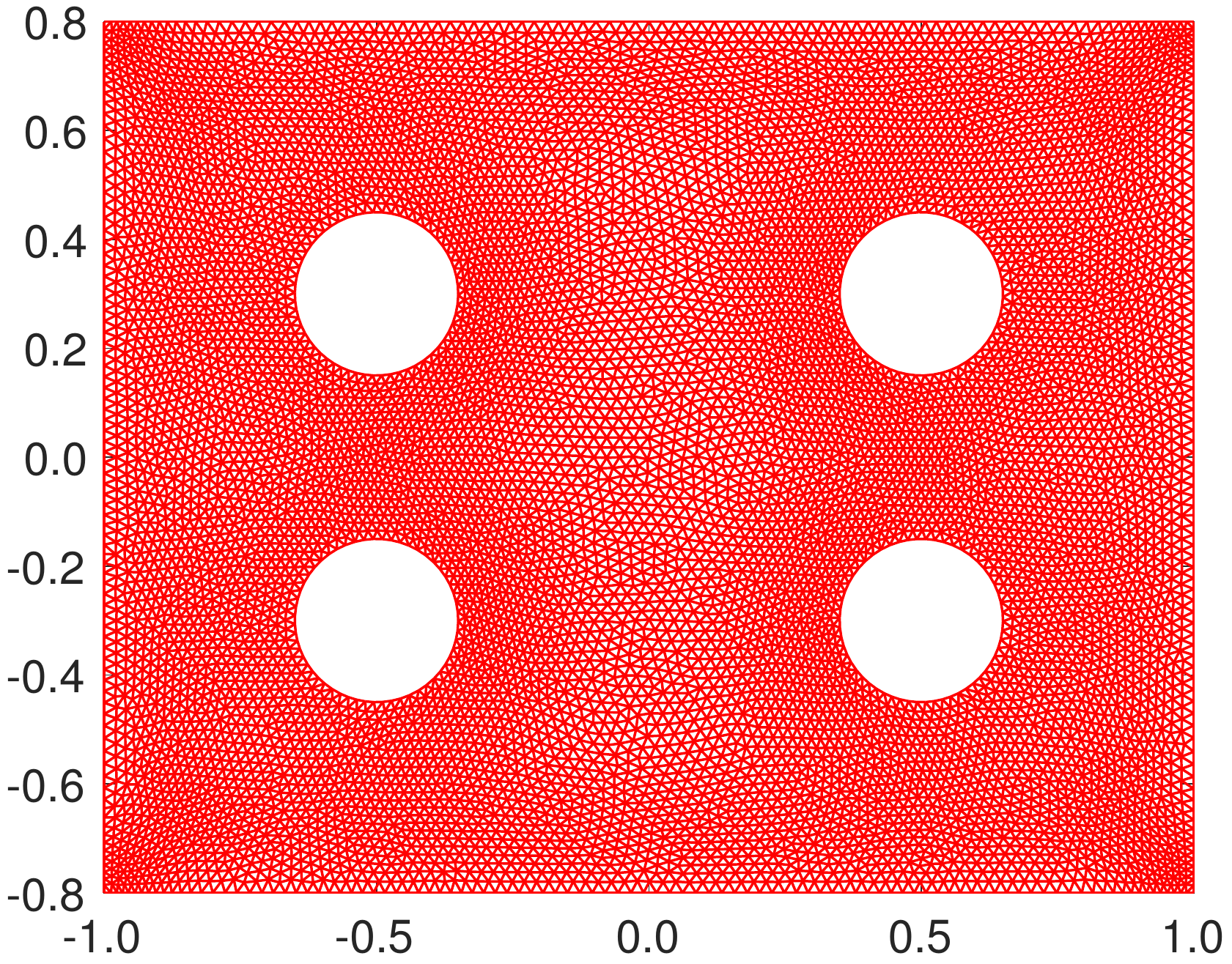}}\qquad
\subfigure[Mesh at $t = 0.304$.]{\includegraphics[width = 0.3\textwidth]{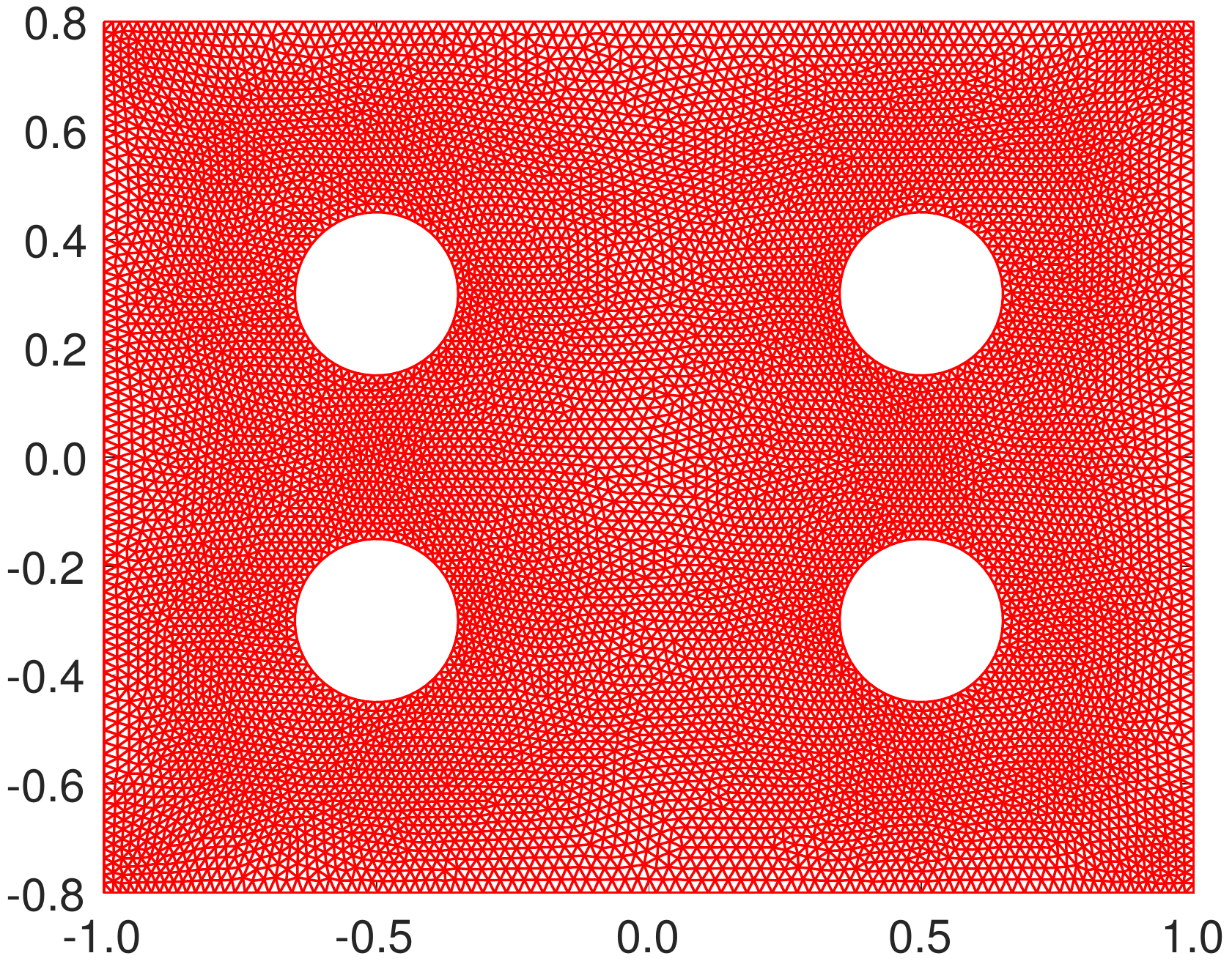}}\qquad
\subfigure[Mesh at $t = 0.309$.]{\includegraphics[width = 0.3\textwidth]{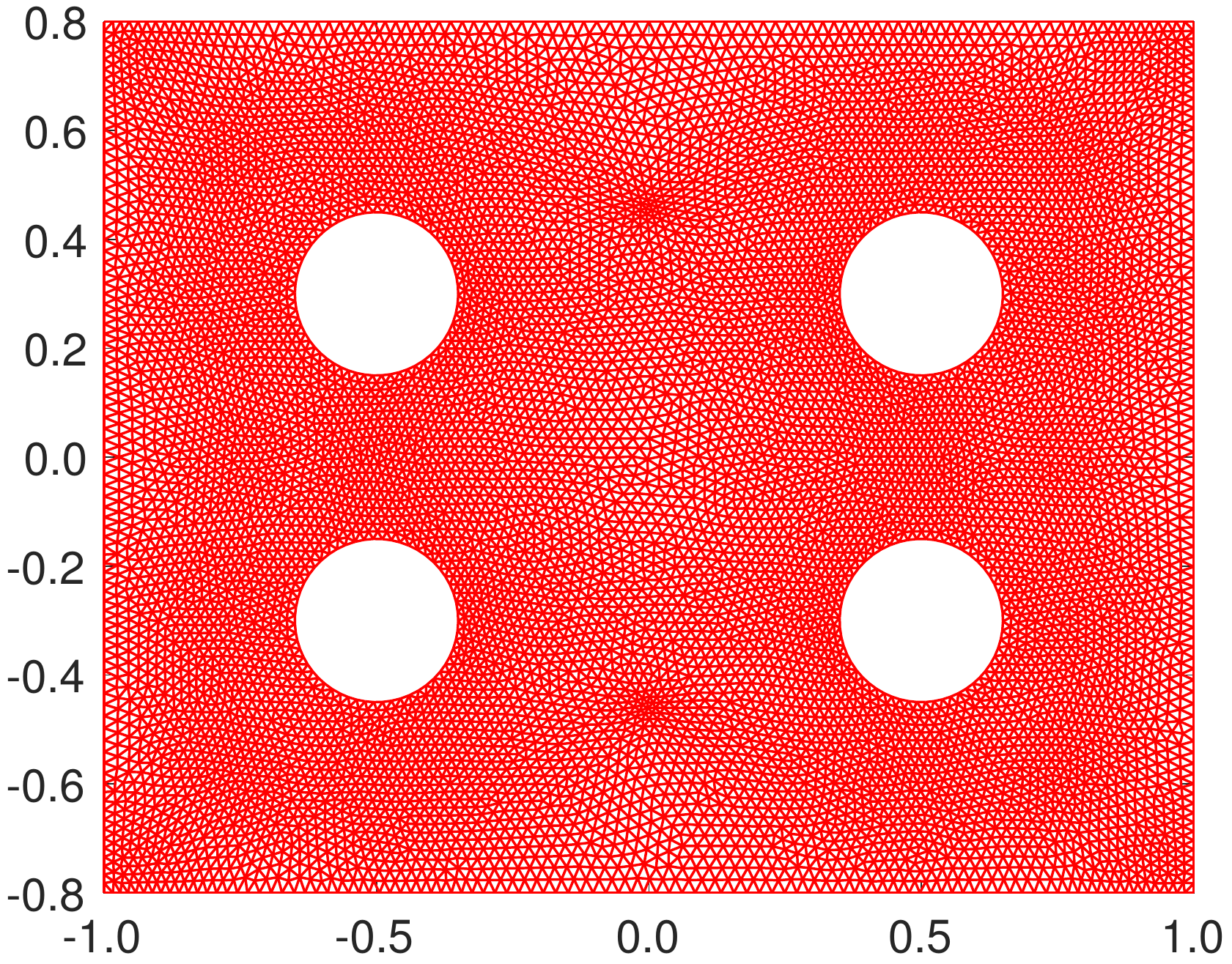}}
\subfigure[Solution at $t = 0.092$.]{\includegraphics[width = 0.3\textwidth]{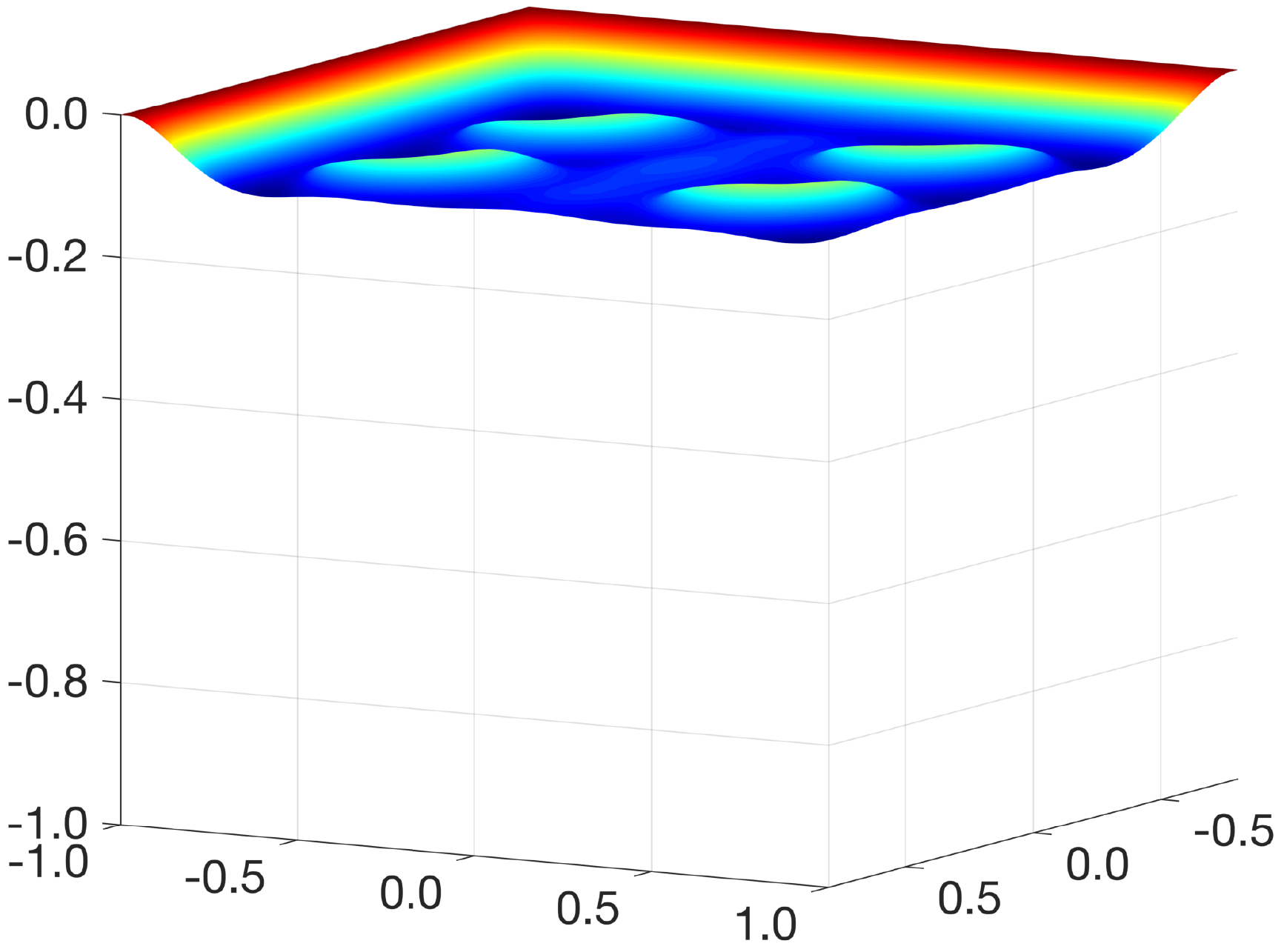}}\qquad
\subfigure[Solution at $t = 0.304$.]{\includegraphics[width = 0.3\textwidth]{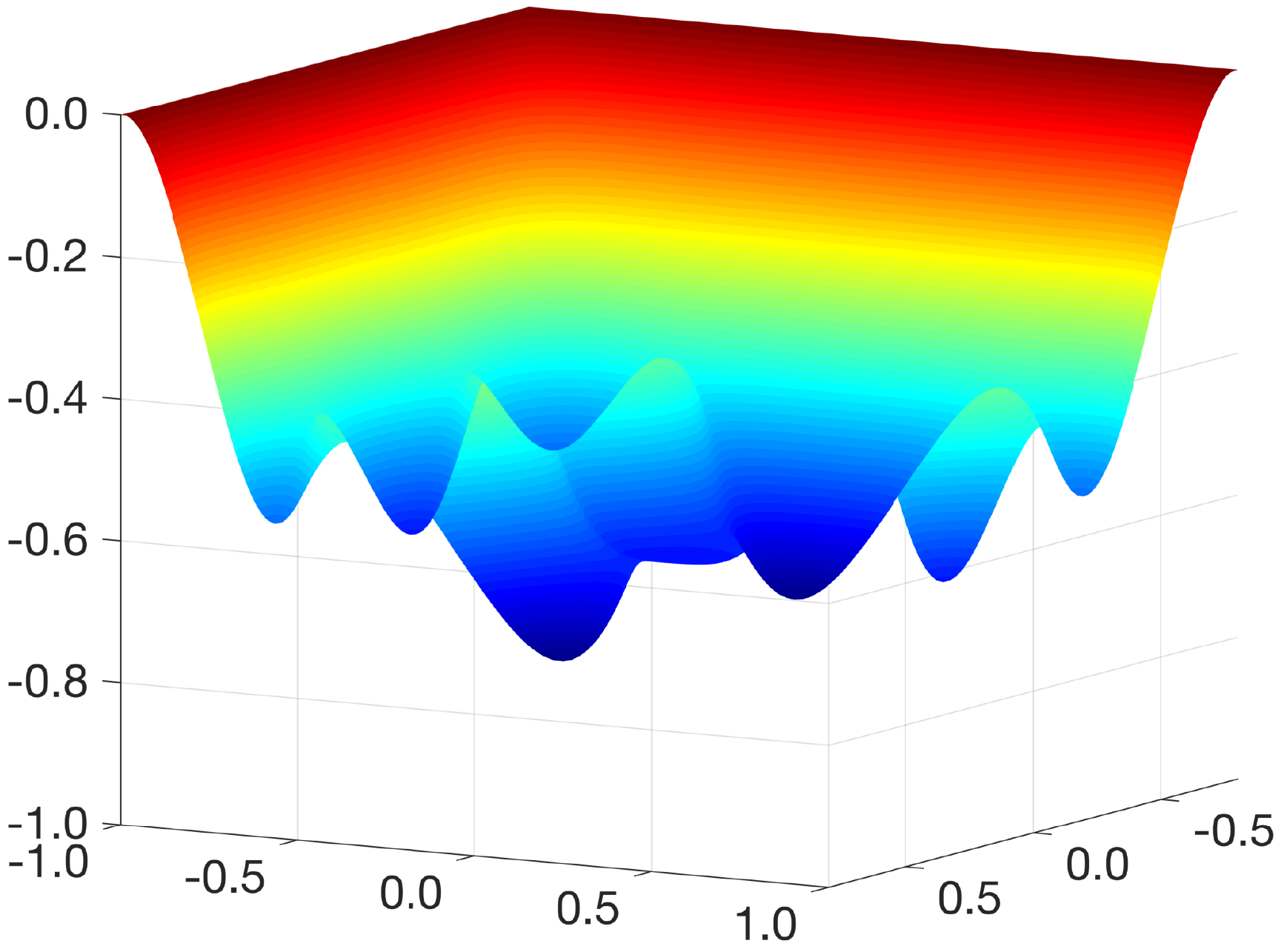}}\qquad
\subfigure[Solution at $t = 0.309$.]{\includegraphics[width = 0.3\textwidth]{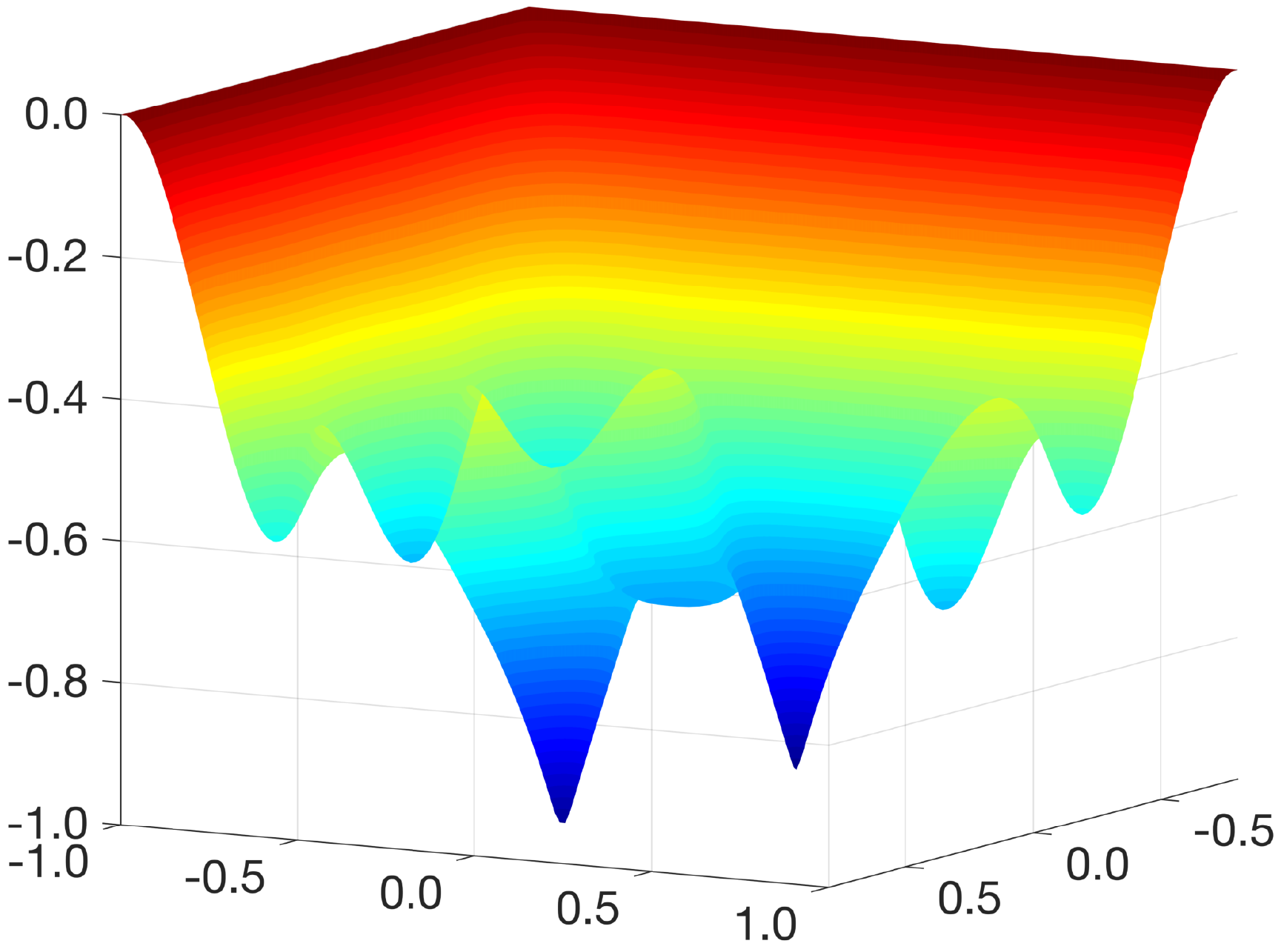}}
\caption{The evolution of the solution for $\eps = 0.02$. The mesh size is $N = 19846$. \label{fh2e-2}}
\end{figure}

\begin{figure}[!ht]
\centering
\subfigure[Mesh at $t = 0.112$.]{\includegraphics[width = 0.3\textwidth]{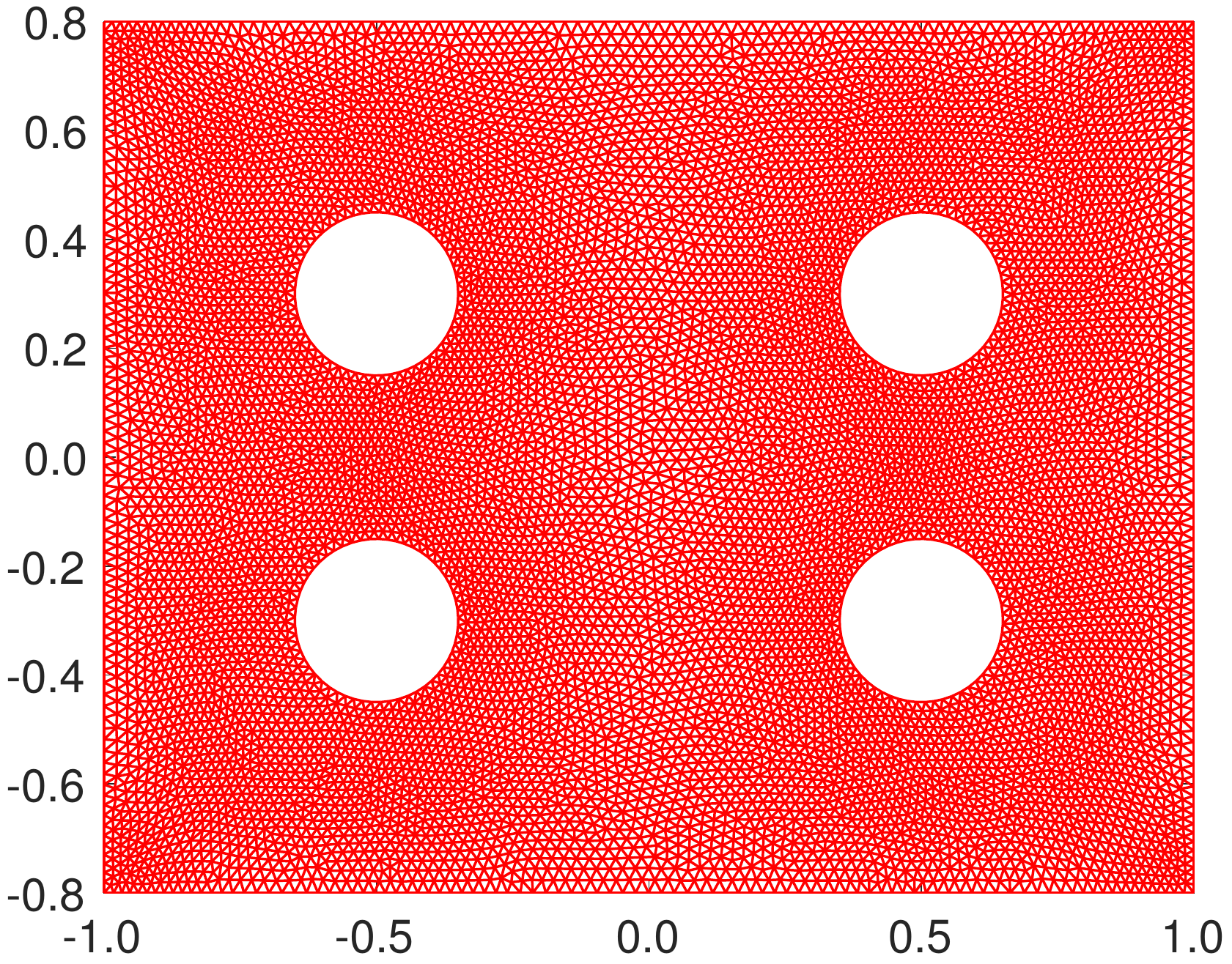}}\qquad
\subfigure[Mesh at $t = 0.282$.]{\includegraphics[width = 0.3\textwidth]{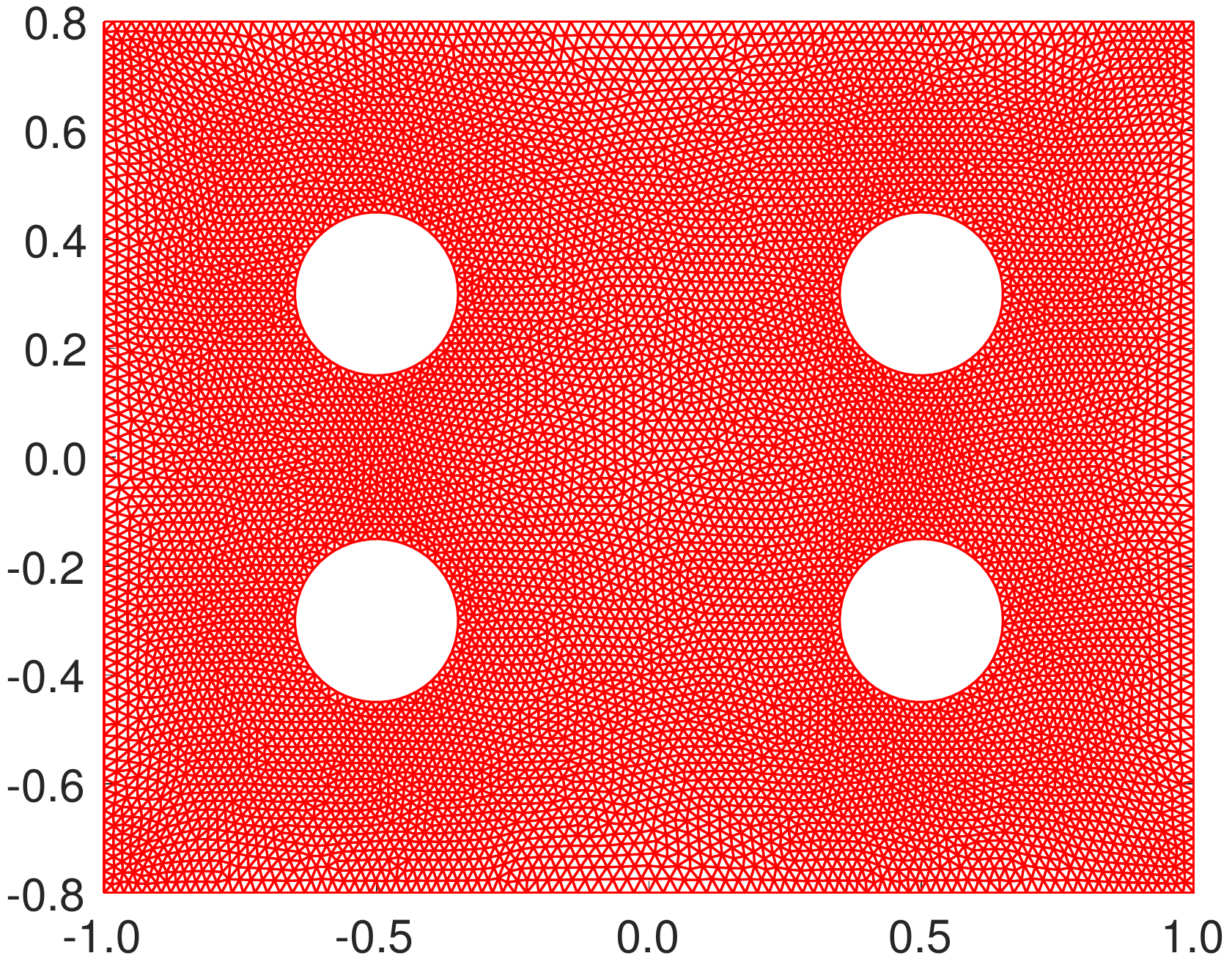}}\qquad
\subfigure[Mesh at $t = 0.307$.]{\includegraphics[width = 0.3\textwidth]{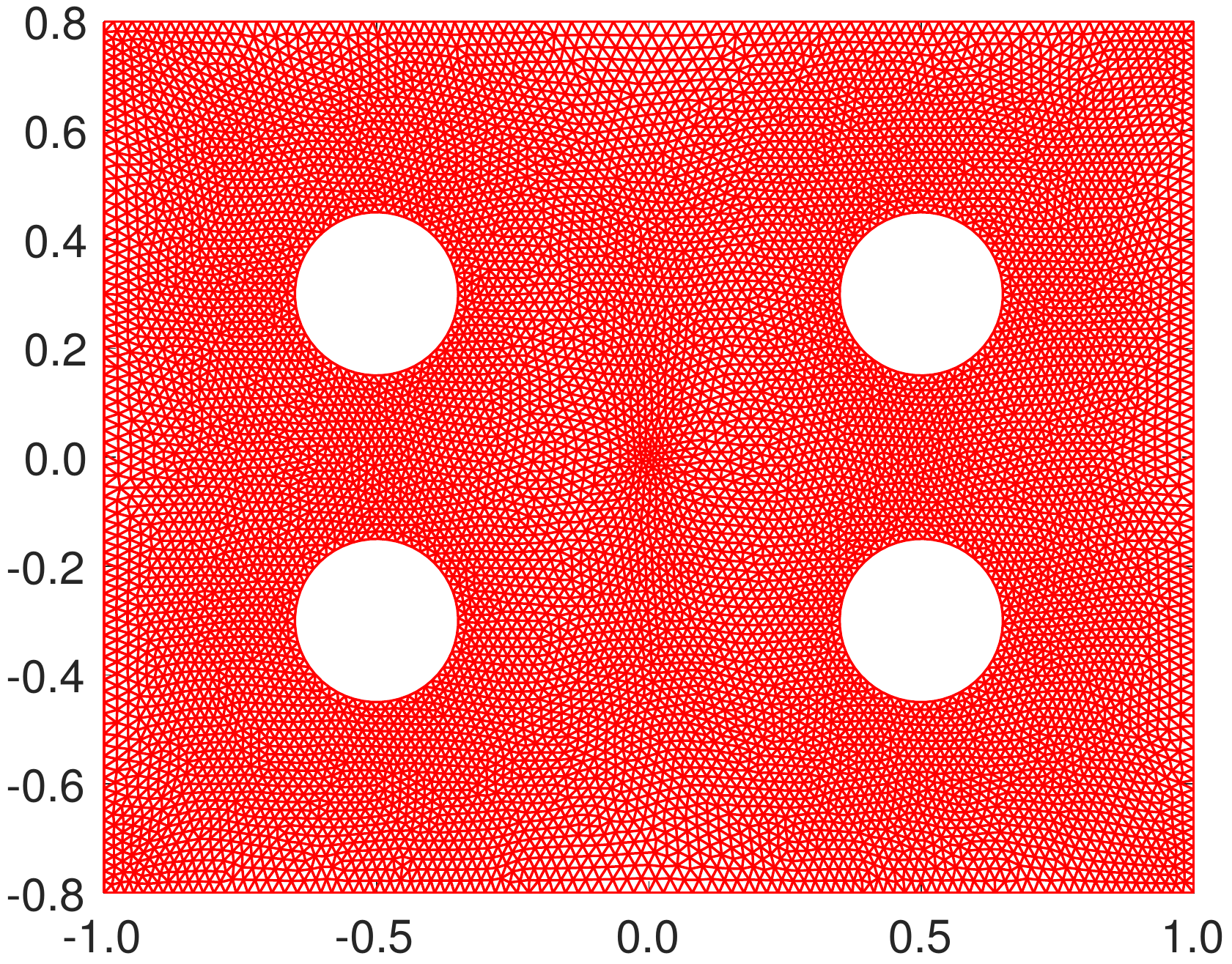}}
\subfigure[Solution at $t = 0.112$.]{\includegraphics[width = 0.3\textwidth]{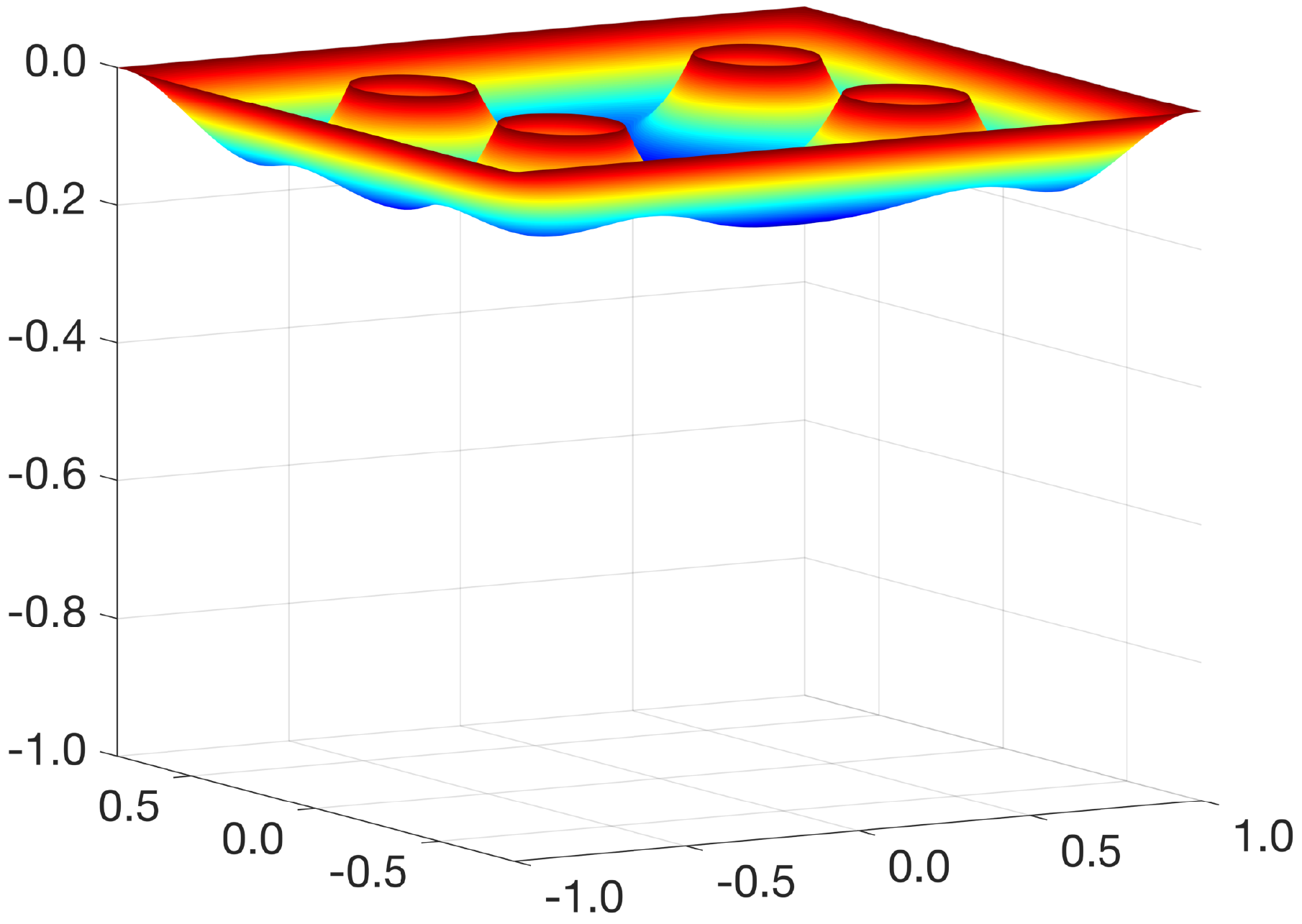}}\qquad
\subfigure[Solution at $t = 0.282$.]{\includegraphics[width = 0.3\textwidth]{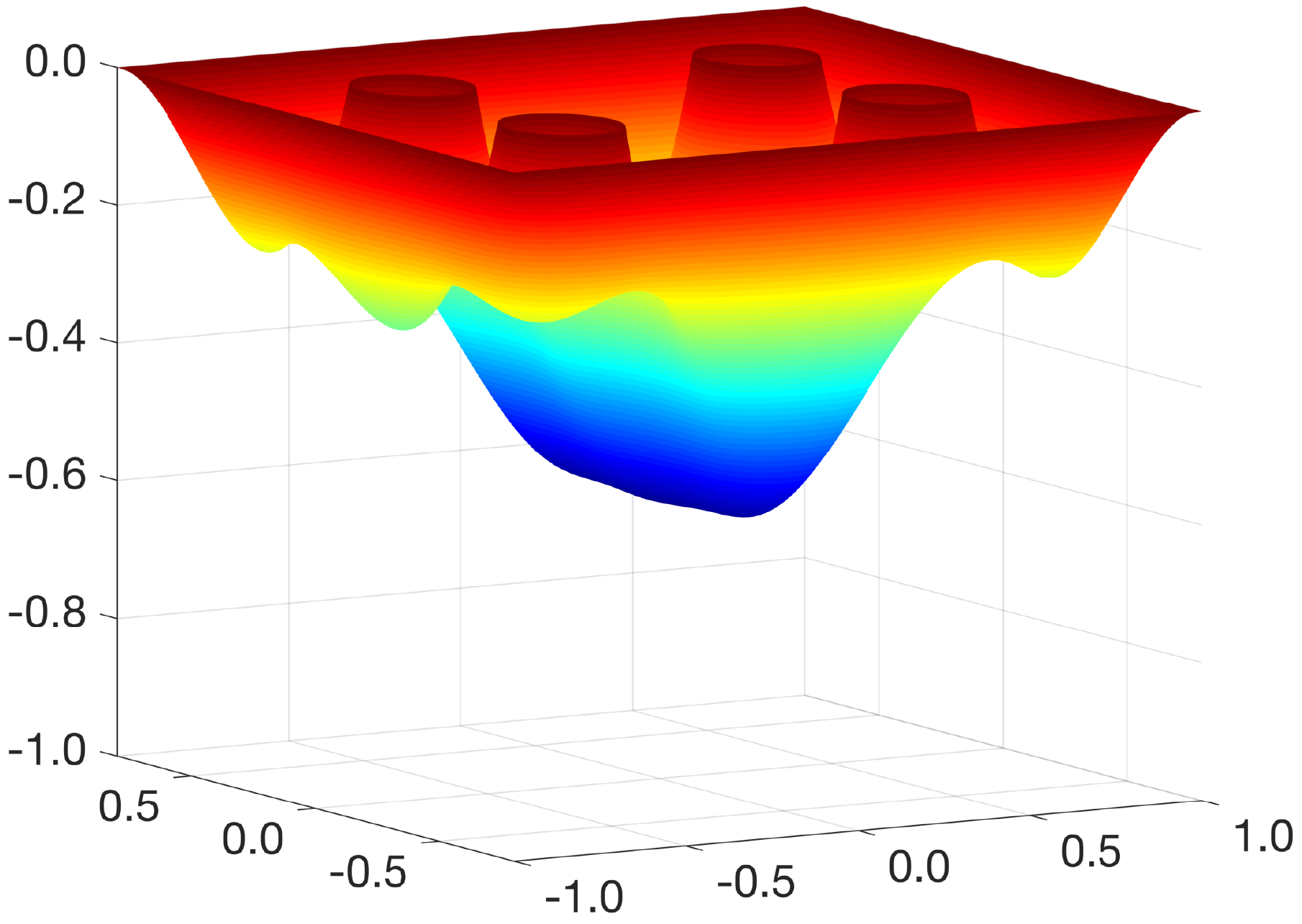}}\qquad
\subfigure[Solution at $t = 0.307$.]{\includegraphics[width = 0.3\textwidth]{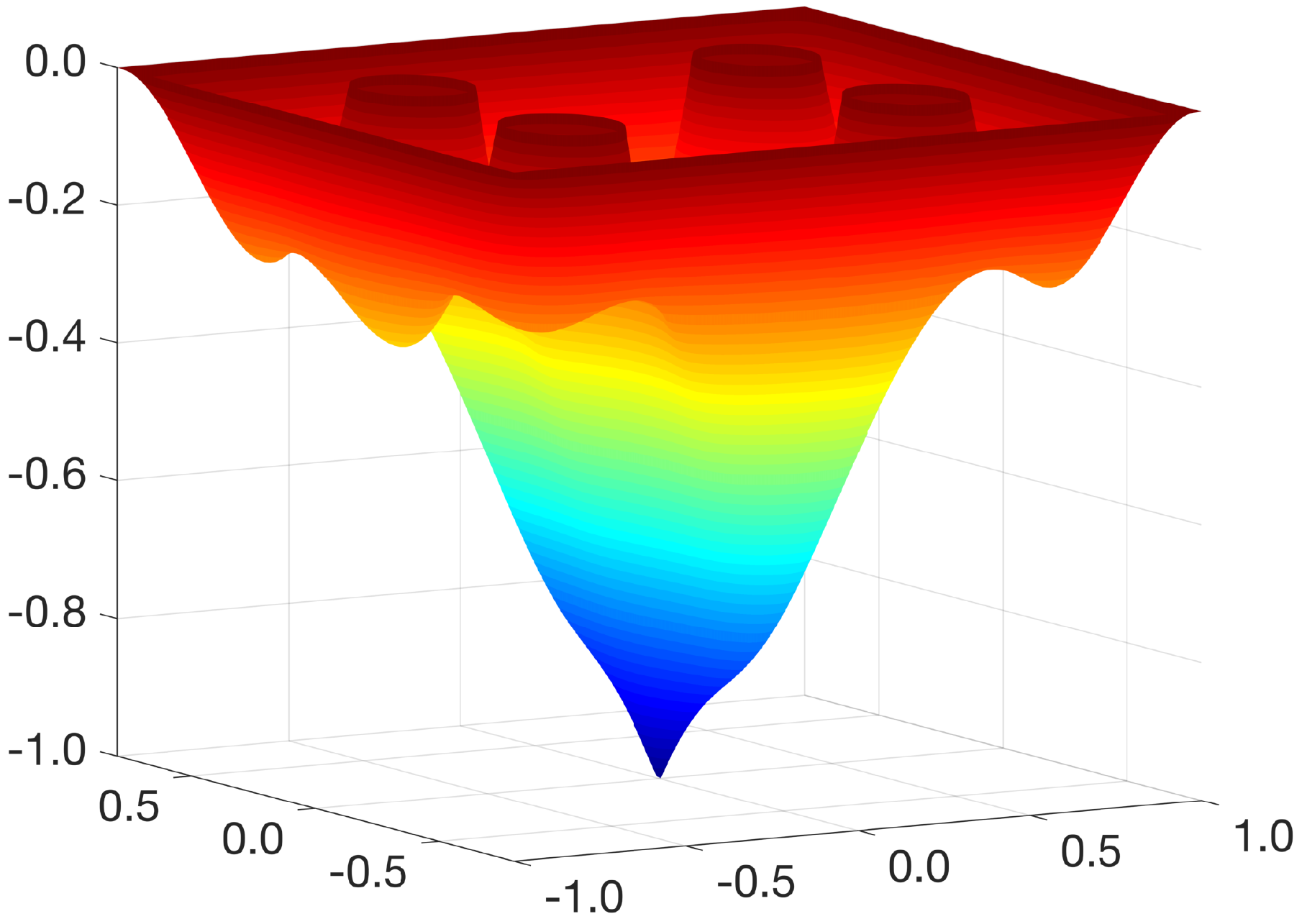}}
\caption{The evolution of the solution for $\eps = 0.04$. The mesh size is $N = 19734$. \label{fh4e-2}}
\end{figure}

We display the skeleton together with the numerically computed singularities obtained for $\eps \in (0.005, 0.1)$ in Fig.~\ref{fourholesk}. For small $\eps$ values, four singularities that are close to the corners are observed. As $\eps$ increases, four singularities merge into two, and eventually into one. For the values $\eps = 0.007$, $\eps = 0.02$, and $\eps = 0.04$, the snapshots of the evolution are presented in Fig.~\ref{fh7e-3}, Fig.~\ref{fh2e-2}, and Fig.~\ref{fh4e-2}. Note that multiple troughs are observed in Fig.~\ref{fh7e-3} and Fig.~\ref{fh2e-2}, while only the smallest touchdown locations are captured.

\vspace{10pt}

%% example 4
\begin{exam}[Asymmetric Domain]
{\em We consider the asymmetric domain given in polar coordinates by}
\begin{equation}\label{eqn:nonparam}
(x,y) = r(\theta)(\cos\theta,\sin\theta), \qquad r(\theta) =  1 + (0.15\sin 2\theta + 0.3\cos 3\theta).
\end{equation}
\end{exam}

\begin{figure}[htbp]
\centering
\includegraphics[width = 0.45\textwidth]{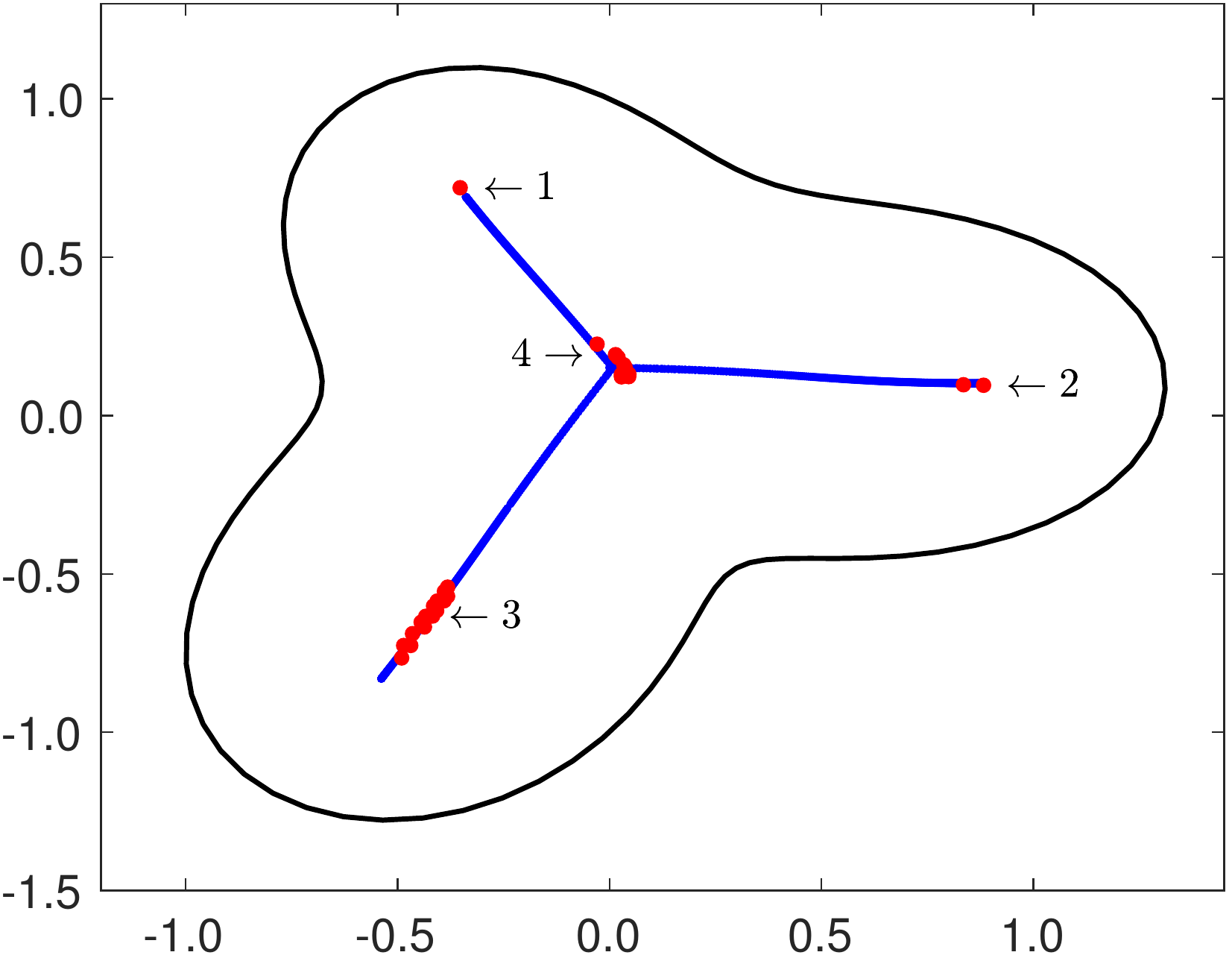}
\caption{Skeleton and touchdown points for the non-symmetric domain \eqref{eqn:nonparam}. The points marked $1-4$ correspond to the first touchdown location for solutions of \eqref{eq:intro} for values $\eps= 0.02,0.024,0.04,0.092$, respectively. \label{parsk}}
\end{figure}

\begin{figure}[!ht]
\subfigure[Mesh at $t = 0.010$.]{\includegraphics[width = 0.3\textwidth]{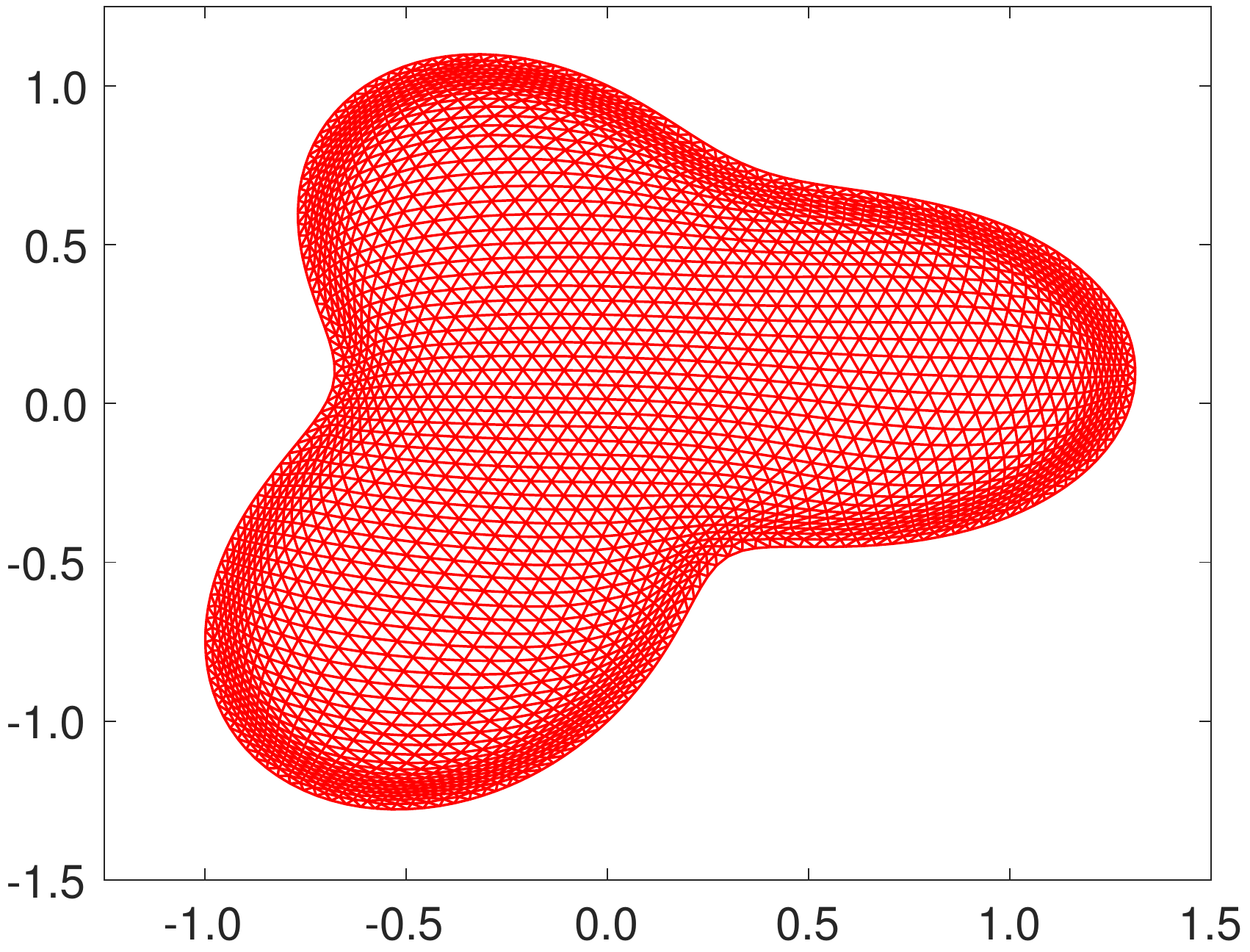}}\qquad
\subfigure[Mesh at $t = 0.232$.]{\includegraphics[width = 0.3\textwidth]{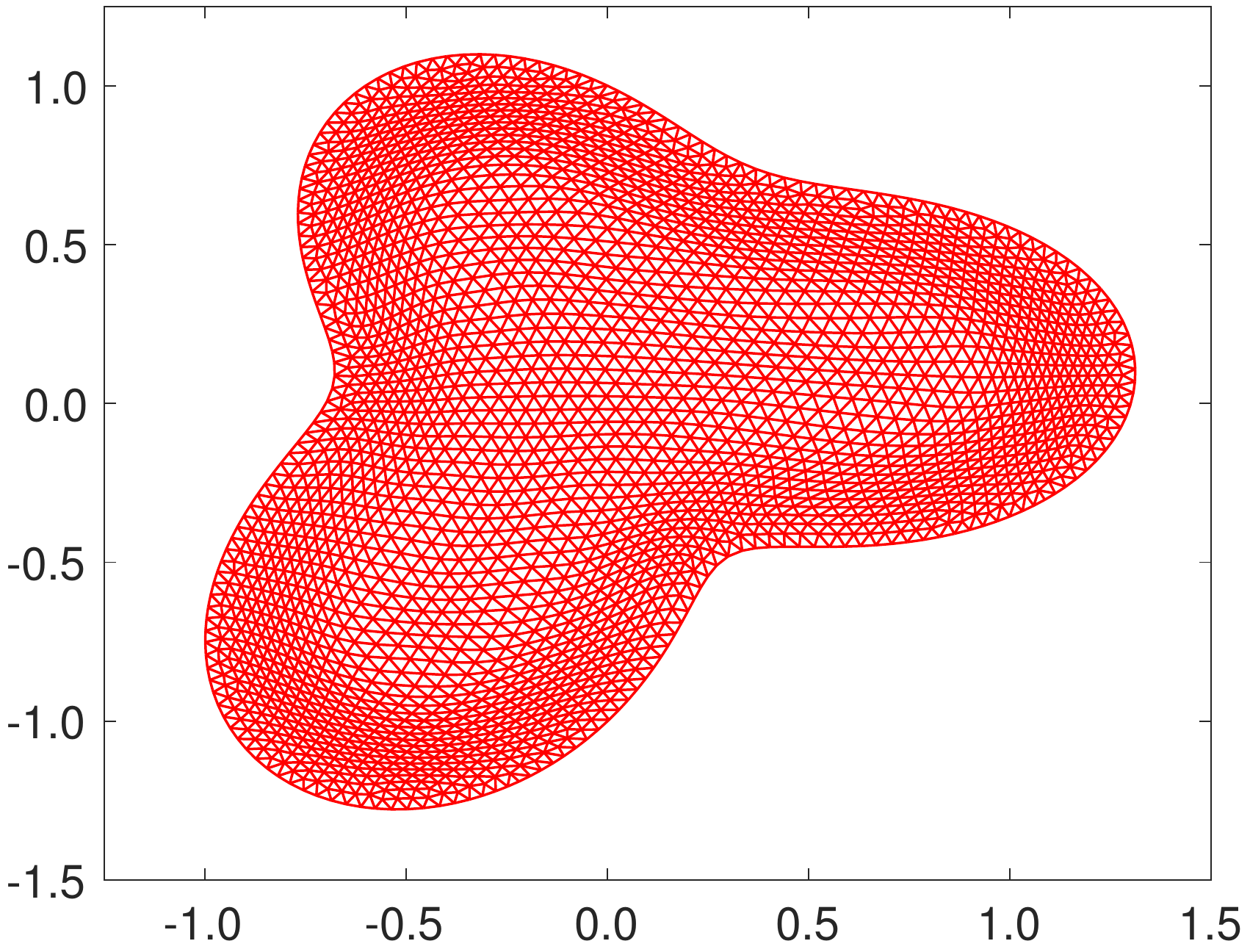}\label{partouchdowneps=0.02b}}\qquad
\subfigure[Mesh at $t = 0.310$.]{\includegraphics[width = 0.3\textwidth]{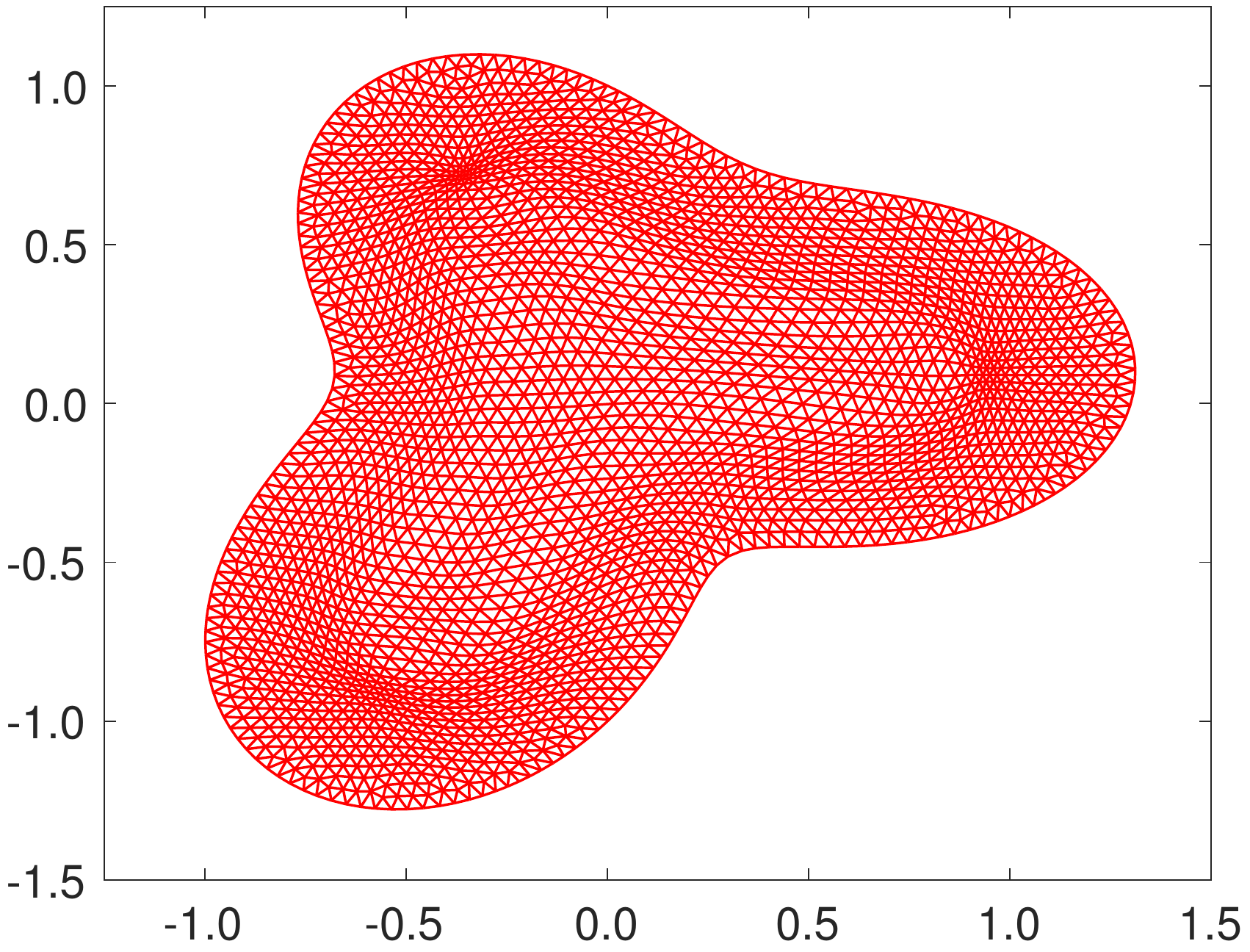}}\\
\subfigure[Profile at $t = 0.010$.]{\includegraphics[width = 0.3\textwidth]{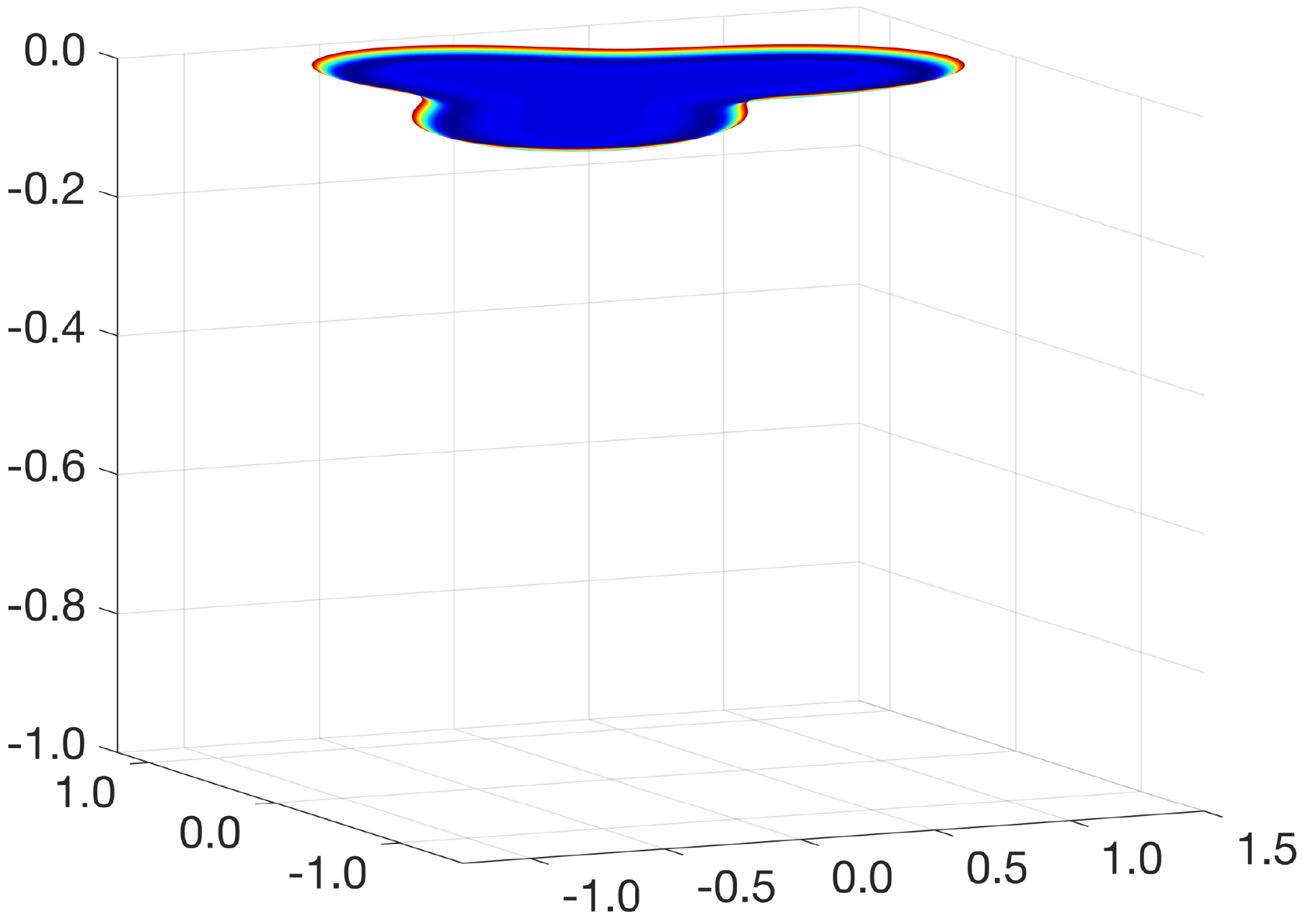}}\qquad
\subfigure[Profile at $t = 0.232$.]{\includegraphics[width = 0.3\textwidth]{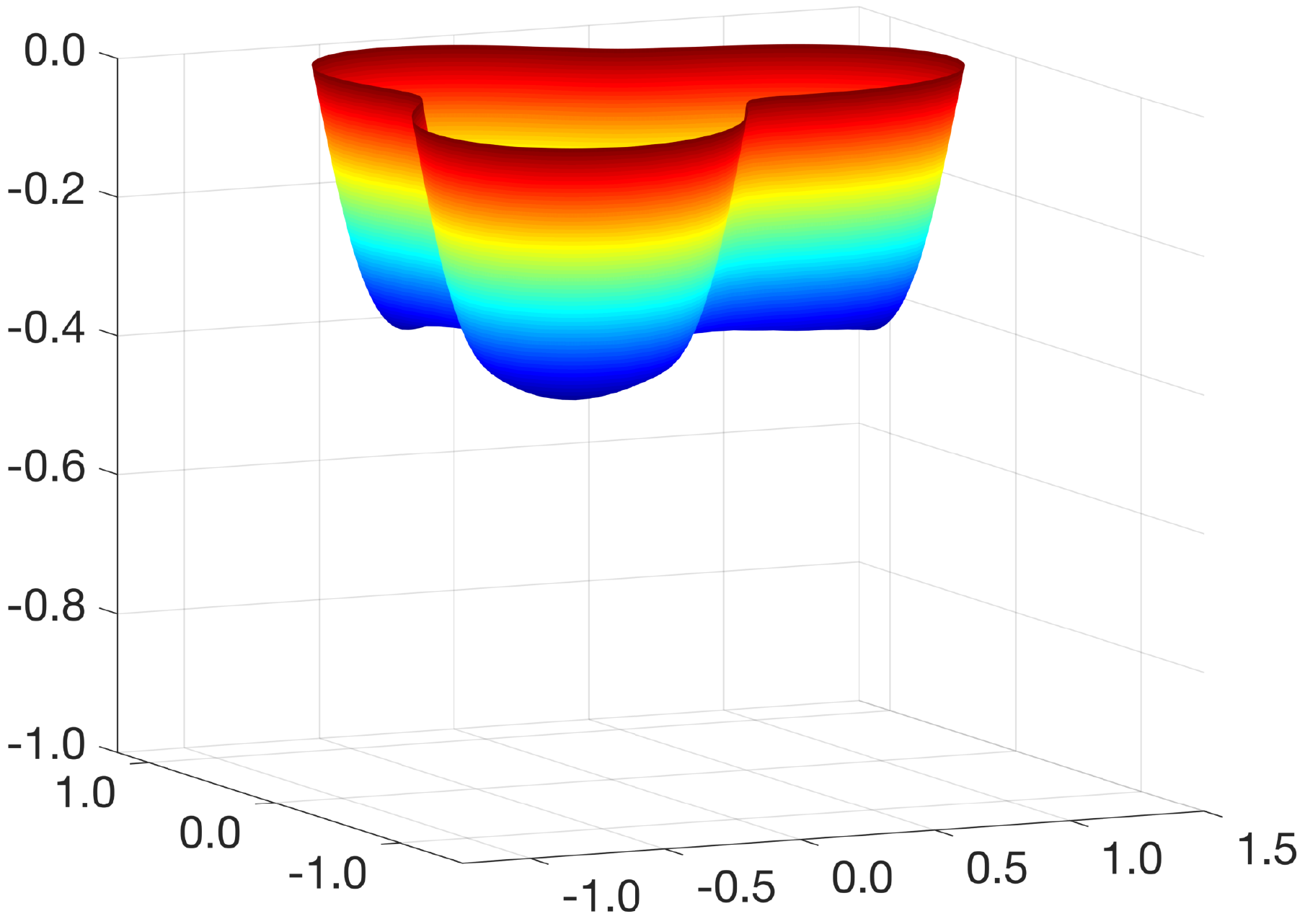}}\qquad
\subfigure[Profile at $t = 0.310$.]{\includegraphics[width = 0.3\textwidth]{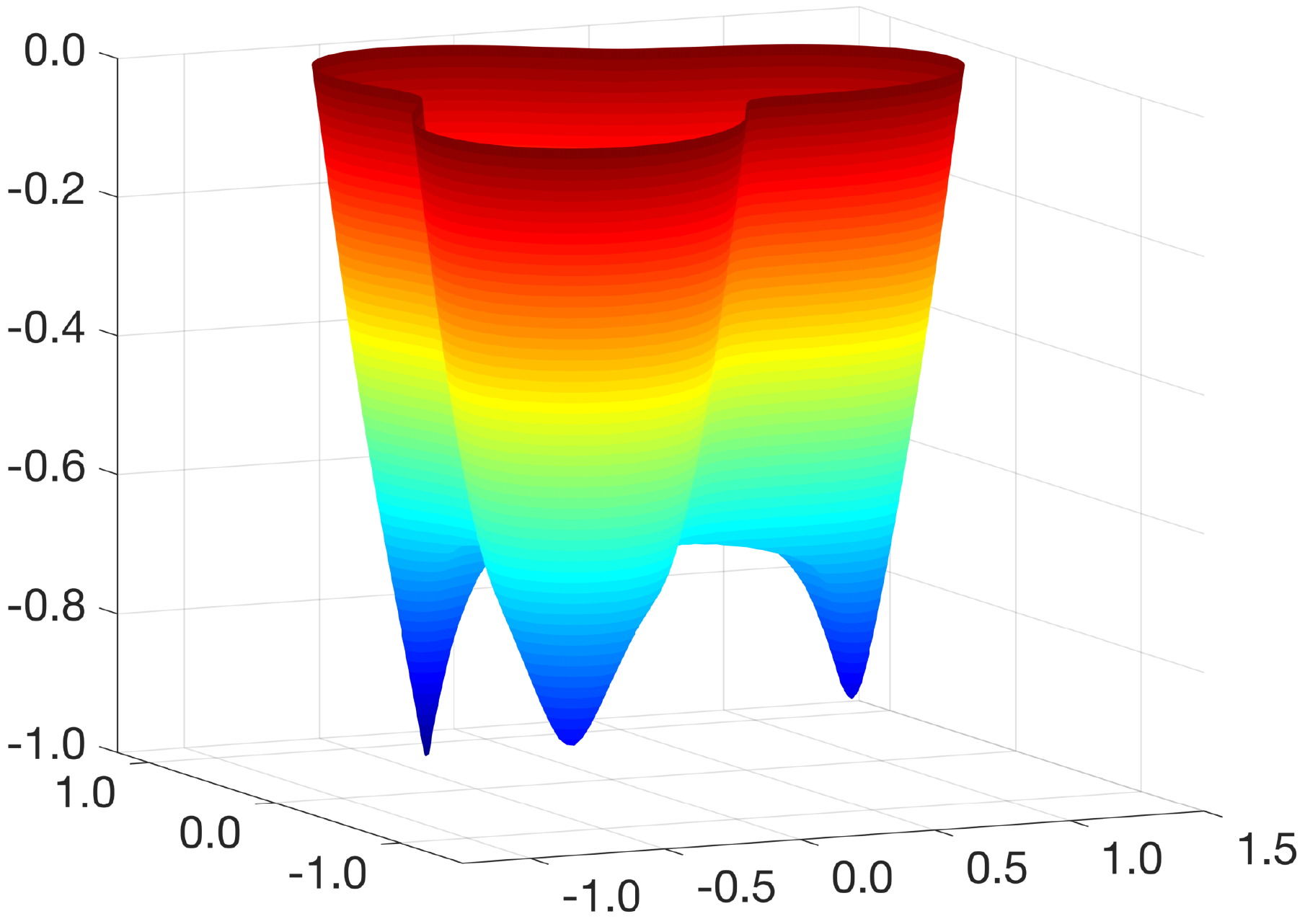}\label{partouchdowneps=0.02_f}}
\caption{Three snapshots of the evolution of the solution of \eqref{eq:intro} and the mesh for $\eps = 0.02$  (Mark 1 in Fig.~\ref{parsk}). The mesh size is $N = 5244$. \label{partouchdowneps=0.02}}
\end{figure}

\begin{figure}[!hb]
\subfigure[Mesh at $t = 0.009$.]{\includegraphics[width = 0.3\textwidth]{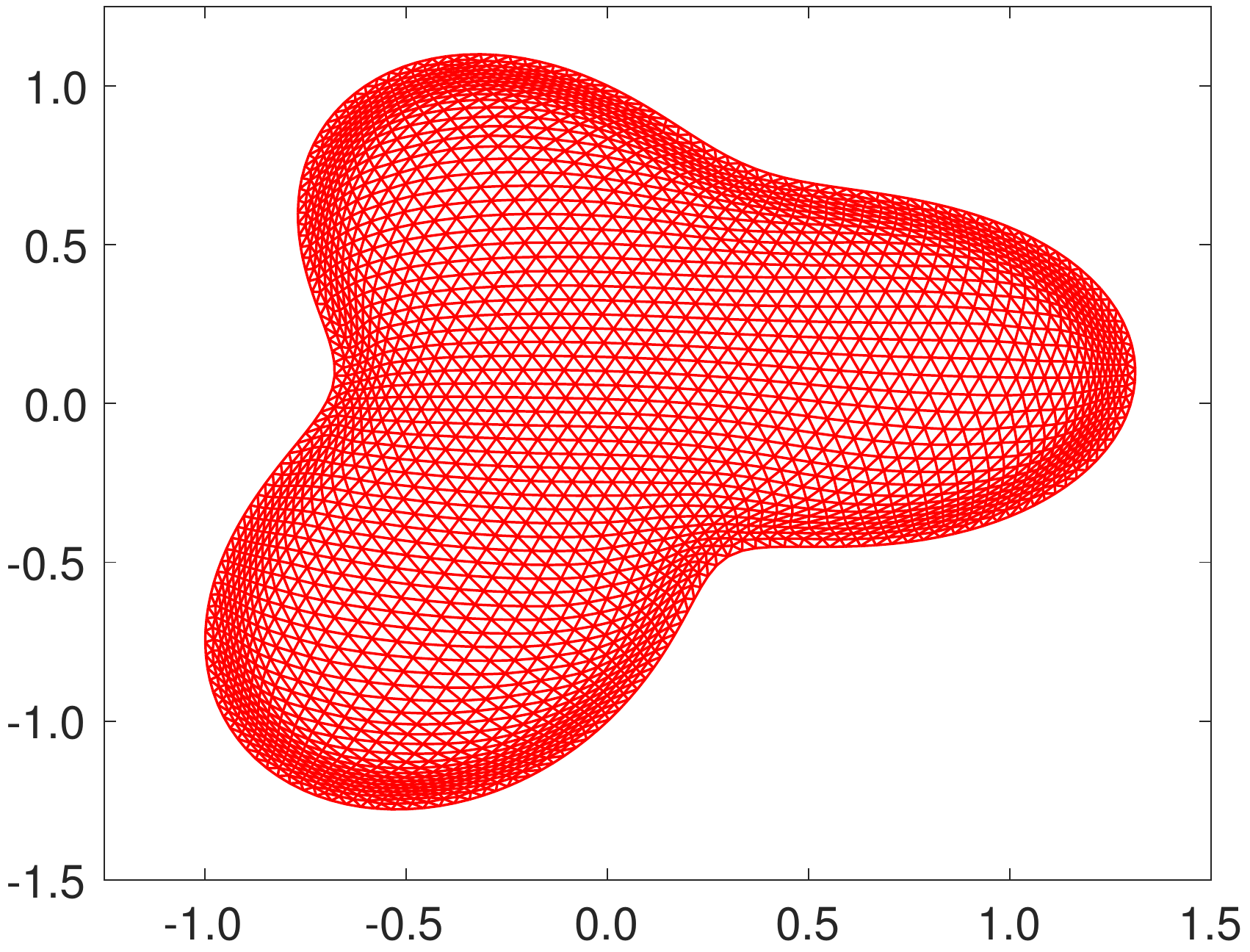}}\qquad
\subfigure[Mesh at $t = 0.190$.]{\includegraphics[width = 0.3\textwidth]{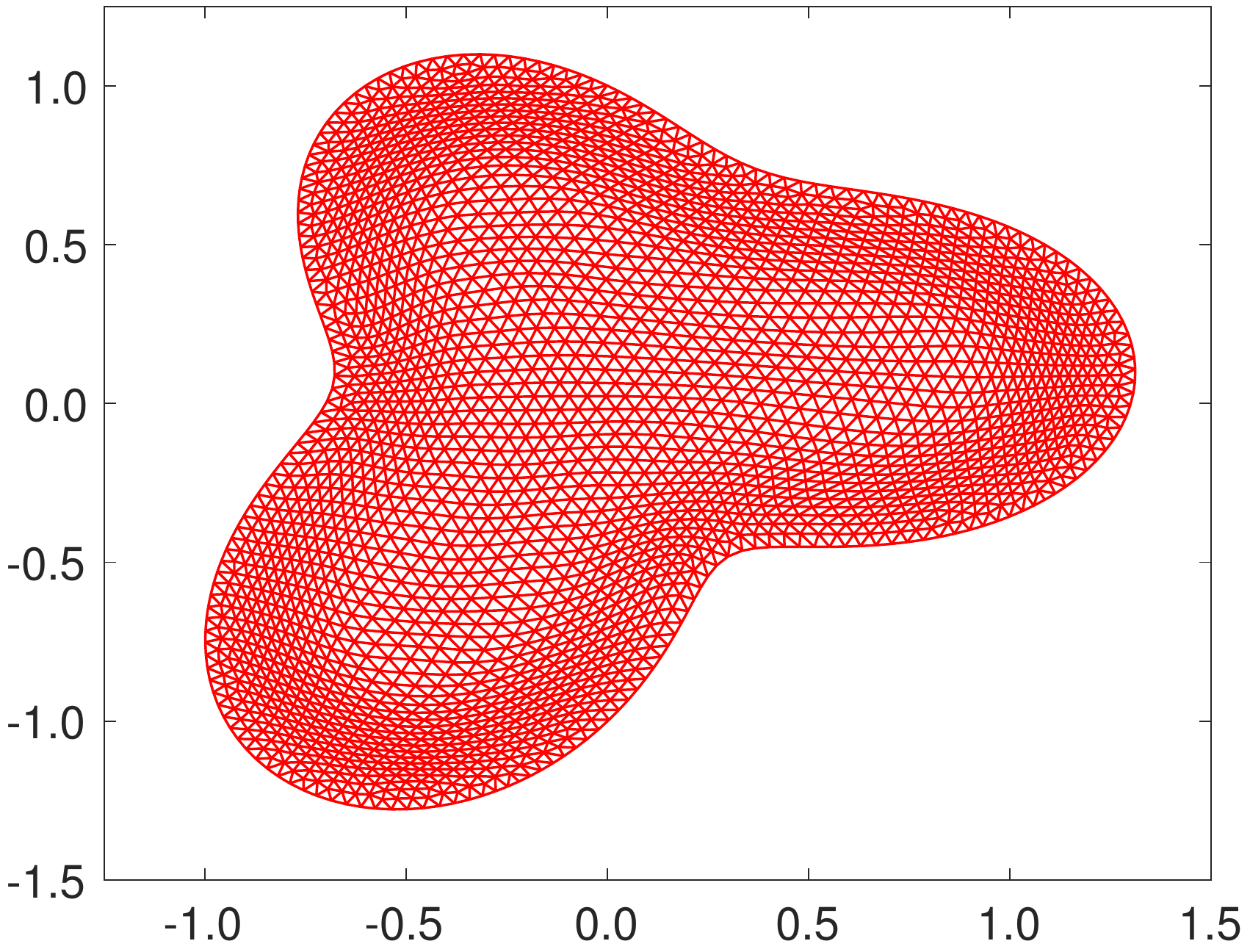}\label{partouchdowneps=0.024b}}\qquad
\subfigure[Mesh at $t = 0.309$.]{\includegraphics[width = 0.3\textwidth]{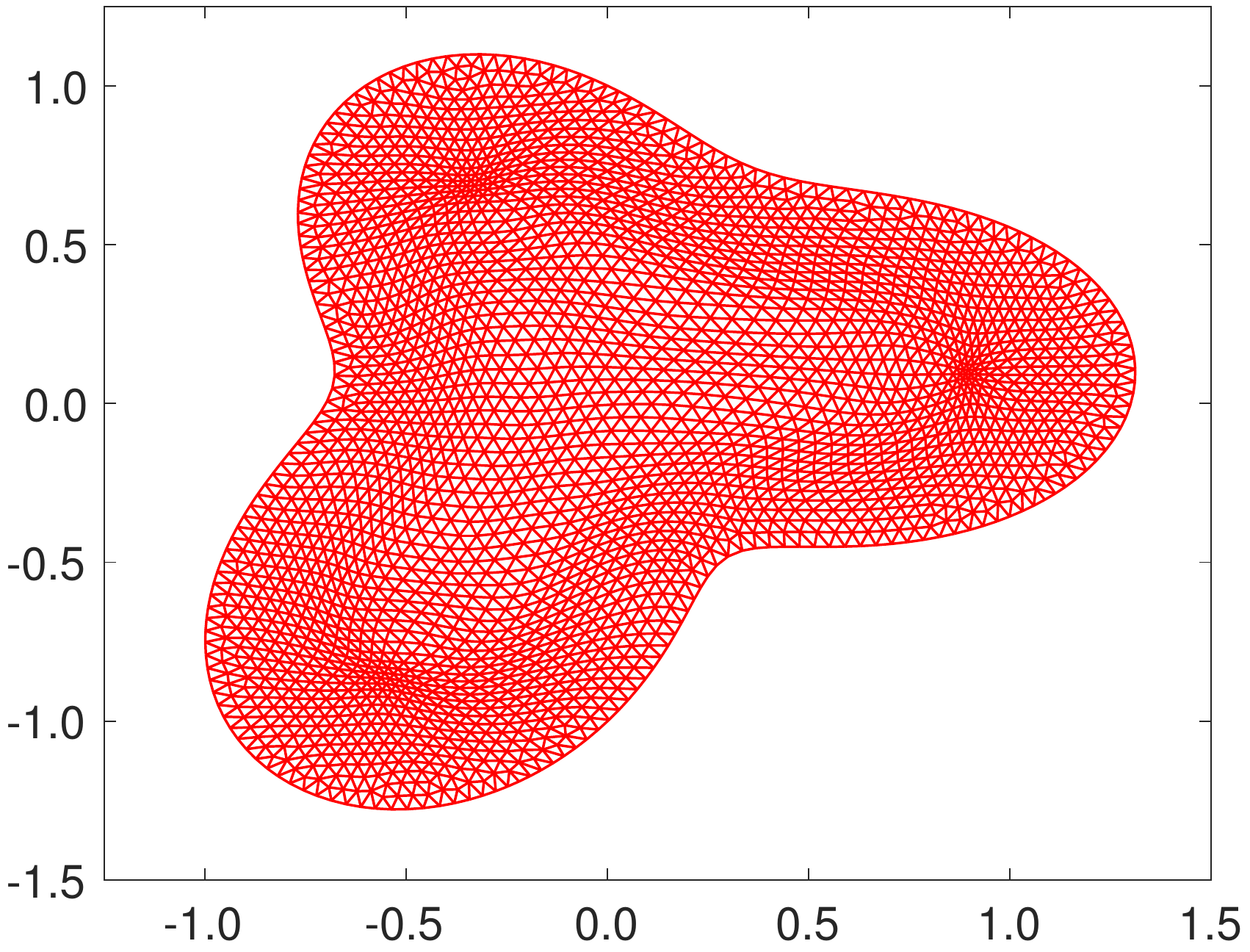}}\\
\subfigure[Profile at $t = 0.009$.]{\includegraphics[width = 0.3\textwidth]{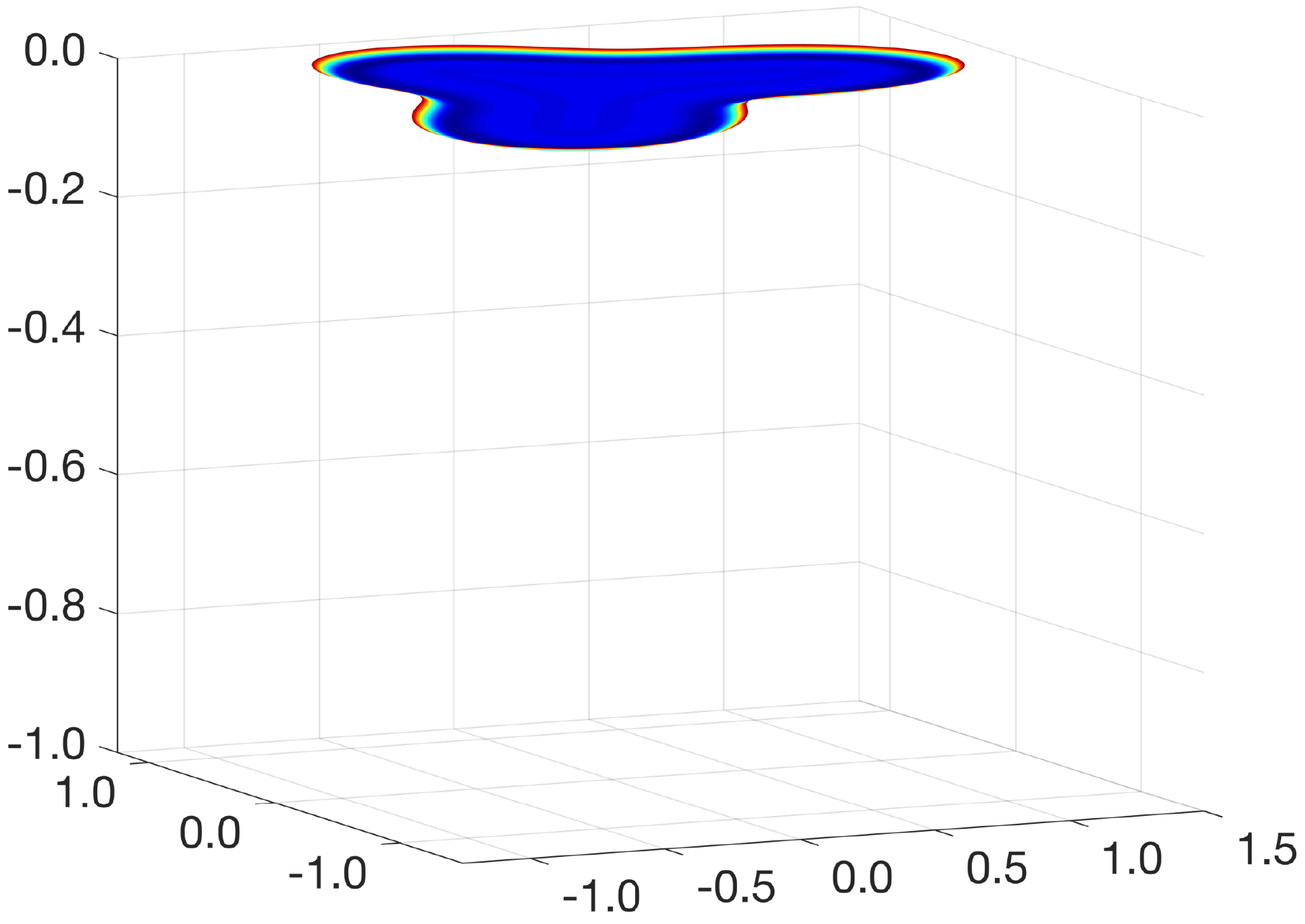}}\qquad
\subfigure[Profile at $t = 0.190$.]{\includegraphics[width = 0.3\textwidth]{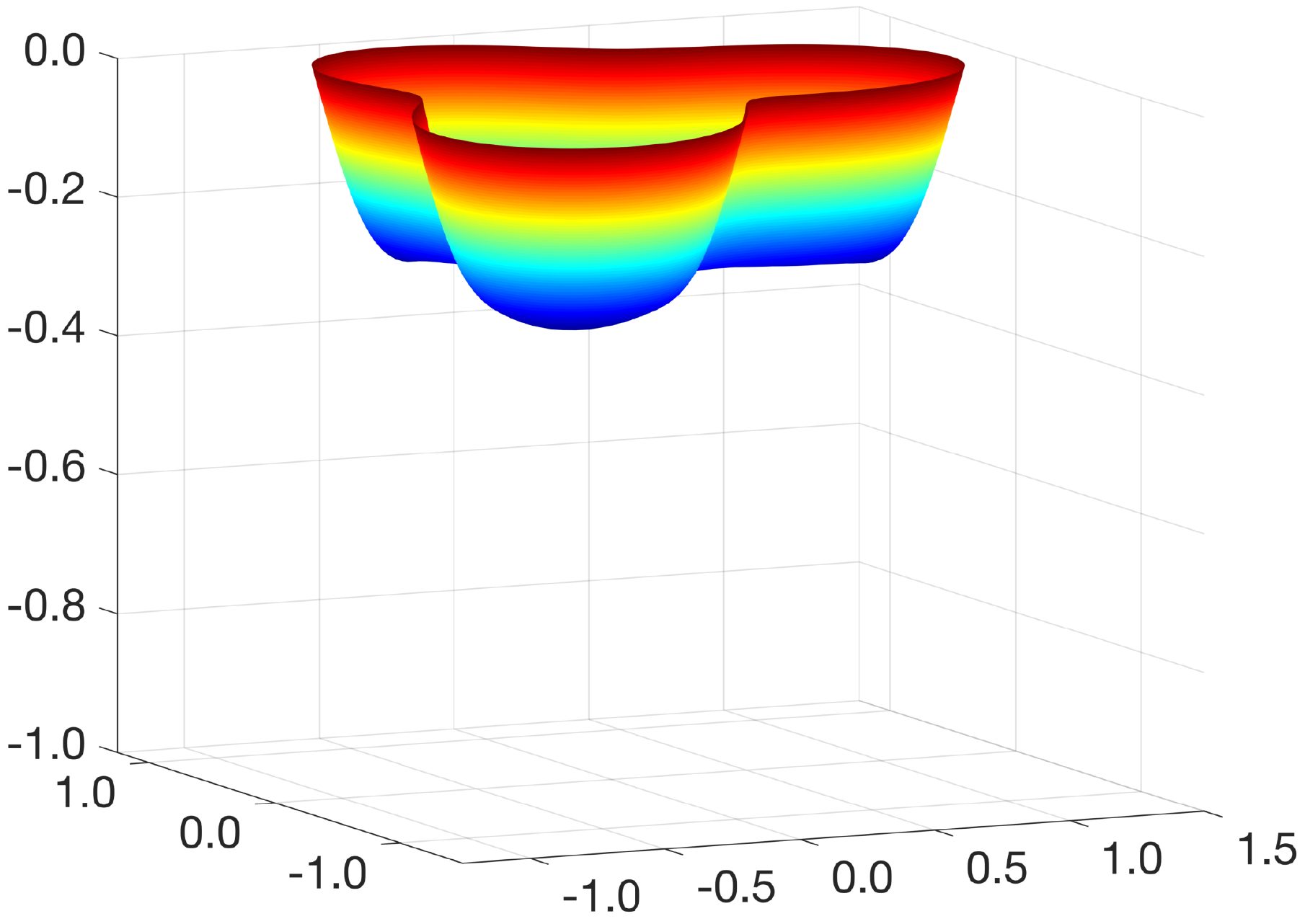}}\qquad
\subfigure[Profile at $t = 0.309$.]{\includegraphics[width = 0.3\textwidth]{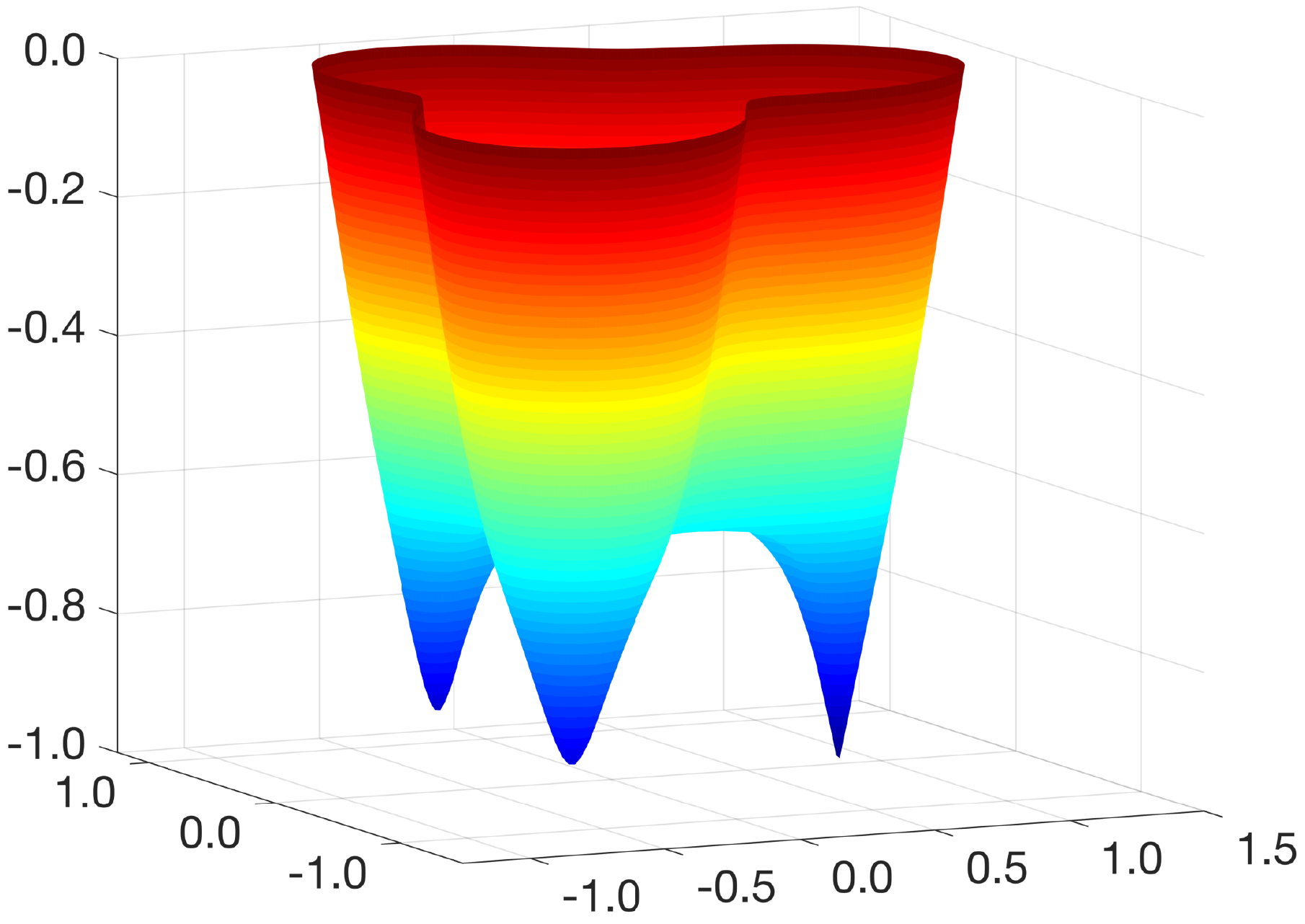} \label{partouchdowneps=0.024_f}}
\caption{Three snapshots of the evolution of the solution of \eqref{eq:intro} and the mesh for $\eps = 0.024$  (Mark 2 in Fig.~\ref{parsk}). The mesh size is $N = 5244$. \label{partouchdowneps=0.024}}
\end{figure}

\begin{figure}[!ht]
\centering
\subfigure[Mesh at $t = 0.009$. ]{\includegraphics[width = 0.3\textwidth]{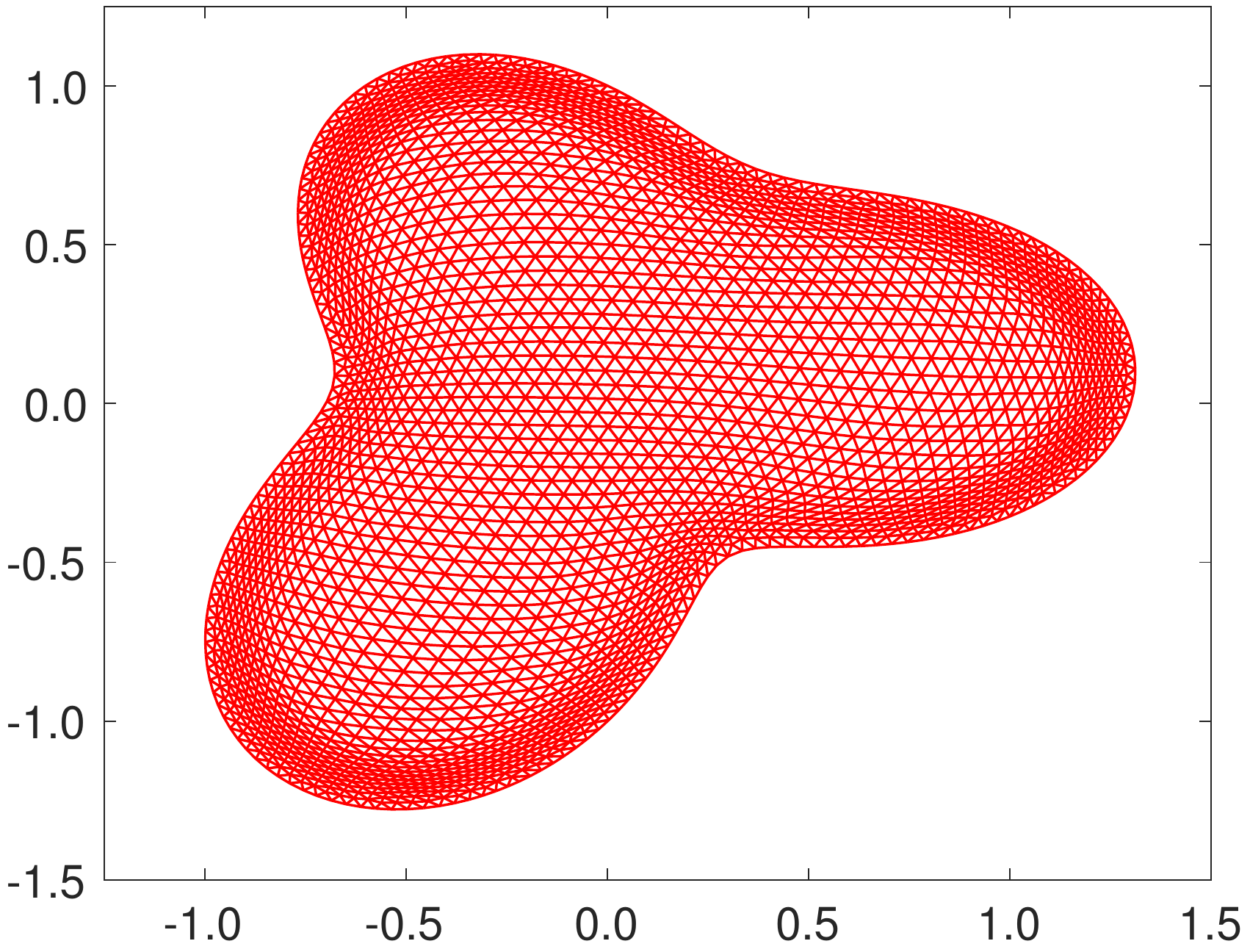}}\qquad
\subfigure[Mesh at $t = 0.280$.]{\includegraphics[width = 0.3\textwidth]{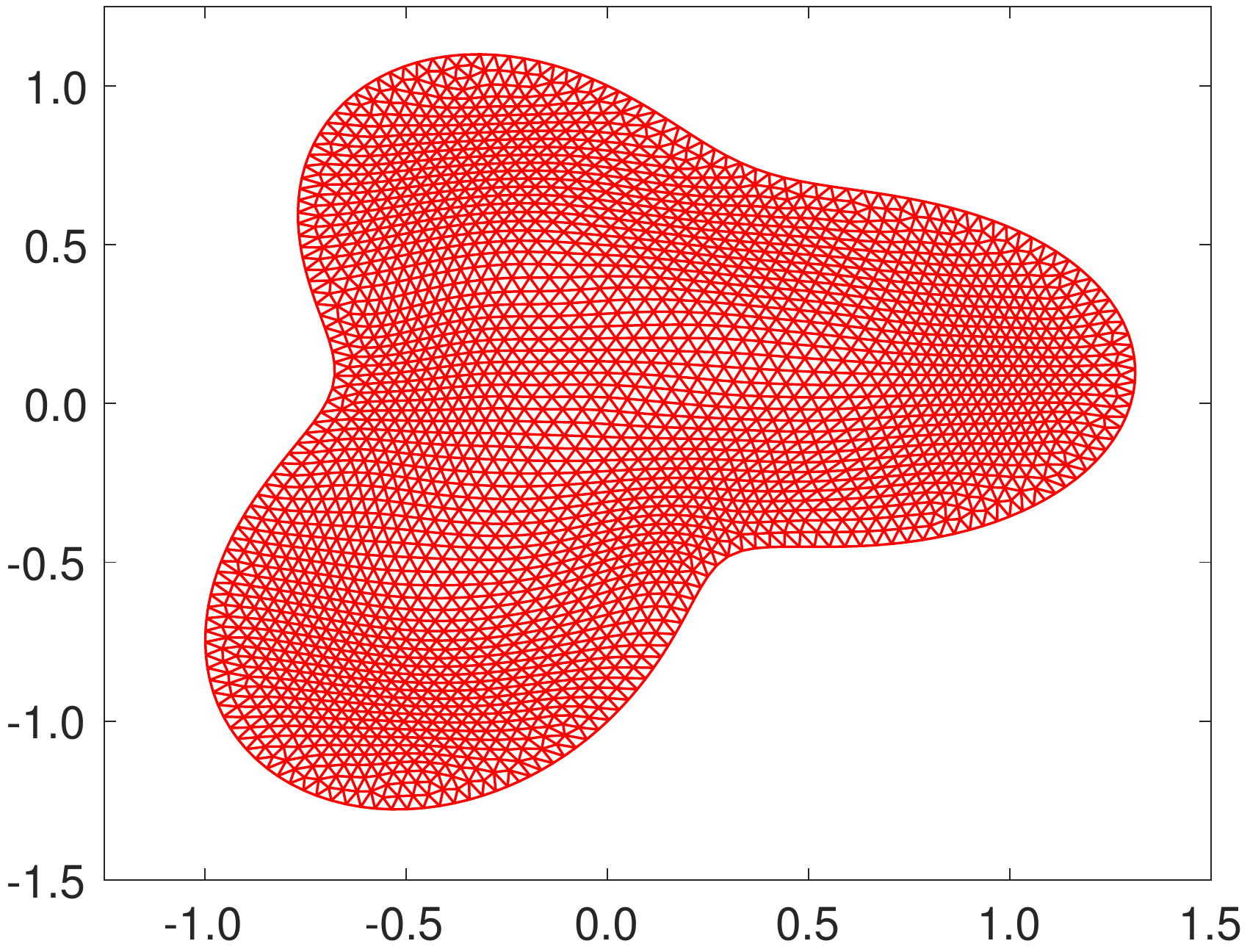}\label{partouchdowneps=0.04b}}\qquad
\subfigure[Mesh at $t = 0.303$.]{\includegraphics[width = 0.3\textwidth]{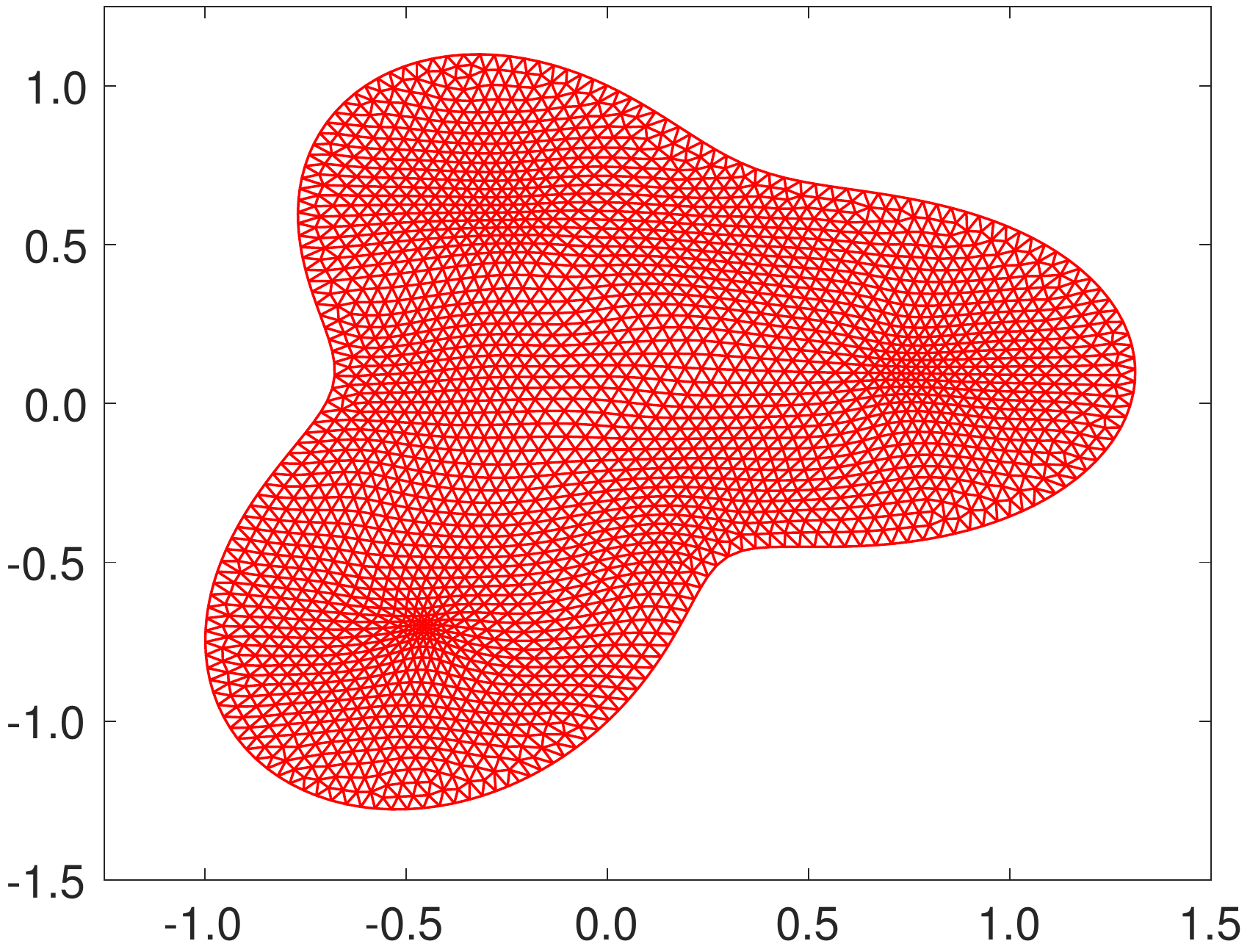}}\\
\subfigure[Solution at $t = 0.009$.]{\includegraphics[width = 0.3\textwidth]{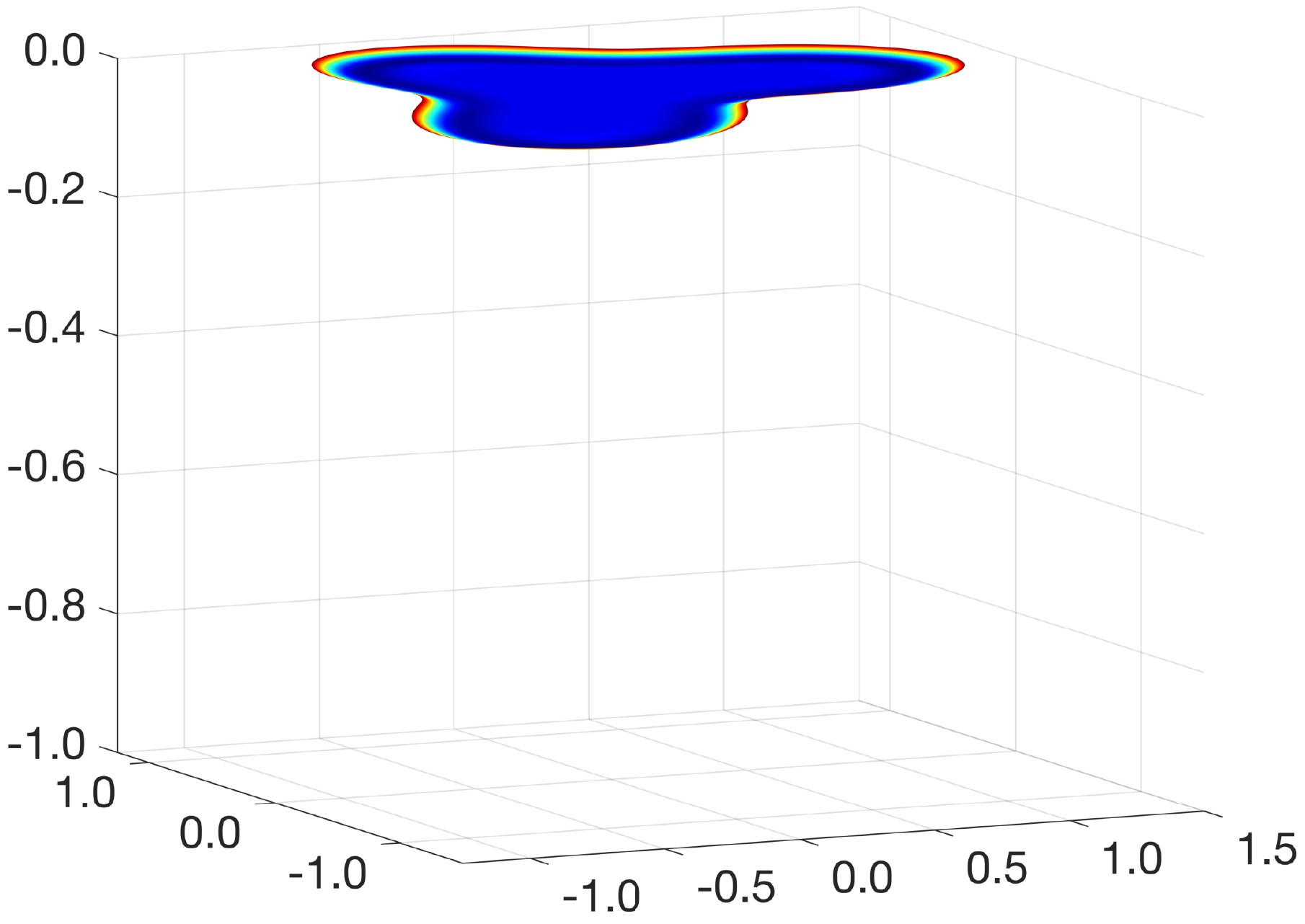}}\qquad
\subfigure[Solution at $t = 0.280$.]{\includegraphics[width = 0.3\textwidth]{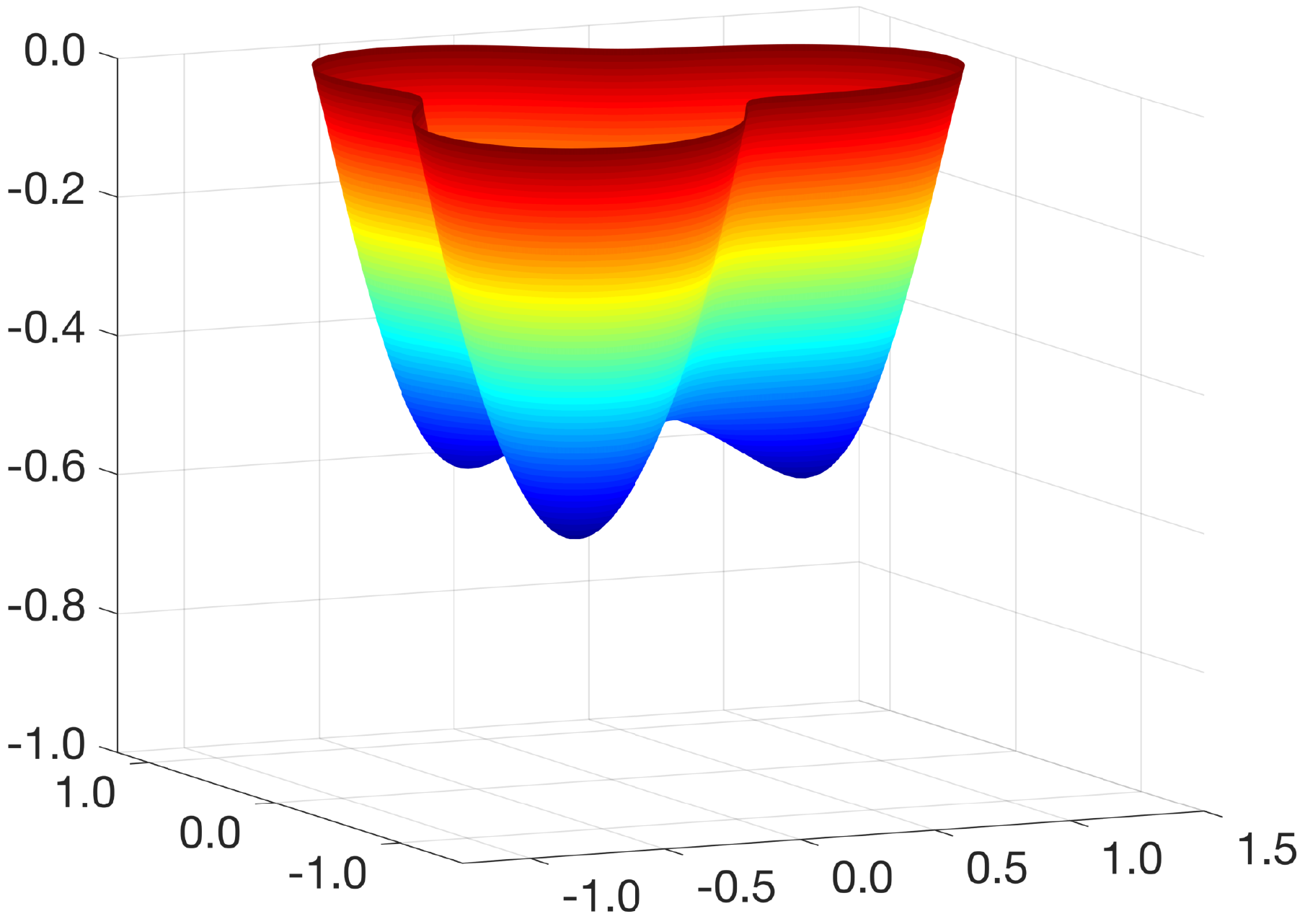}}\qquad
\subfigure[Solution at $t = 0.303$.]{\includegraphics[width = 0.3\textwidth]{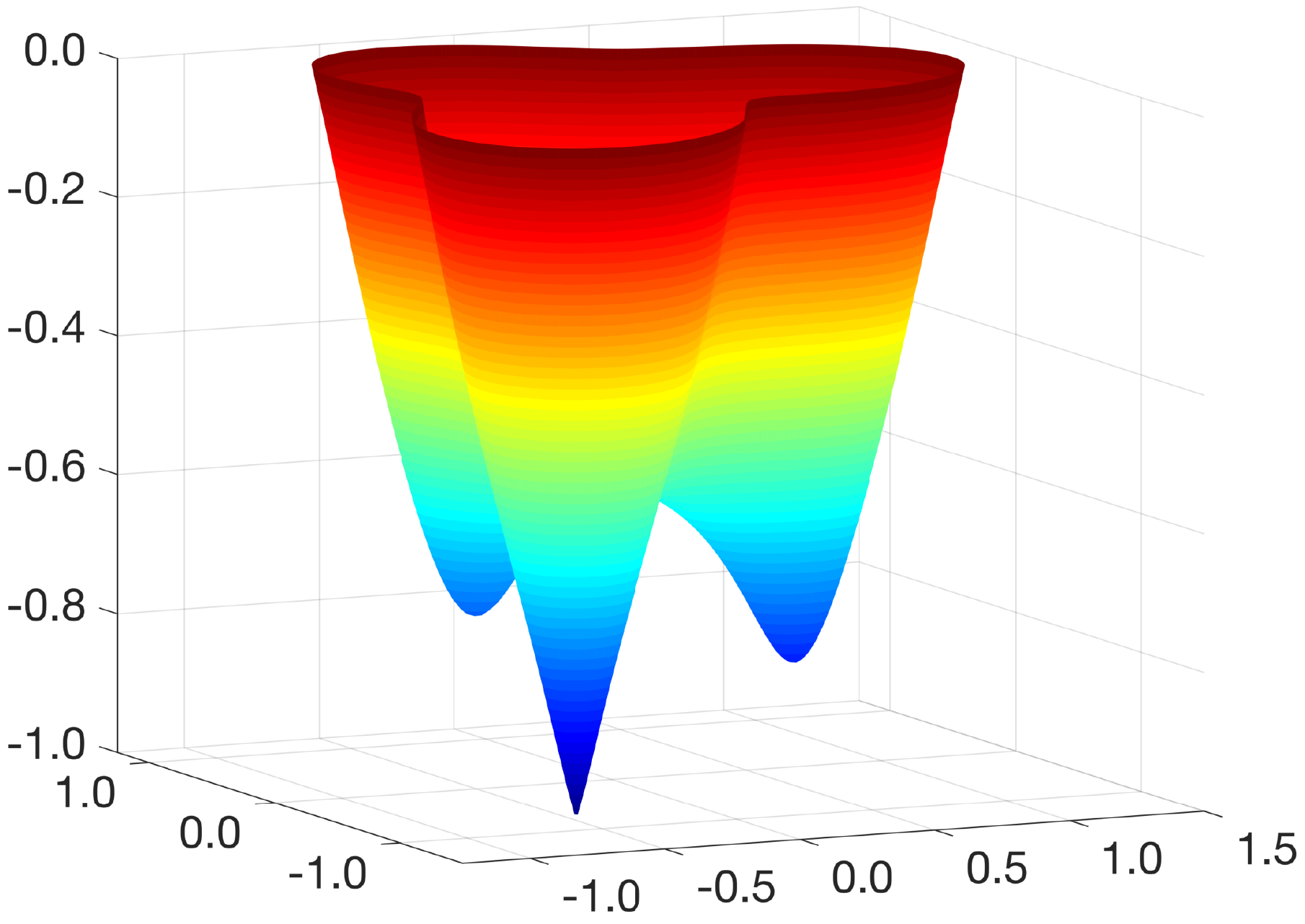}\label{partouchdowneps=0.04_f}}
\caption{The evolution of the solution of \eqref{eq:intro} and the mesh for $\eps = 0.04$. The mesh size is $N = 5244$.\label{partouchdowneps=0.04}}
\end{figure}

\begin{figure}[!hb]
\centering
\subfigure[Mesh at $t = 0.004$. ]{\includegraphics[width = 0.3\textwidth]{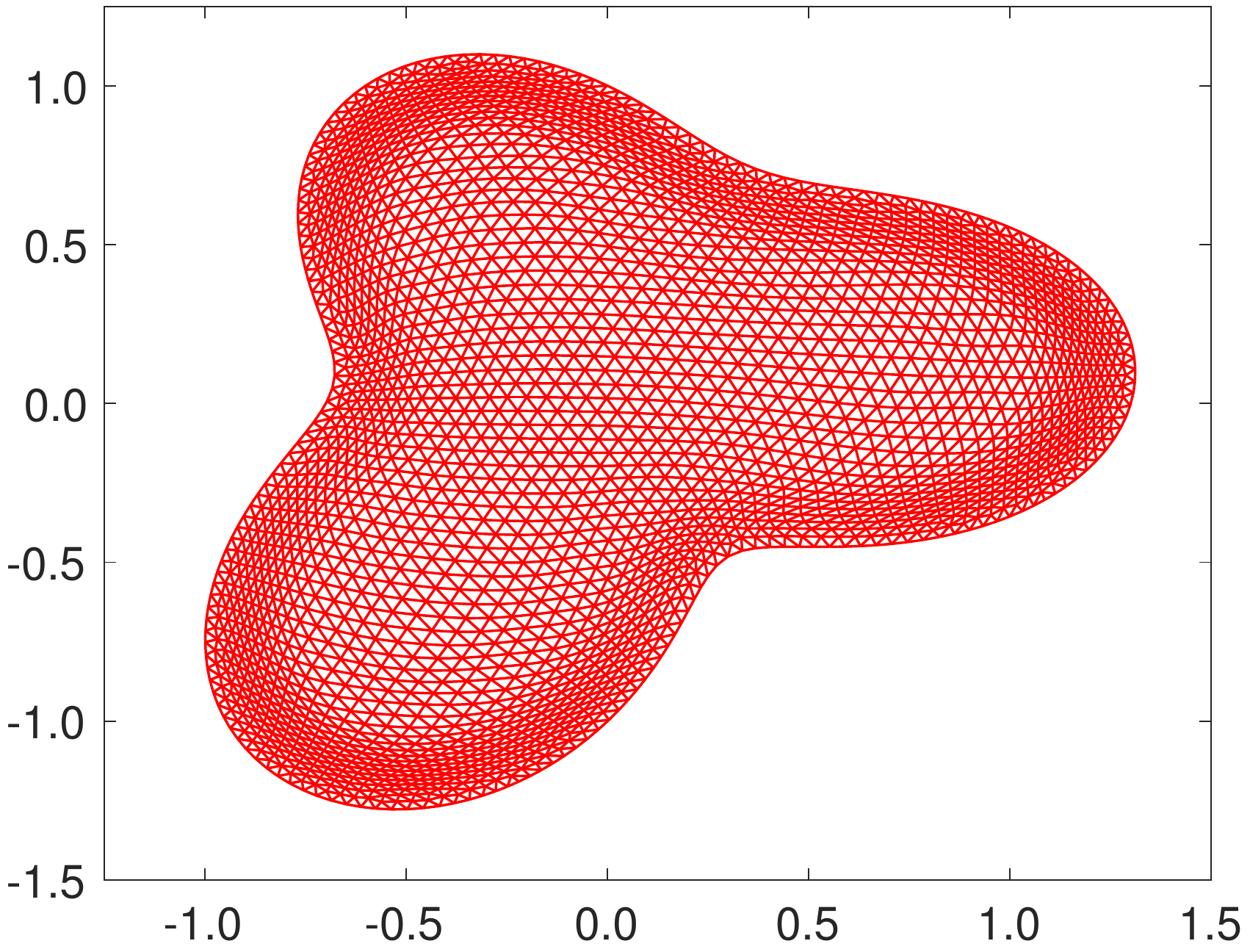}}\qquad
\subfigure[Mesh at $t = 0.250$.]{\includegraphics[width = 0.3\textwidth]{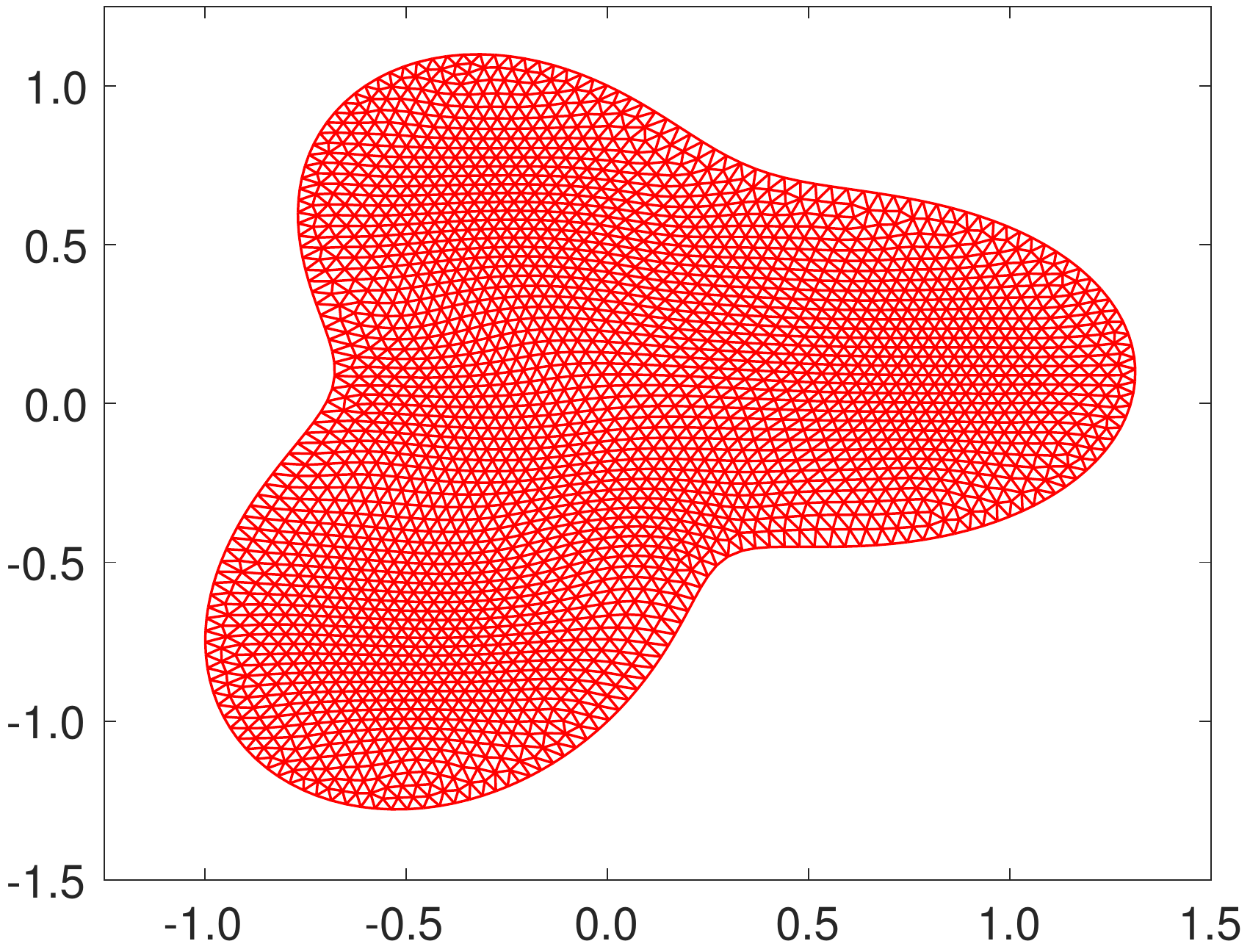}\label{partouchdowneps=0.092b}}\qquad
\subfigure[Mesh at $t = 0.303$.]{\includegraphics[width = 0.3\textwidth]{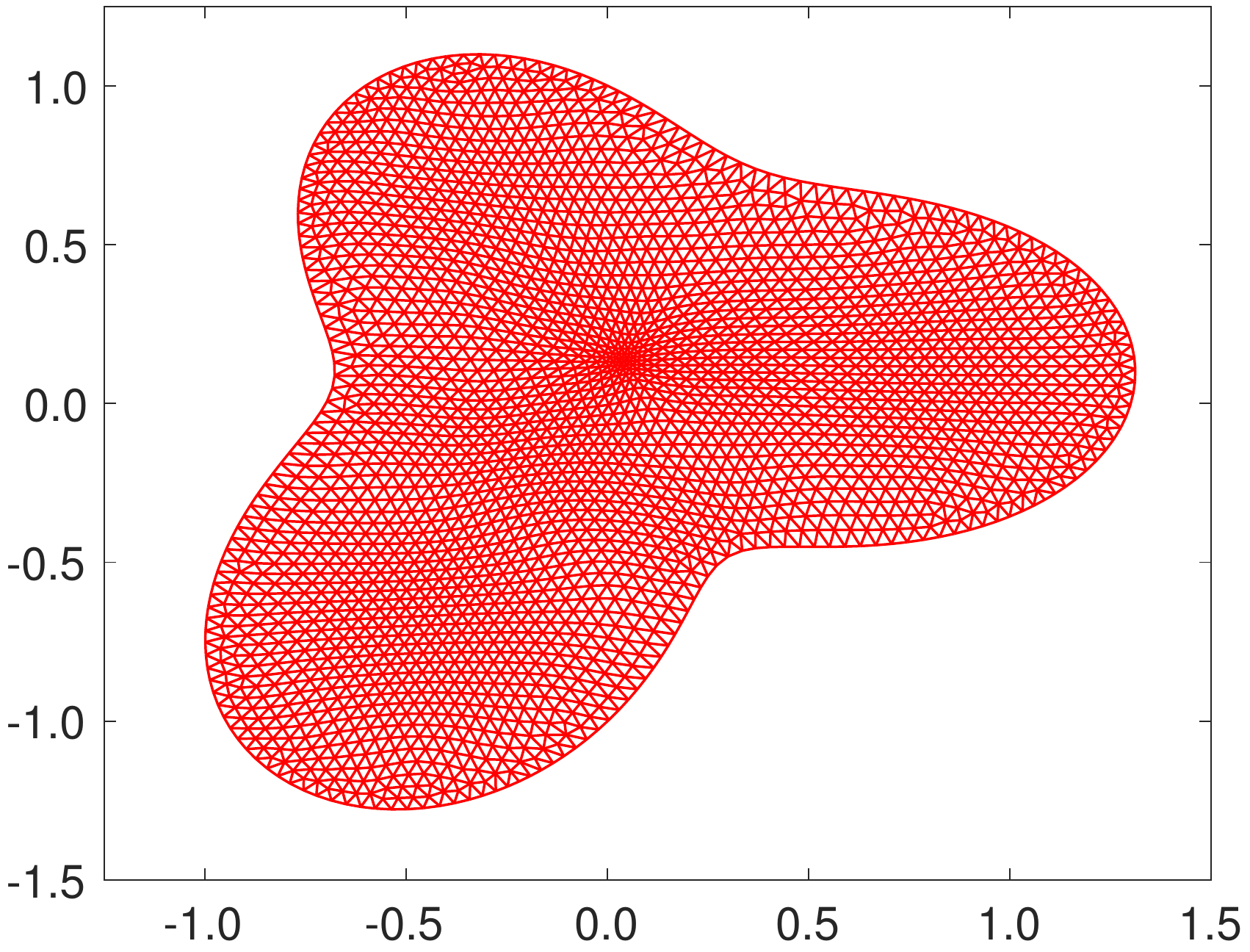}}\\
\subfigure[Solution at $t = 0.004$.]{\includegraphics[width = 0.3\textwidth]{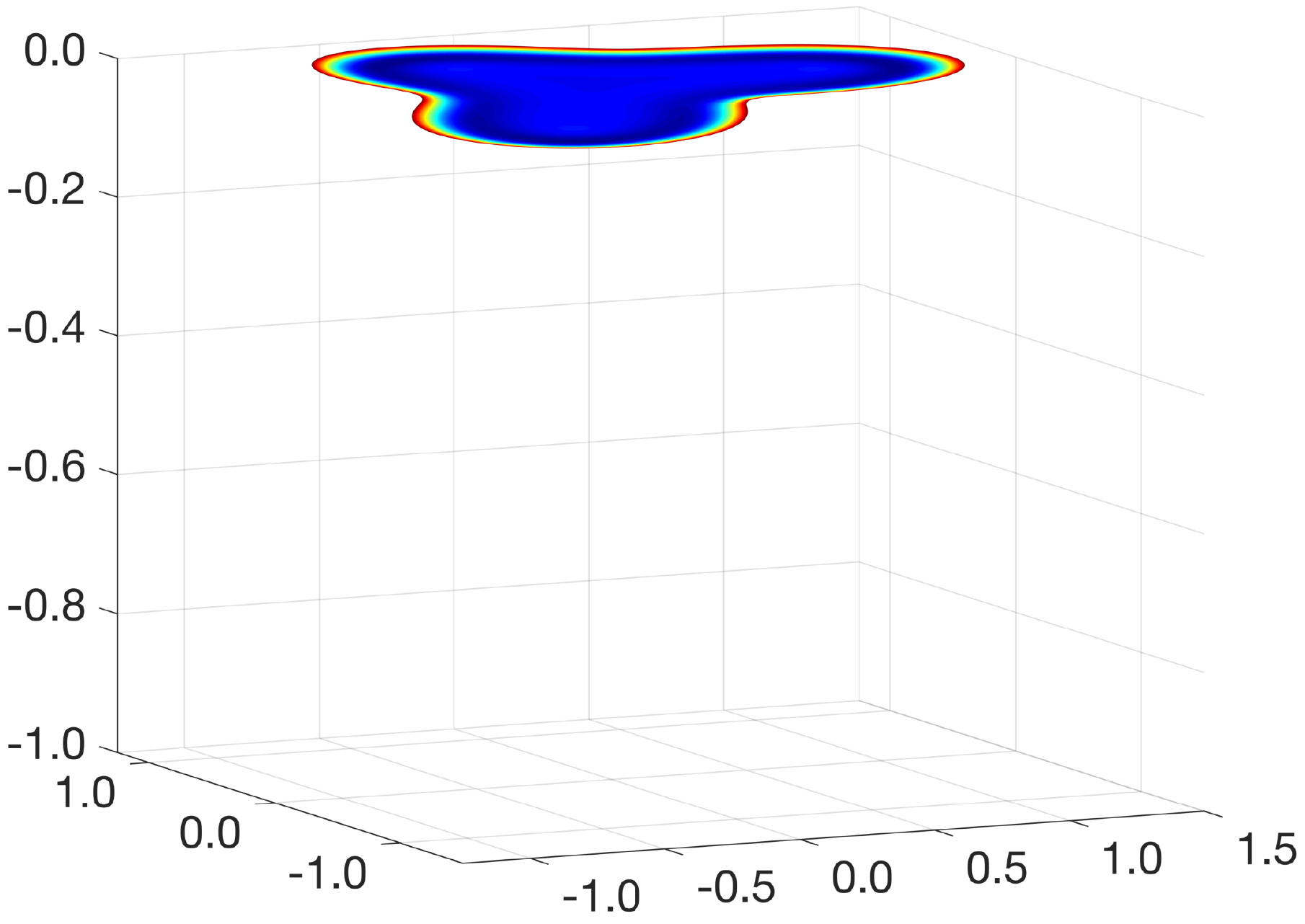}}\qquad
\subfigure[Solution at $t = 0.250$.]{\includegraphics[width = 0.3\textwidth]{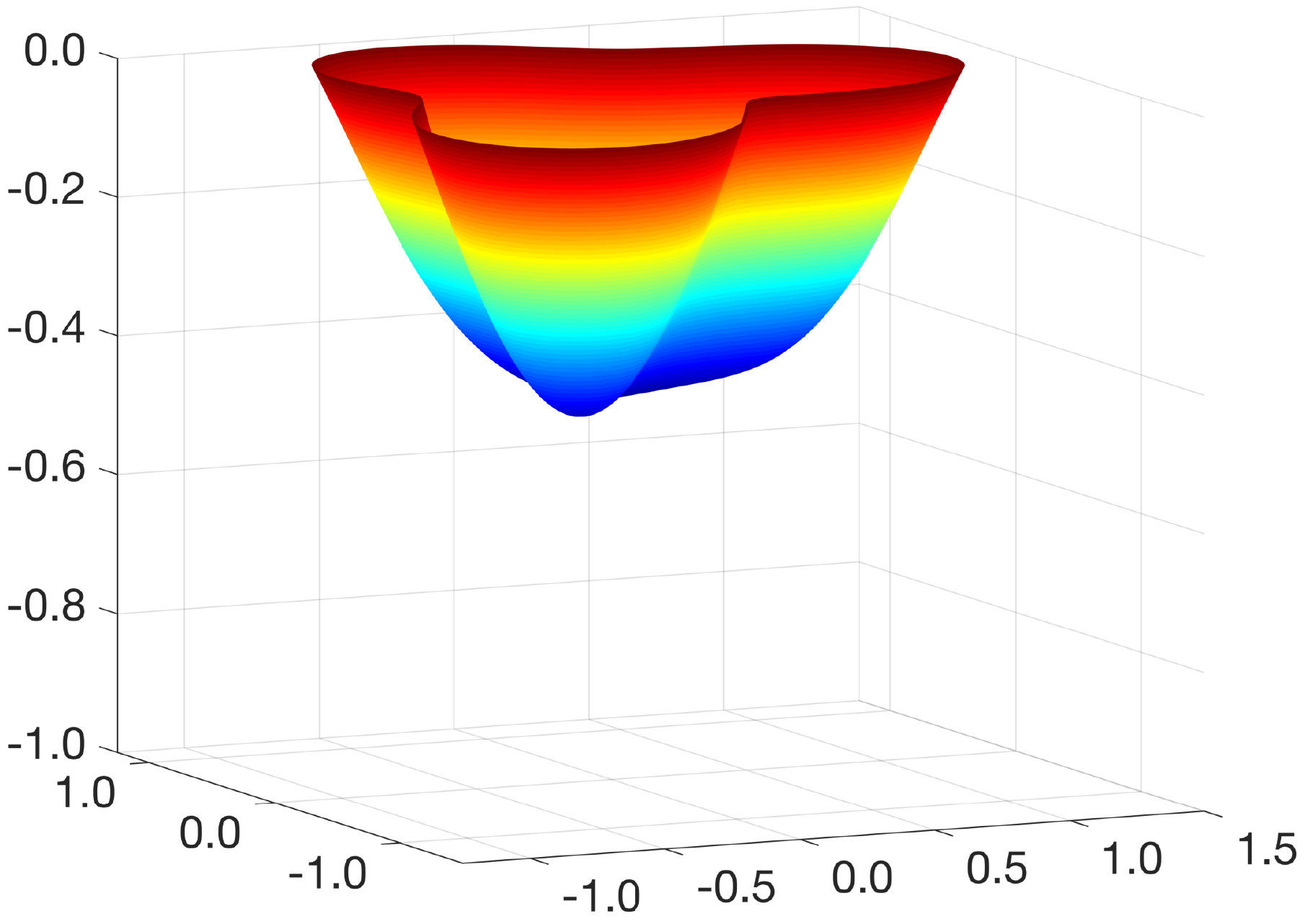}}\qquad
\subfigure[Solution at $t = 0.303$.]{\includegraphics[width = 0.3\textwidth]{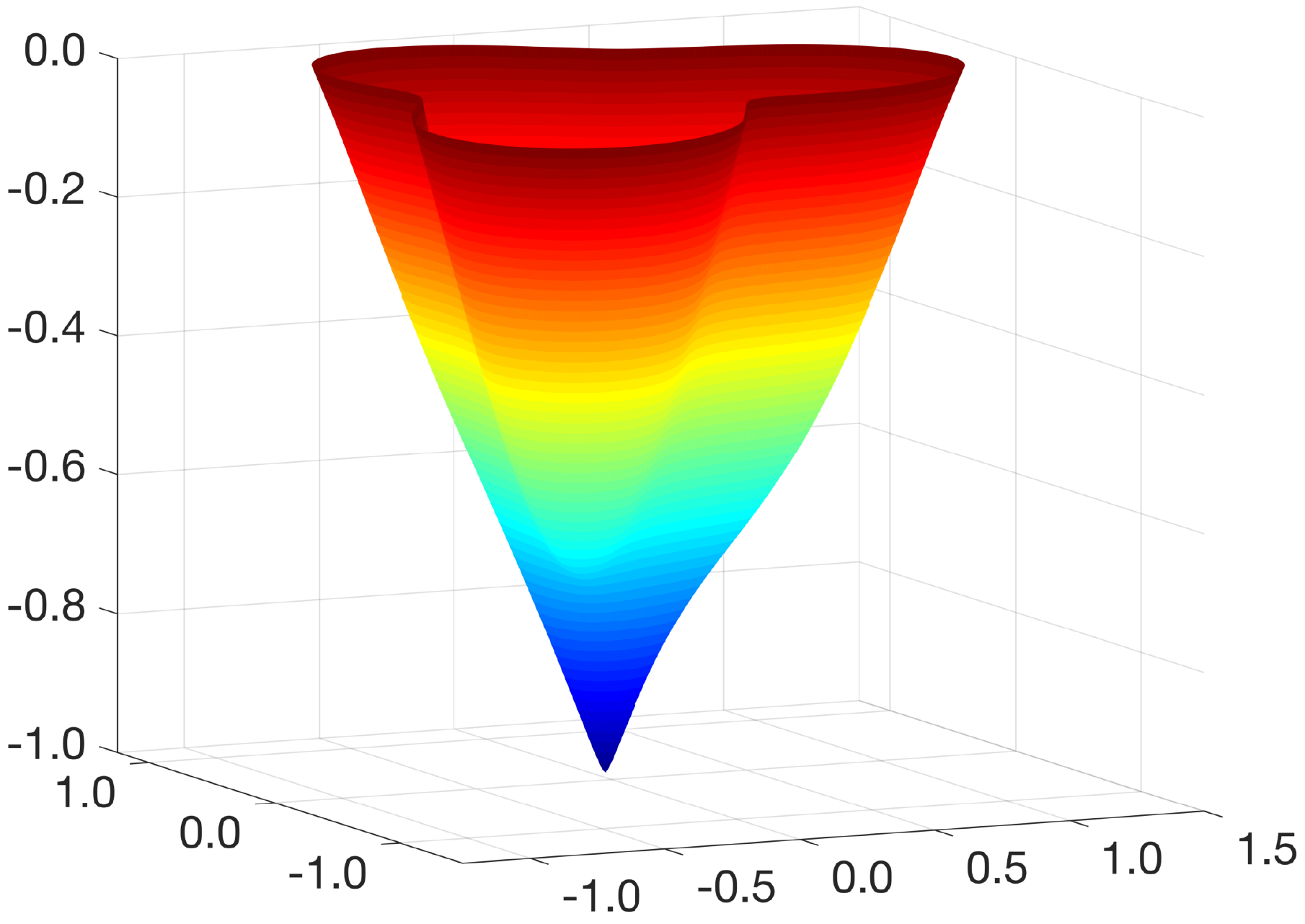}}
\caption{The evolution of the solution of \eqref{eq:intro} and the mesh for $\eps = 0.092$. The mesh size is $N = 5244$.\label{partouchdowneps=0.092}}
\end{figure}

In Fig.~\ref{parsk} the skeleton $\skel$ for the domain is displayed along with the first touchdown locations that arise from the parameter values $\eps \in (0.02,0.1)$. As $\eps$ increases, the singularity moves along each of the three branches of the skeleton before becoming fixed near the center of the domain. As with Example \ref{examrec}, we see that the track of the first touchdown location does not vary continuously with $\eps$. For the parameter values $\eps = 0.02\, (\text{Mark 1}), 0.024\, (\text{Mark 2}), 0.04\, (\text{Mark 3}), 0.092\, (\text{Mark 4})$ marked on Fig.~\ref{parsk}, we show snapshots of the evolution of the solution in Figs.~\ref{partouchdowneps=0.02}-\ref{partouchdowneps=0.092}, respectively.

\begin{figure}[!hb]
\centering
\subfigure[Profile at touchdown, $\eps= 0.02$. Mesh size is $N = 5244$. ]{\includegraphics[width = 0.3\textwidth]{ProfileNewPar_eps_2e_02no3.pdf}}\qquad
\subfigure[Profile at touchdown, $\eps= 0.024$. Mesh size is $N = 5244$.]{\includegraphics[width = 0.3\textwidth]{ProfileNewPar_eps_24e_3no3.pdf}}\qquad
\subfigure[Profile at touchdown, $\eps= 0.04$. Mesh size is $N = 5244$.]{\includegraphics[width = 0.3\textwidth]{ProfileNewPar_eps_4e_02no3.pdf}}\\
\subfigure[Profile at touchdown, $\eps= 0.02$. Mesh size is $N = 11955$.]{\includegraphics[width = 0.3\textwidth]{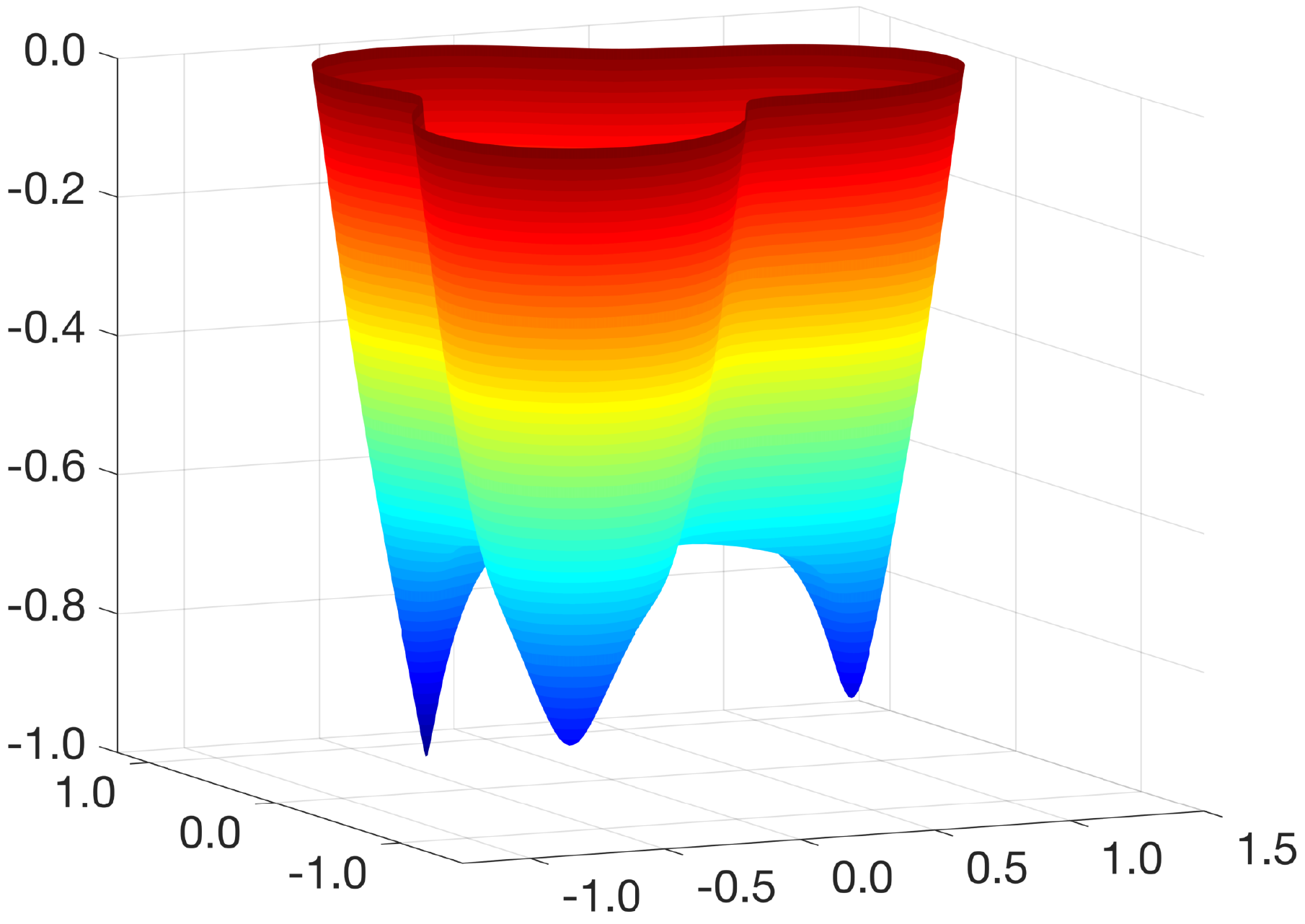}}\qquad
\subfigure[Profile at touchdown, $\eps= 0.024$. Mesh size is $N = 11955$.]{\includegraphics[width = 0.3\textwidth]{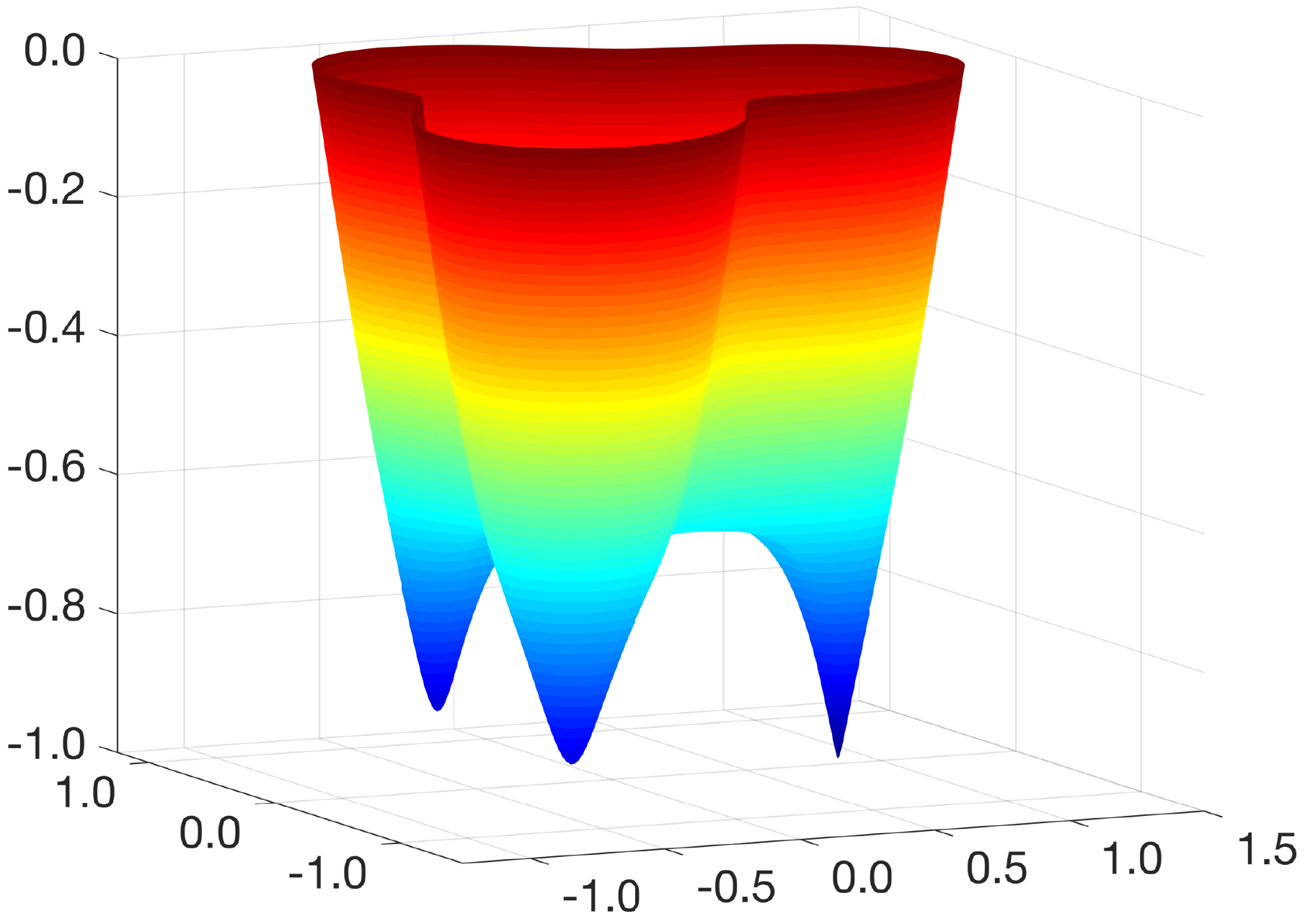}}\qquad
\subfigure[Profile at touchdown, $\eps= 0.04$. Mesh size is $N = 11955$.]{\includegraphics[width = 0.3\textwidth]{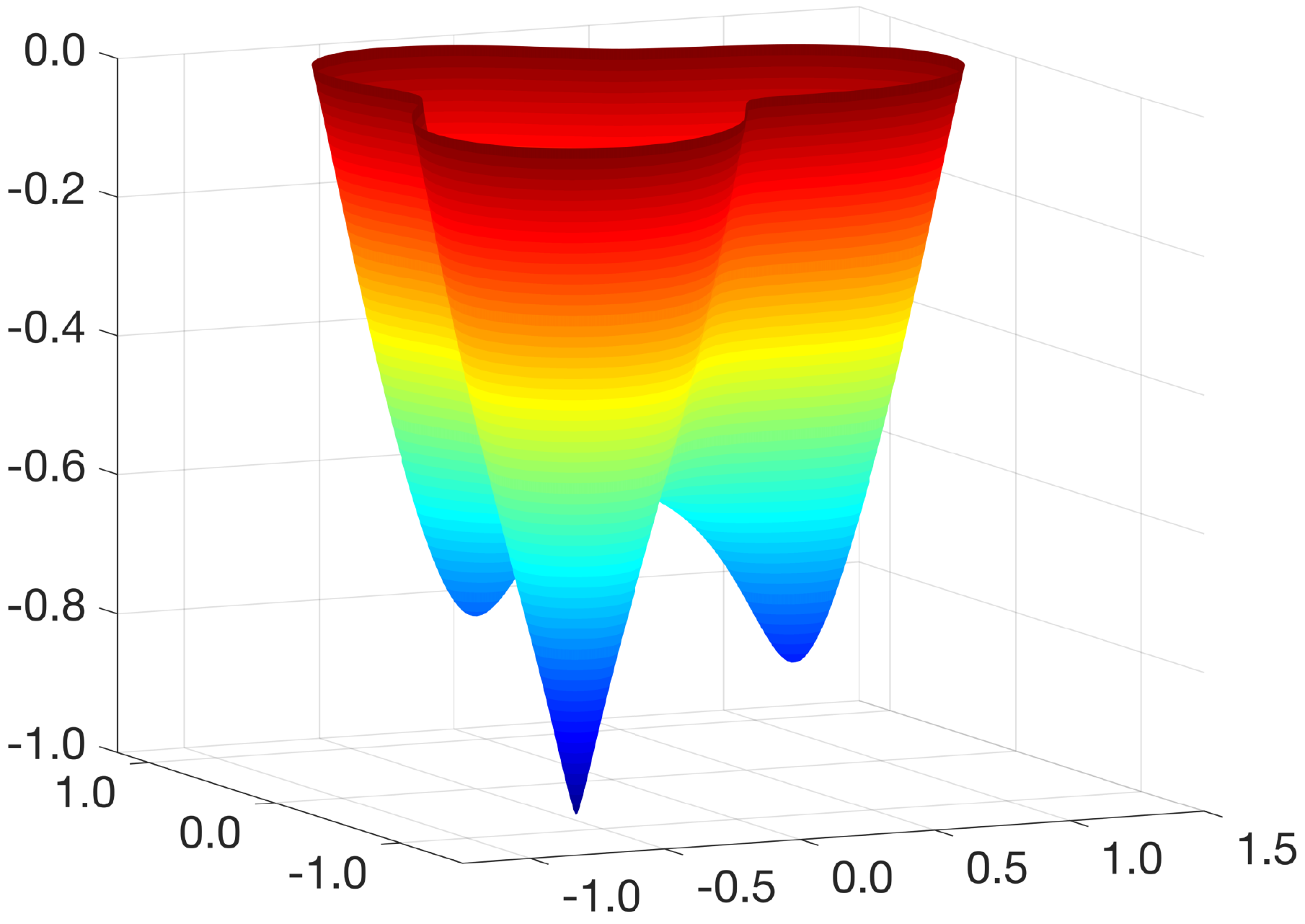}}
\caption{ The top row shows profiles pictures of $\eps = 0.02$, $\eps = 0.024$, $\eps = 0.04$ with mesh size $N = 5244$. The associated profiles obtained by mesh size $N = 11955$ are presented in the second row.\label{par12touchdownProfile}}
\end{figure}

\begin{figure}[!hb]
\centering
\subfigure[Mesh at touchdown, $\eps= 0.02$. Mesh size is $N = 5244$. ]{\includegraphics[width = 0.3\textwidth]{MeshNewPar_eps_2e_02no3.pdf}}\qquad
\subfigure[Mesh at touchdown, $\eps= 0.024$. Mesh size is $N = 5244$.]{\includegraphics[width = 0.3\textwidth]{MeshNewPar_eps_24e_3no3.pdf}}\qquad
\subfigure[Mesh at touchdown, $\eps= 0.04$. Mesh size is $N = 5244$.]{\includegraphics[width = 0.3\textwidth]{MeshNewPar_eps_4e_02no3.pdf}}\\
\subfigure[Mesh at touchdown, $\eps= 0.02$. Mesh size is $N = 11955$.]{\includegraphics[width = 0.3\textwidth]{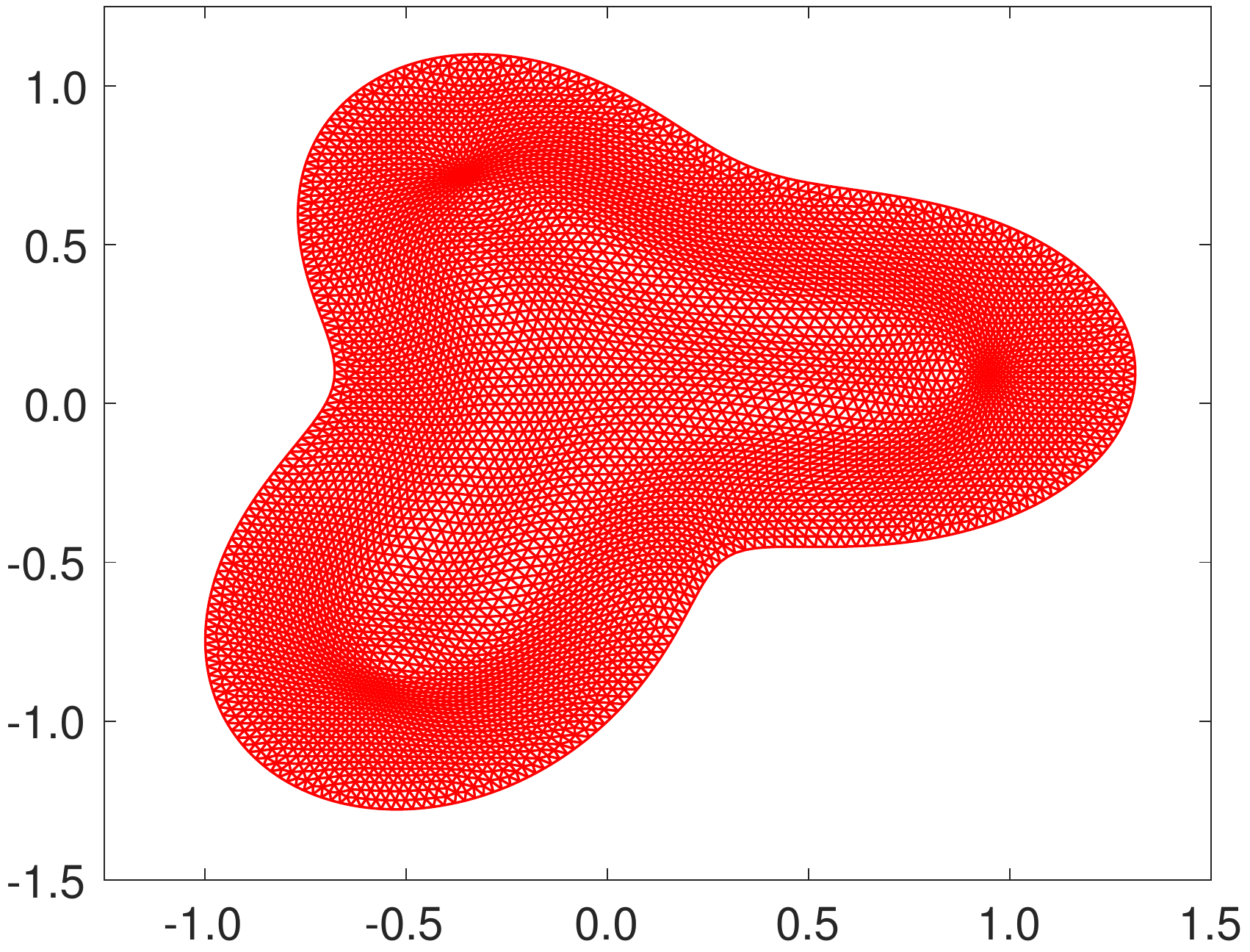}}\qquad
\subfigure[Mesh at touchdown, $\eps= 0.024$. Mesh size is $N = 11955$.]{\includegraphics[width = 0.3\textwidth]{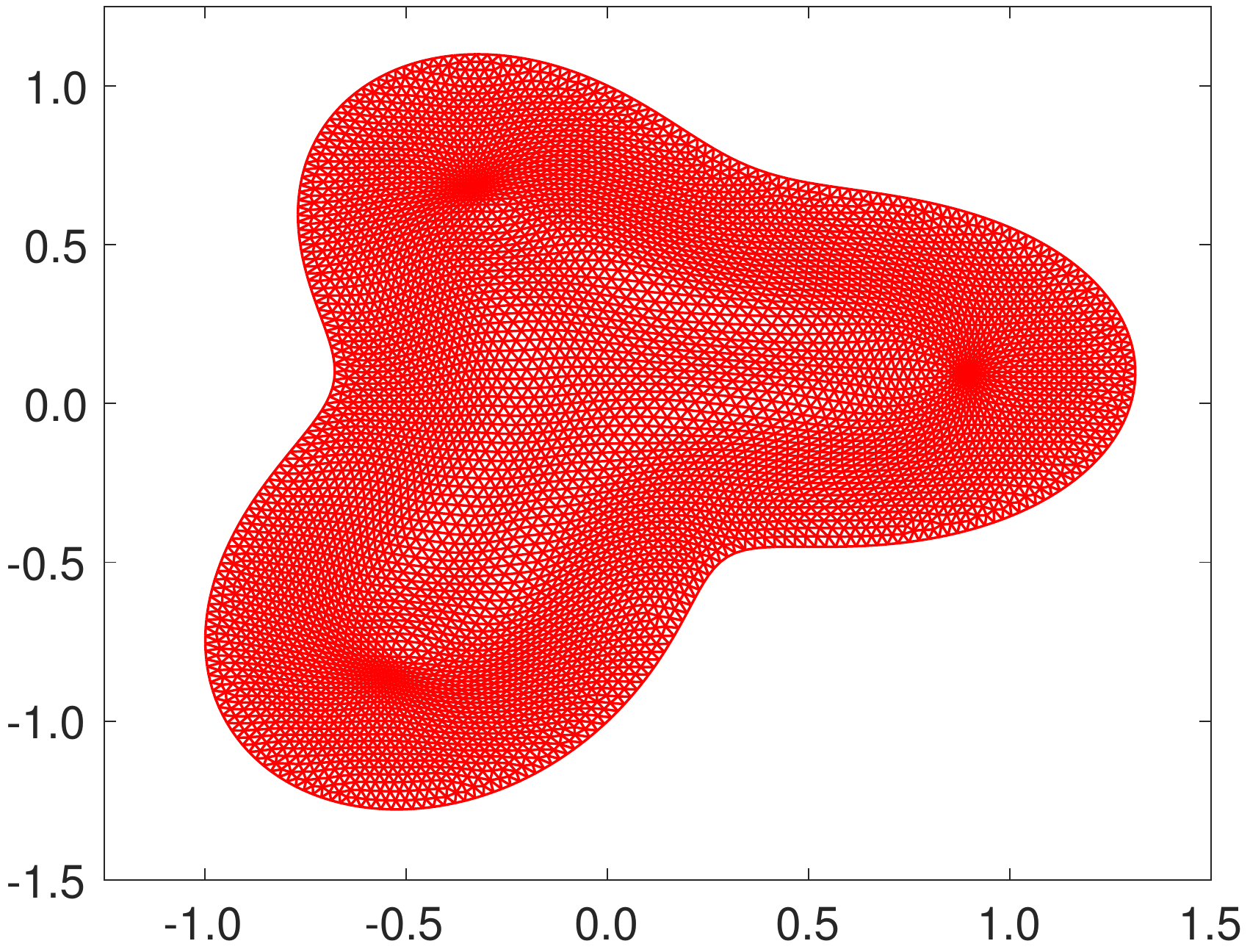}}\qquad
\subfigure[Mesh at touchdown, $\eps= 0.04$. Mesh size is $N = 11955$.]{\includegraphics[width = 0.3\textwidth]{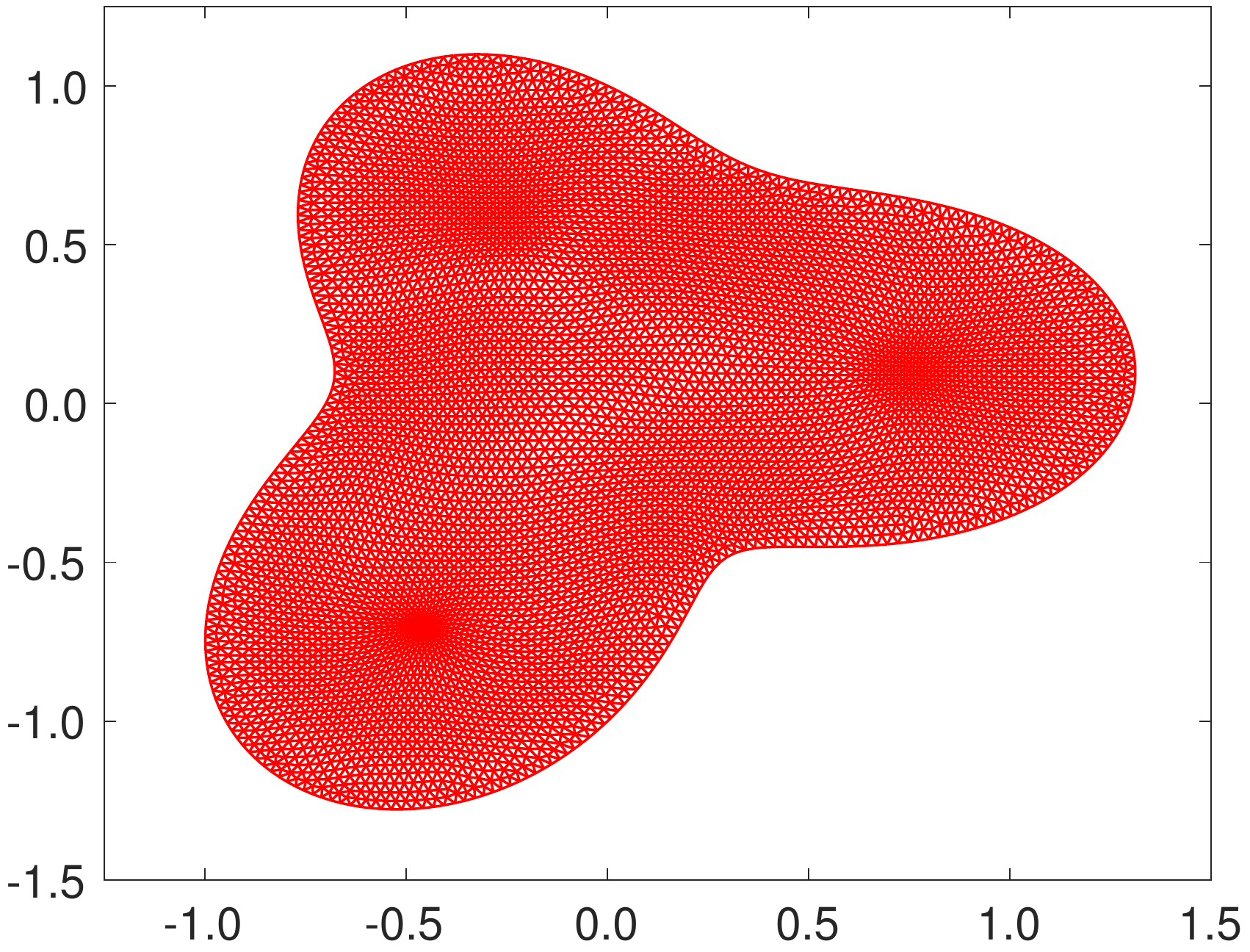}}
\caption{ The top row shows mesh pictures of $\eps = 0.02$, $\eps = 0.024$, $\eps = 0.04$ with mesh size $N = 5244 $. The associated profiles obtained by mesh size $N = 11955$ are presented in the second row.
\label{par12touchdownMesh}}
\end{figure}

In Fig.~\ref{partouchdowneps=0.02_f} the solution of \eqref{eq:intro} close to singularity is shown for $\eps=0.02$ and three distinct troughs in the solution are clear. The trough with the lowest value, and the one that contacts first, is centered at Mark 1 in Fig.~\ref{parsk}. For the slightly increased parameter value $\eps=0.024$ the solution close to touchdown is shown in Fig.~\ref{partouchdowneps=0.024_f}. At this value, the qualitative solution features look very similar, however, the center of the lowest trough is now shifted to Mark 2 on a separate branch of $\skel$. At the value $\eps = 0.04$, with final profile shown in Fig.~\ref{partouchdowneps=0.04_f}, the lowest point has shifted again and is now centered on the third arm of $\skel$ at Mark 3 on Fig.~\ref{parsk}. These observations suggest that in this asymmetric case simultaneous two point touchdown can occur for particular fixed values of $\eps$ in the ranges $(0.02,0.024)$ and $(0.024,0.4)$. The final profiles for $\eps = 0.02$, $0.024$, and $0.04$ obtained with mesh sizes $N = 5244$ and $N = 11955$ are shown in Fig.~\ref{par12touchdownProfile}. The associated mesh pictures are shown in Fig.~\ref{par12touchdownMesh}. We can see that the allocation of the singularities is robust with respect to grid refinement.

At the larger value $\eps =0.092$, the snapshots in Fig.~\ref{partouchdowneps=0.092} show that the three peaks merge very quickly in the evolution of the solution. By the time the solution of \eqref{eq:intro}  is close to singularity at this value of $\eps$, the solution has only one trough which is centered close to the geometric center of the domain at Mark 4 in Fig.~\ref{parsk}.

In each of the intermediary evolution plots in Figs.~\ref{partouchdowneps=0.02b}-\ref{partouchdowneps=0.092b}, the mesh is seen to be accurately capturing the firestorm set $\omega(t)$. In this solution regime, the adaptive algorithm allocates grid resolution to the vicinity of $\partial\Omega$ in order to capture this expanding boundary layer. The third snapshot of the solution is taken very close to touchdown and the mesh has adapted to increase the resolution in the vicinity of the forming singularities. This shows the numerical method capturing multiple types of dynamic fine scale solution modalities.

\clearpage

%\begin{figure}[htbp]
%\centering
%\subfigure[Mesh at $t= 0.009$.]{\includegraphics[width = 0.3\textwidth]{par1mesh1.pdf}}\qquad
%\subfigure[Mesh at $t= 0.2864$.]{\includegraphics[width = 0.3\textwidth]{par1mesh2.pdf}}\qquad
%\subfigure[Mesh at $t= 0.3101$.]{\includegraphics[width = 0.3\textwidth]{par1mesh3.pdf}}\\
%\subfigure[Solution at $t= 0.009$.]{\includegraphics[width = 0.3\textwidth]{par1u1.pdf}}\qquad
%\subfigure[Solution at $t= 0.2864$.]{\includegraphics[width = 0.3\textwidth]{par1u2.pdf}}\qquad
%\subfigure[Solution at $t= 0.3101$.]{\includegraphics[width = 0.3\textwidth]{par1u3.pdf}}
%\caption{ The evolution of the solution for $\eps = 0.052$. The mesh size $N = 5082$. \label{p52}}
%\end{figure}

\clearpage

\section{Conclusions}\label{sec:discussion}

This paper has studied the influence of geometry and parameter values on the location of singularities in a fourth-order PDE system modeling microscopic elastic-electrostatic deflections. We have developed a precision numerical tool for exploring the contact sets in these elastic deformations which is an important problem in the design of micro-electro mechanical systems. Specifically, we have developed an adaptive moving mesh PDE method which dynamically relocates the mesh points to provide additional resolution in spatial regions with fine scale solution behavior. This method can automatically detect and resolve different types of dynamic features such as sharp interfaces and multiple forming singularities. The method can also accommodate the complex geometries and topological defects common in the design of real MEMS devices (cf.~Fig.~\ref{fig:intro_diagram_B}).

To complement this numerical tool, we have used an asymptotic analysis to obtain the skeleton - a reduced representation of the domain which gives an estimate of the potential singularity locations. This analysis also reveals that the sensitive dependence of the contact set on the equation parameters and the shape of the domain $\Omega$ is due to a non-monotone boundary layer profile. The superposition of the solution along rays emanating from $\partial\Omega$ lowers the value of the solution at certain points in the domain which then become more likely to be singularity locations. We find that the skeleton gives a good qualitative description of the possible contact sets. The quantitative accuracy of the skeleton is variable and in particular we find it to be diminished in non-simply connected domains. For engineers seeking to prevent the two surfaces coming into physical contact through the placement of deflection limiters \cite{Krylov2010}, the skeleton provides a good estimate of the points at which these should be centered.

The range of potential applications for this new numerical method is significant. Adaptive methods for higher-order PDE systems such as \eqref{eq:intro} have received somewhat less attention in comparison to their second-order counterparts, particularly in higher dimensions (cf.~\cite{MovCol42007} for one-dimension). The methods developed here are applicable to a wide variety of other higher-order PDE systems that are central features in many topical applications such as rock folding \cite{Dodwell2012a}, ion bombardment lithography \cite{Perkinson2016}, thin films dynamics \cite{LBGK2007,Witelski2017} and pattern formation \cite{LG2013,Gavish2011,Dai2014}.

\vspace{10pt}

\section*{Acknowledgments}
A.~E.~Lindsay was supported by NSF grant DMS-1516753.  K.~L.~DiPietro was supported by NSF grant DGE-1313583. R.~D.~Haynes was support by a NSERC Discovery Grant (Canada).


\vspace{10pt}

\begin{thebibliography}{10}

%\bibliography{MMbib}{}
%\bibliographystyle{plain}
%\bibliographystyle{abbrv}
%

\bibitem{Bai94a}
M.~J. Baines.
\newblock {\em Moving Finite Elements}.
\newblock Oxford University Press, Oxford, 1994.

\bibitem{Baines-2011}
M.~J. Baines, M.~E. Hubbard, and P.~K. Jimack.
\newblock Velocity-based moving mesh methods for nonlinear partial differential
  equations.
\newblock {\em Comm. Comput. Phys.}, 10:509--576, 2011.

\bibitem{BatraReview}
R.~C. Batra, M.~Porfiri, and D.~Spinello.
\newblock Review of modeling electrostatically actuated microelectromechanical systems.
\newblock {\em Smart Mater. Struct.}, 16:R23--R31, 2007.

\bibitem{BGW}
C.~J.~Budd, V.~Galaktionov, and J.~F. Williams.
\newblock Self-similar blow-up in higher-order semilinear parabolic equations.
\newblock {\em SIAM J. Appl. Math.}, 64:1775--1809, 2004.

\bibitem{BHR09}
C.~J. Budd, W.~Huang, and R.~D. Russell.
\newblock Adaptivity with moving grids.
\newblock {\em Acta Numer.}, 18:111--241, 2009.

\bibitem{Budd2006}
C.~J.~Budd and J.~F.~Williams.
\newblock Parabolic {M}onge-{A}mp\'{e}re methods for blow-up problems in several spatial dimensions,
\newblock {\em J. Phys. A}, 39:5425--5444, 2006.

\bibitem{PMA2009}
C. J. Budd and J. F. Williams.
\newblock Moving Mesh Generation Using the Parabolic {M}onge-�{�A}mp\'{e}re Equation,
\newblock {\em SIAM J. Sci. Comput.}, 31:3438--3465, 2009.

\bibitem{Budd2010}
C. J. Budd and J. F. Williams.
\newblock How to adaptively resolve evolutionary singularities in differential equations with symmetry
\newblock {\em J. Eng. Math.}, 66:217--236, 2010.

\bibitem{Ceniceros2001}
H, D. Ceniceros and Thomas Y. Hou.
\newblock An efficient dynamically adaptive mesh for potentially singular solutions,
\newblock {\em J. Comput. Phys.}, 172:609--639, 2001.

\bibitem{LD2017} K.~L.~DiPietro and A.~E.~Lindsay.
Monge-Amp\'{e}re simulation of fourth order PDEs in two dimensions with
application to elastic-electrostatic contact problems.
{\em J. Comput. Phys.}, 349:328--350, 2017.

\bibitem{Dai2014}
S. Dai and K. Promislow.
\newblock Geometric evolution of bilayers under the functionalized Cahn{\textendash}Hilliard equation.
\newblock {\em Proc. R. Soc. Lond. Ser. A Math. Phys. Eng. Sci.}, 469:20120505 (20 pp.),  2013.

\bibitem{Dodwell2012a} T. J. Dodwell, M. A. Peletier, C. J. Budd, and G. W. Hunt.
\newblock Self-similar voiding solutions of a single layered model of folding rocks.
{\em SIAM J. Appl. Math.}, 72:444--463, 2012.

\bibitem{Friedman88}
A.~Friedman and L.~Oswald.
\newblock The blow-up time for higher order semilinear parabolic equations with
  small leading coefficients.
\newblock {\em J. Diff. Eq.}, 75:239--263, 1988.

\bibitem{GALAK}
V.~Galaktionov.
\newblock Five types of blow-up in a semilinear fourth-order reaction-diffusion
  equation: an analytical-numerical approach.
\newblock {\em Nonlinearity}, 22:1695--1741, 2009.

\bibitem{Galaktionov2002}
V.~A. Galaktionov and J.-L. V{\'a}zquez.
\newblock The problem of blow-up in nonlinear parabolic equations.
\newblock {\em Dis. Cont. Dyn. Sys.}, 8:399--433, 2002.

\bibitem{Gavish2011}
N. Gavish, G. Hayrapetyan, K. Promislow, and L. Yang.
\newblock Curvature driven flow of bi-layer interfaces.
\newblock {\em Phys. D}, 240:675--693, 2011.

\bibitem{LG2013}
K. B. Glasner and A. E. Lindsay.
\newblock The stability and evolution of curved domains arising from one-dimensional localized patterns.
\newblock {\em SIAM J. Appl. Dyn. Sys.}, 12:650--673, 2013.

\bibitem{Montijano2004}
S.~Gonz{\'a}lez-Pinto, J.~I. Montijano, and S.~P{\'e}rez-Rodr{\'{\i}}guez.
\newblock Two-step error estimators for implicit {R}unge-{K}utta methods
  applied to stiff systems.
\newblock {\em ACM Trans. Math. Software}, 30:1--18, 2004.

\bibitem{HW96}
E.~Hairer and G.~Wanner.
\newblock {\em Solving Ordinary Differential Equations. {II}}, Volume~14 of
  {\em Springer Series in Computational Mathematics}.
\newblock Springer-Verlag, Berlin, second edition, 1996.
\newblock Stiff and differential-algebraic problems.

\bibitem{Hua01b}
W.~Huang.
\newblock Variational mesh adaptation: isotropy and equidistribution.
\newblock {\em J. Comput. Phys.}, 174:903--924, 2001.

\bibitem{Hua05b}
W.~Huang.
\newblock Metric tensors for anisotropic mesh generation.
\newblock {\em J. Comput. Phys.}, 204:633--665, 2005.

\bibitem{HK2014}
W.~Huang and L.~Kamenski.
\newblock A geometric discretization and a simple implementation for
  variational mesh generation and adaptation.
\newblock {\em J. Comput. Phys.}, 301:322--337, 2015.

\bibitem{HK2015}
W.~Huang and L.~Kamenski.
\newblock On the mesh nonsingularity of the moving mesh {PDE} method.
\newblock {\em Math. Comp.}, in press.
\newblock (DOI: 10.1090/mcom/3271)

\bibitem{HR11}
W.~Huang and R.~D. Russell.
\newblock {\em Adaptive Moving Mesh Methods}.
\newblock Springer, New York, 2011.
\newblock Applied Mathematical Sciences Series, Vol. 174.

\bibitem{Krylov2010}
S.~Krylov and N.~Dick.
\newblock Dynamic stability of electrostatically actuated initially curved
  shallow micro beams.
\newblock {\em Continuum Mech. Therm.}, 22:445--468, 2010.

\bibitem{Lulinsky2011} S. Krylov, B. R. Ilic, and S. Lulinsky.
Bistability of curved microbeams actuated by fringing electrostatic fields.
{\em Nonlinear Dyn.},  66:403--426, 2011.

\bibitem{Craighead2008} S. Krylov, B. R. Ilic, D. Schreiber, S. Seretensky, and H. Craighead.
The pull-in behavior of electrostatically actuated bistable microstructures.
{\em J. Micromech. Microeng.}, 18:055026, 2008.

\bibitem{Lin2007}
F.~Lin and Y.~Yang.
\newblock Nonlinear non-local elliptic equation modelling electrostatic actuation.
\newblock {\em Proc. R. Soc. London Ser. A}, 463:1323--1337, 2007.

\bibitem{Capacitor2012}
A.~Lindsay and J.~Lega.
\newblock Multiple quenching solutions of a fourth order parabolic pde with a
  singular nonlinearity modeling a {MEMS} capacitor.
\newblock {\em SIAM J. Appl. Math.}, 72:935--958, 2012.

\bibitem{Equilibrium2014}
A.~Lindsay, J.~Lega, and K.~Glasner.
\newblock Regularized model of post-touchdown configurations in electrostatic
  {MEMS}: Equilibrium analysis.
\newblock {\em Physica D}, 280/281:95--108, 2014.

\bibitem{selfsimilar2014}
A.~E. Lindsay.
\newblock An asymptotic study of blow up multiplicity in fourth order parabolic
  partial differential equations.
\newblock {\em Dis. Cont. Dyn. Sys. Ser. B}, 19:189--215, 2014.

\bibitem{Lindsay2016}
A.~E. Lindsay.
\newblock Regularized model of post-touchdown configurations in electrostatic
  {MEMS}: bistability analysis.
\newblock {\em J. Eng. Math.}, 99:65--77, 2016.

\bibitem{Lindsay2013}
A.~E. Lindsay, J.~Lega, and F.~J. Sayas.
\newblock The quenching set of a {MEMS} capacitor in two-dimensional geometries.
\newblock {\em J. Nonlinear Sci.}, 23:807--834, 2013.

\bibitem{Malandain1998}
G.~Malandain and S.~Fern{\'a}ndez-Vidal.
\newblock Euclidean skeletons.
\newblock {\em Image Vis. Comput.}, 16:317--327, 1998.

\bibitem{Pelesko2002}
J.~A. Pelesko.
\newblock Mathematical modeling of electrostatic {MEMS} with tailored
  dielectric properties.
\newblock {\em SIAM J. Appl. Math.}, 62:888--908, 2002.

\bibitem{Perkinson2016} J. C. Perkinson, M. J. Aziz, M. P. Brenner, and M. Holmes-Cerfon.
\newblock Designing steep, sharp patterns on uniformly ion-bombarded surfaces.
\newblock {\em Proc. Nat. Acad. Sci.}, 113:11425--11430, 2016.

\bibitem{PB} J.~A.~Pelesko and D.~H.~Bernstein.
{\em Modeling MEMS and NEMS},
Chapman Hall and CRC Press, 2002.

\bibitem{Philippin15} G. A. Philippin.
Blow-up phenomena for a class of fourth-order parabolic problems.
{\em Proc. Amer. Math. Soc.}, 143:2507--2513, 2015.

\bibitem{Qiu2004} J. Qiu, J. H. Lang, and A. H. Slocum.
A curved-beam bistable mechanism.
{\em J. Microelectromech. Syst.}, 13:137--146, 2004.

\bibitem{MovCol42007}
 R. D. Russell, J. F. Williams, and X. Xu.
\newblock MOVCOL4: A Moving Mesh Code for Fourth-Order Time-Dependent Partial Differential Equations
\newblock {\em SIAM J. Sci. Comput.}, 29:197--220, 2007.

\bibitem{sandia} Courtesy of Sandia National Laboratories, SUMMiT(TM) Technologies,
{\tt www.mems.sandia.gov}.

\bibitem{Tan05}
T.~Tang.
\newblock Moving mesh methods for computational fluid dynamics flow and transport.
\newblock In {\em Recent Advances in Adaptive Computation (Hangzhou, 2004)},
  Volume 383 of {\em AMS Contemporary Mathematics}, pages 141--173. Amer. Math.
  Soc., Providence, RI, 2005.

\bibitem{Tsai07} N.-C. Tsai and C.-Y. Sue.
Review of MEMS-based drug delivery and dosing systems.
{\em Sens. Actuators A Phys.}, 134:555--564, 2007.

\bibitem{Beek2012}
J.~T.~M. van Beek and R.~Puers.
\newblock A review of {MEMS} oscillators for frequency reference and timing
  applications.
\newblock {\em J. Micromech. Microeng.}, 22:013001, 2012.

\bibitem{Witelski2017} H. Ji and T. P. Witelski.
\newblock Finite-time thin film rupture driven by modified evaporative loss.
\newblock {\em Phys. D}, 342:1--15, 2017

\bibitem{Witelski00}
T.~P. Witelski and A.~J. Bernoff.
\newblock Dynamics of three-dimensional thin film rupture.
\newblock {\em Phys. D}, 147:155--176, 2000.

\bibitem{LBGK2007}
H.-W. Lu,  K. Glasner, A. L. Bertozzi, and C.-J. Kim.
\newblock A diffuse-interface model for electrowetting drops in a {H}ele-{S}haw cell.
\newblock {\em J. Fluid Mech.}, 590:411--��435 2007.


\bibitem{ZHANG2014}
R.~Zhang and L.~Cai.
\newblock On the semi linear equations of electrostatic {NEMS} devices.
\newblock {\em Z. Angew. Math. Phys.}, 65:1207--1222, 2014.

\end{thebibliography}
\end{document}